\documentclass[notitlepage,leqno,11pt]{article}

\usepackage{latexsym}
\usepackage{amsmath}
\usepackage{amsfonts}
\usepackage{amssymb}
\usepackage{hyperref}
\renewcommand{\a }{\alpha }
\renewcommand{\b }{\beta }
\renewcommand{\d}{\delta }
\newcommand{\D }{\Delta }

\newcommand{\e }{\varepsilon }
\newcommand{\g }{\gamma}

\newcommand{\G }{\Gamma }
\renewcommand{\l }{\lambda }
\renewcommand{\L }{\Lambda }

\newcommand{\n }{\nabla }
\newcommand{\var }{\varphi }

\newcommand{\s }{\sigma }
\newcommand{\Sig }{\Sigma}
\renewcommand{\t }{\tau }
\renewcommand{\th }{\theta }
\renewcommand{\o }{\omega }
\renewcommand{\O }{\Omega }
\newcommand{\z }{\zeta }
\newcommand{\zn }{{\zeta_{n+1}}}
\newcommand{\oy }{\overline{y}}
\newcommand{\ov}{\overline}
\newcommand{\un}{\underline}
\newcommand{\U}{\Upsilon}
\newcommand{\pa}{\partial}
\newcommand{\Pti}{\Phi_{\tilde{I}-2}}

\newcommand{\be}{\begin{equation}}
\newcommand{\ee}{\end{equation}}
\newenvironment{pf}{\noindent{\sc Proof}.\enspace}{\rule{2mm}{2mm}\medskip}
\newenvironment{pfn}{\noindent{\sc Proof}}{\rule{2mm}{2mm}\medskip}

\newtheorem{lemma}{Lemma}[section]

\newcommand{\R}{\mathbb{R}}

\newcommand{\N}{\mathbb{N}}

\newcommand{\oa}{\ov{a}}
\newcommand{\ob}{\ov{b}}
\newcommand{\no}{\noindent}
\newcommand{\ms}{\medskip}

\author{Fethi MAHMOUDI and Andrea MALCHIODI $^1$ }

\date{}

\title{Concentration on  minimal submanifolds for a singularly perturbed  Neumann problem}

\begin{document}

\newtheorem{lem}{Lemma}[section]
\newtheorem{pro}[lem]{Proposition}
\newtheorem{thm}[lem]{Theorem}
\newtheorem{rem}[lem]{Remark}
\newtheorem{cor}[lem]{Corollary}
\newtheorem{df}[lem]{Definition}
\newtheorem{rems}[lem]{Remarks}

\maketitle

\begin{center}

{\small Sissa, Via Beirut 2-4, 34014 Trieste, Italy}

\end{center}

\footnotetext[1]{E-mail addresses: mahmoudi@sissa.it,
malchiod@sissa.it}

\

\

\

\noindent {\sc abstract}. We consider the equation $- \e^2 \D u + u
= u^p$ in $\O \subseteq \R^N$, where $\O$ is open, smooth and
bounded, and we prove concentration of solutions along
$k$-dimensional minimal submanifolds of $\pa \O$, for $N \geq 3$ and
for $k \in \{1, \dots, N-2\}$. We impose Neumann boundary
conditions, assuming $1<p <\frac{N-k+2}{N-k-2}$ and $\e \to 0^+$.
This result settles in full generality a phenomenon previously
considered only in the particular case $N = 3$ and $k = 1$.
\begin{center}

\bigskip\bigskip

\noindent{\it Key Words:} Singularly Perturbed Elliptic Problems,
Differential geometry,  Local Inversion, Fourier Analysis.

\bigskip

\centerline{\bf AMS subject classification: 35B25, 35B34, 35J20,
35J60, 53A07}

\end{center}

\section{Introduction}\label{s:intr}

In this paper we study concentration phenomena for the problem
\begin{equation}\tag{$P_\e$}\label{eq:pe}
  \begin{cases}
    -\e^2 \D u + u = u^p & \text{ in } \O, \\
    \frac{\partial u}{\partial \nu} = 0 & \text{ on } \partial \O,
    \\ u > 0 & \text{ in } \O,
  \end{cases}
\end{equation}
where $\O$ is a smooth bounded domain of $\R^N$, $p > 1$, and where
$\nu$ denotes the unit normal to $\partial \O$. Given a smooth
embedded non-degenerate minimal submanifold $K$ of $\pa \O$, of
dimension $k \in \{ 1, \dots, N-2 \}$, we prove existence of
solutions of \eqref{eq:pe} concentrating along $K$. Since the
solutions we find have a specific asymptotic profile, which is
described below, a natural restriction on $p$ is imposed, depending
on the dimension $N$ and $k$, namely $p < \frac{N-k+2}{N-k-2}$.

\

\

Problem \eqref{eq:pe} or some of its variants (including the
presence of non-homogeneous terms, different boundary conditions,
etc.) arise in several contexts, as the Nonlinear Schr\"odinger
Equation or from modeling reaction-diffusion systems, see for
example \cite{am}, \cite{fw}, \cite{ni} and references therein. A
typical phenomenon one observes is the existence of solutions which
are sharply concentrated near some subsets of their domain.

Concerning reaction-diffusion systems, this phenomenon is related to
the so-called Turing's instability, \cite{tu}. According to this
principle, reaction-diffusion systems whose reactants have very
different diffusivities might generate stable non-trivial patterns.
This is indeed more likely to happen when more reactants are present
since, as shown in \cite{ch}, \cite{mat}, scalar reaction-diffusion
equations in a convex domain admit only constant stable equilibria.

A well-know system is the following one
\begin{equation}\label{eq:gm}\tag{$GM$}
   \begin{cases}
     \mathcal{U}_t = d_1 \D \mathcal{U} - \mathcal{U} +
     \frac{\mathcal{U}^p}{\mathcal{V}^q} & \text{ in } \O \times (0,+\infty),
     \\ \mathcal{V}_t = d_2 \D \mathcal{V} - \mathcal{V} +
     \frac{\mathcal{U}^r}{\mathcal{V}^s} & \text{ in } \O \times (0,+\infty),
   \\ \frac{\partial \mathcal{U}}{\partial \nu} =
   \frac{\partial \mathcal{V}}{\partial \nu} = 0 & \text{ on }
   \partial \O \times (0,+\infty),
   \end{cases}
\end{equation}
introduced in \cite{gm} to describe some biological experiment. The
functions $\mathcal{U}$ and $\mathcal{V}$ represent the densities of
some chemical substances, the numbers $p, q, r, s$ are non-negative
and such that $0 < \frac{p-1}{q} < \frac{r}{s+1}$, and it is assumed
that the diffusivities $d_1$ and $d_2$ satisfy $d_1 \ll 1 \ll d_2$.
In the stationary case of \eqref{eq:gm}, as explained in \cite{ni},
\cite{nty}, when $d_2 \to + \infty$ the function $\mathcal{V}$ is
close to a constant (being nearly harmonic and with zero normal
derivative at the boundary), and therefore the equation satisfied by
$\mathcal{U}$ is similar to \eqref{eq:pe}, with $\e^2 = d_1$.

\

The typical concentration behavior of solutions $u_\e$ to
\eqref{eq:pe} is via a scaling of the variables in the form $u_\e(x)
\sim u_0 \left( \frac{x-Q}{\e} \right)$, where $Q$ is some point of
$\ov\O$, and where $u_0$ is a solution of the problem
\begin{equation}\label{eq:u0}
    - \D u_0 + u_0 = u_0^p \qquad \qquad \text{ in } \R^N \quad
    \hbox{ (or in } \R^N_+ = \{ (x_1, \dots, x_N) \in \R^N \; : \;
    x_N > 0\}),
\end{equation}
the domain depending on whether $Q$ lies in the interior of $\O$ or
at the boundary; in the latter case Neumann conditions are imposed.

When $p < \frac{N+2}{N-2}$ (and indeed only if this inequality is
satisfied), problem \eqref{eq:u0} admits positive radial solutions
which decay to zero at infinity. Solutions of \eqref{eq:pe} with
this profile are called {\em spike-layers}, since they are highly
concentrated near some point of $\ov{\O}$. There is an extensive
literature regarding this type of solutions, beginning from the
papers \cite{lnt}, \cite{nt91}, \cite{nt93}. Indeed their structure
is very rich, and there are also solutions with multiple peaks, both
at the boundary and at the interior of $\O$. We refer for example to
the papers \cite{dy}, \cite{dfw}, \cite{gr2}, \cite{gui},
\cite{guiw2}, \cite{guiww}, \cite{liyy}, \cite{ln}, \cite{w97}.

\

In recent years, some new types of solutions have been constructed:
they indeed concentrate at sets of positive dimension and their
profile consists of solutions of \eqref{eq:u0} which do not decay to
zero at infinity. In \cite{malm}, \cite{malm2} it has been shown
that given any smooth bounded domain $\O \subseteq \R^N$, $N \geq
2$, and any $p > 1$, there exists a sequence $\e_j \to 0$ such that
$(P_{\e_j})$ possesses solutions concentrating at $\partial \O$
along this sequence. Their profile is a solution of \eqref{eq:u0}
(for $N = 1$) on the half real line which tends to zero at infinity
and which satisfies the condition $u'_0(0) = 0$. This function can
also be trivially extended as a cylindrical solution to
\eqref{eq:u0} on the whole $\R^N_+$ .

Later in \cite{mal2} it has been proved that, if $\O$ is a smooth
bounded set of $\R^3$, if $p > 1$ and if $h$ is a closed, simple
non-degenerate geodesic on $\pa \O$, then there exists again a
sequence $(\e_j)_j$ converging to zero such that $(P_{\e_j})$ admits
solutions $u_{\e_j}$ concentrating along $h$ as $j$ tends to
infinity. In this case the profile of $u_{\e_j}$ is a decaying
solution of \eqref{eq:u0} in $\R^2_+$, again extended to a
cylindrical solution in higher dimension.

These are examples of a phenomenon which has been conjectured to
hold in more general cases: in fact it is expected that, under
generic assumptions, if $\O \subseteq \R^N$ and if $k$ is an integer
between $1$ and $N - 1$, there exist solutions of \eqref{eq:pe}
concentrating along $k$-dimensional sets when $\e$ tends to zero.
While the case $k = N - 1$ has been tackled in \cite{malm2}, the
goal of the present paper is to consider $k \leq N - 2$, and to
 prove this conjecture under rather mild assumptions on the limit
set. Before stating our main theorem we introduce some preliminary
notation.

Given a smooth $k$-dimensional manifold $K$ of $\pa \O$, and given
any $q \in K$ we can choose a system of coordinates $(\oy,\z)$ in
$\O$ orthonormal at $q$ and such that $(\oy, 0)$ are coordinates on
$K$, and with the property that
\begin{equation}\label{eq:onsysq}
    \frac{\pa}{\pa \oy_a}|_q \in T_q K, \quad a = 1, \dots, k; \qquad
    \quad \frac{\pa}{\pa \z_i}|_q \in T_q \pa \O, \quad i = 1, \dots, n;
    \qquad \quad \frac{\pa}{\pa \zn}|_q = \nu(q),
\end{equation}
where we have set $n = N - k - 1$. Our main theorem is the
following: we refer to Section \ref{s:geom} for the geometric
terminology.

\begin{thm}\label{t:m}
Let $\O \subseteq \R^N$, $N \geq 3$, be a smooth and bounded domain,
and let $K \subseteq \pa \O$ be a compact embedded non-degenerate
minimal submanifold of dimension $k \in \{ 1, \dots, N-2 \}$. Then,
if $p \in \left( 1, \frac{N-k+2}{N-k-2} \right)$, there exists a
sequence $\e_j \to 0$ such that $(P_{\e_j})$ admits positive
solutions $u_{\e_j}$ concentrating along $K$ as $j \to \infty$.
Precisely there exists a positive constant $C$, depending on $\O, K$
and $p$ such that for any $x \in \O$ $u_{\e_j}(x) \leq C e^{-\frac{
dist(x,K)}{C \e_j}}$; moreover for any $q \in K$, in a system of
coordinates $(\oy,\z)$ satisfying \eqref{eq:onsysq}, for any integer
$m$ one has $u_{\e_j}(0,\e_j \cdot)
\stackrel{C^m_{loc}(\R^{n+1}_+)}{\longrightarrow} w_0(\cdot)$, where
$w_0 : \R^{n+1}_+ \to \R$ is the unique radial solution of
\begin{equation}\label{eq:w0}
  \begin{cases}
    - \D u + u = u^p & \hbox{ in } \R^{n+1}_+, \\
    \frac{\partial u}{\partial \nu} = 0 & \hbox{ on } \partial
    \R^{n+1}_+, \\
    u > 0, u \in H^1(\R^{n+1}_+).
  \end{cases}
\end{equation}
\end{thm}

\

\begin{rems} (a) Differently from the previous papers concerning
the case $N = 3$ and $k = 1$, or concentration at the whole $\pa
\O$, we require an upper bound on $p$ depending on $N$ and $k$. This
condition is rather natural, since \eqref{eq:w0} is solvable if and
only if $p < \frac{N-k+2}{N-k-2}$, see \cite{bl}, \cite{po},
\cite{str} and in this case the solution is radial and unique (up to
a translation), see \cite{gnn}, \cite{kwo}. In any case, our
assumptions allow supercritical exponents as well.

(b) As for the results in \cite{mal2}, \cite{malm} and \cite{malm2},
existence is proved only along a sequence $\e_j \to 0$ (actually
with our proof it can be obtained for $\e$ in a sequence of
intervals $(a_j, b_j)$ approaching zero, but not for any small
$\e$). This is caused by a resonance phenomenon we are going to
discuss below, explaining the ideas of the proof. This resonance is
peculiar of multidimensional spike-layers, see also \cite{dkw}, and
other geometric problems, see \cite{mmp}, \cite{mp}. In some cases,
when some symmetry is present, it is possible to get rid of this
resonance phenomenon working in spaces of invariant functions. We
refer for example to the papers \cite{amn, amn2, bd, bape, dy2, dap,
mnw, mopa}.
\end{rems}

\

\noindent We can describe the resonance phenomenon, which causes the
main difficulty in proving Theorem \ref{t:m}, in the following way.
By the change of variables $x \mapsto \e x$, we are reduced to
consider the problem
\begin{equation}\tag{$\tilde{P}_\e$}\label{eq:tpe}
  \begin{cases}
    - \D u + u = u^p & \text{ in } \O_\e, \\
    \frac{\partial u}{\partial \nu} = 0 & \text{ on } \partial \O_\e,
    \\ u > 0 & \text{ in } \O_\e,
  \end{cases}
\end{equation}
where $\O_\e = \frac 1 \e \O$. As for \eqref{eq:onsysq}, given
$\hat{q} \in K_\e := \frac 1 \e K$, we can choose scaled coordinates
$(y,\z)$ on $\O_\e$ such that $\pa_{y_a}|_{\hat q} \in T_{\hat q}
K_\e$, $\pa_{\z_i}|_{\hat q} \in T_{\hat q} \pa \O_\e$ and
$\pa_{\z_{n+1}}|_{\hat q} = \nu(\hat q)$. Then, letting
$\tilde{u}_\e$ denote the scaling of $u_\e$ to $\O_\e$,  we have
that, in a plane through $\hat{q}$ normal to $K_\e$, $\tilde{u}_\e$
behaves like $\tilde{u}_\e (0,\z) = u_\e(0, \e \z) \simeq w_0(\z)$.
This amounts to the fact that $\tilde{u}_\e(x) \simeq
w_0\left(dist(x,K_\e)\right)$, $x \in \O_\e$, and therefore
$\tilde{u}_\e$ has a fixed profile in the directions perpendicular
to the expanding domain $K_\e$. Since the function
$w_0\left(dist(x,K_\e)\right)$ can be considered as an {\em
approximate solution} to \eqref{eq:tpe}, it is natural to use local
inversion arguments near this function in order to find true
solutions. For this purpose it is necessary to understand the
spectrum of the linearization of \eqref{eq:tpe} at approximate
solutions.

For simplicity, let us assume for the moment that $K$ is
$(N-2)$-dimensional, namely that its codimension in $\pa \O$ is
equal to $1$, as in \cite{mal2}. Then, letting $\tilde{\nu}$ denote
the normal to $K$ in $\pa \O$, we can parameterize naturally a
neighborhood of $K_\e$ as a product of the form $K_\e \times \left(-
\frac \d \e, \frac \d \e\right)$, where $\d$ is a small positive
number, via the exponential map in $\pa \O_\e$
\begin{equation}\label{eq:pardo}
    (y,s) \mapsto \exp^{\pa \O_\e}_y (s \tilde\nu); \qquad \qquad (y,s)
  \in K_\e \times \left(- \frac \d \e, \frac \d \e\right).
\end{equation}
Similarly, if $\nu(y,s)$ is the inner unit normal to $\pa \O_\e$ at
the image of $(y,s)$ under the above map, we can parameterize a
neighborhood of $K_\e$ in $\O_\e$ with a product $K_\e \times
\left(- \frac \d \e, \frac \d \e\right) \times \left( 0, \frac \d
\e\right)$ by
$$
  (y,s,t) \mapsto \exp^{\pa \O_\e}_y (s \tilde\nu) + t \nu(y,s);
  \qquad \qquad (y,s,t) \in K_\e \times \left(- \frac \d \e, \frac
  \d \e\right) \times \left( 0, \frac \d \e\right).
$$
When $\e$ tends to zero, the standard Euclidean metric of $\O_\e$
becomes closer and closer (on the above set) to the product of the
metric of $K_\e$ and that of $\R^2$ (parameterized by the variables
$s$ and $t$ as cartesian coordinates). Therefore, since the set
$\left(- \frac \d \e, \frac \d \e\right) \times \left( 0, \frac \d
\e\right)$ converges to $\R^2_+ = \left\{ (s,t) \in \R^2 \; : \; t >
0 \right\}$, in a first approximation we get that the linearization
of \eqref{eq:tpe} at $\tilde{u}_\e$ is
\begin{equation}\label{eq:modlin}
    \begin{cases}
    - \D_{K_\e} u - \pa^2_{ss} u - \pa^2_{tt} u + u - p w_0(\z) u =
    0 & \hbox{ in } K_\e \times \R^2_+, \\
    \frac{\partial u}{\partial \nu} = 0 & \hbox{ on } K_\e \times
    \partial \R^{2}_+.
  \end{cases}
\end{equation}
The spectrum of this linear operator can be evaluated almost
explicitly. Referring to Section \ref{s:mod} for details (see also
\cite{mal2}, Proposition 2.9 for  the case $N = 3$), here we just
give some qualitative description of its properties.

Given an arbitrary function $u \in H^1(K_\e \times \R^2_+)$, we can
decompose it in Fourier modes in the variables $\z = (s,t)$ as
$$
  u(y,\z) = \sum_j \phi_j(\e y) u_j(\z).
$$
Here $\phi_j$ are the eigenfunctions of the Laplace-Beltrami
operator on $K$, namely $- \D_K \phi_j = \rho_j \phi_j$, $j = 0, 1,
2, \dots$, where the eigenvalues $(\rho_j)_j$ are counted with their
multiplicities.

If $u$ is an eigenfunction (with respect to the duality induced by
the space $H^1(K_\e \times \R^2_+)$) of the linear operator in
\eqref{eq:modlin} with corresponding eigenvalue $\l$, then it can be
shown (see Section \ref{s:mod} for details) that the functions $u_j$
satisfy the equation
\begin{equation}\label{eq:lina}
    \left\{
    \begin{array}{ll}
      (1-\l) \left[ - \D u_j + (1 + \a) u_j \right]
      - p w_0^{p-1} u_j = 0 & \hbox{ in } \R^2_+, \\
      \frac{\pa u_j}{\pa t} = 0 & \hbox{ on } \pa  \R^2_+,
    \end{array}
  \right.
\end{equation}
where $\a = \e^2 \rho_j$. It is known that when $\a = 0$ the latter
problem admits a negative eigenvalue $\eta_0$ (with eigenfunction
$w_0$), a zero eigenvalue $\s_0$ (with eigenfunction $\pa_s w_0$),
while all the other eigenvalues are positive. This structure is due
to the fact that $w_0$ is a mountain-pass solution of \eqref{eq:w0}
(so its Morse index is at most $1$), and the presence of a kernel
derives from the fact that this equation is invariant by translation
in the $s$ variable. When $\a$ is positive instead, it turns out
that the first eigenvalue $\eta_\a$ of \eqref{eq:lina} and the
second one $\s_\a$ are strictly increasing functions of $\a$ with
positive derivative, and tend to $1$ as $\a \to + \infty$; moreover,
the eigenfunctions corresponding to $\eta_\a$ (resp. $\s_\a$) are
radial (resp. odd in $s$) for every value of $\a$. In particular,
there exists $\ov{\a} > 0$ such that $\eta_{\ov{\a}} = 0$, so when
$\e^2 \rho_j$ is close to $\ov{\a}$ we obtain some small eigenvalues
of the original linearized problem \eqref{eq:modlin}.

{From} the monotonicity in $\a$ and from the Weyl's asymptotic
formula for $\rho_j$, it follows that the eigenvalues of the
operator in \eqref{eq:modlin} are, roughly, either of the form
$\eta_0 + \e^2 j^{\frac{2}{N-2}}$ for some $j \in \N$, or of the
form $\e^2 l^{\frac{2}{N-2}}$ for some $l \in \N$, or have a uniform
positive bound from below.

In the case of general codimension it is not possible to decompose a
neighborhood of $K$ (in $\pa \O$) as for \eqref{eq:pardo}, but
instead one has to model it on the {\em normal bundle} of $K_\e$ in
$\O_\e$, see Subsection \ref{ss:model} for details. Considering the
corresponding approximate linearized operator, one can prove that
its eigenvalues are now, roughly either of the form $\eta_{\e^2
\rho_j} \simeq \eta_0 + \e^2 j^{\frac{2}{k}}$, or of the form
$\s_{\e^2 \o_l} \simeq \e^2 l^{\frac{2}{k}}$, $j, l \in \N$, or,
again, have a uniform positive bound from below. Here $(\rho_j)_j$
still represent the eigenvalues of the Laplace-Beltrami operator on
$K$, while the numbers $(\o_l)_l$ stand for the eigenvalues of the
{\em normal Laplacian} of $K$ (considered as a submanifold of $\pa
\O$), see Section \ref{s:geom} for its definition and the
corresponding Weyl's asymptotic formula. We are interested in
particular in the following two features of the spectrum:

{\bf 1) resonances:} there are two kinds of eigenvalues which can
approach zero. First of all, those of the form $\eta_{\a}$ when $\a$
is close to $\ov{\a}$. This happens when $\e^2 j^{\frac{2}{k}}
\simeq \ov{\a}$, namely when $j \simeq \e^{-k}$; furthermore, the
average distance between two consecutive such eigenvalues is of
order $\e^2 j^{\frac 2k - 1} \simeq j^{-1} \simeq \e^{k}$. The other
resonant eigenvalues are of the form $\s_\a \simeq \a$ for $\a$
close to zero, namely when $\a = \e^2 l^{\frac{2}{k}}$ and $l$ is
sufficiently small (compared to, say, some negative power of $\e$).
Hence the distance from zero of the smallest eigenvalues of this
type is of order $\e^2$. Indeed, an accurate expansion in $\e$, see
Subsection \ref{ss:ae}, yields that this distance is bounded from
below by a multiple of $\e^2$ when $K$ is a non-degenerate minimal
submanifold.

{\bf 2) eigenfunctions:} as for the case of codimension $1$, it
turns out that the eigenfunctions corresponding to the $\eta_\a$'s
are of the form $\phi_j(\e y) u_j(\z)$, where $u_j$ is radial in the
variable $\z$ ($\z$ represent here some orthonormal coordinates in
the normal bundle of $K_\e$). The function $\phi_j$ instead
oscillates faster and faster as $\e$ tends to zero, since $j$ is of
order $\e^{-k}$. On the other hand it is possible to show, see
Subsection \ref{ss:model}, that the eigenfunctions corresponding to
the $\s_\a$'s are products $v_l(|\z|) \langle \z, \var_l \rangle_N$,
where $\langle \cdot, \cdot \rangle_N$ is the scalar product in $N
K$, and where $\var_l$ is a section of the normal bundle $N K_\e$,
and precisely an eigenfunction (scaled in $\e$) of the normal
Laplacian of $K$. Since the resonant modes correspond to low indices
$l$, $\var_l$ does not oscillate {\em as fast} as the resonant
$\phi_j$'s.

\

So far we considered an approximate operator, because in
\eqref{eq:modlin} we assumed a splitting of the metric into a
product. Since we expect to deal with small eigenvalues, a careful
analysis of the approximate solutions is needed (to apply local
inversion arguments), and also a refined understanding of the small
eigenvalues with the corresponding eigenfunctions.

Therefore we first try to obtain approximate solutions as accurate
as possible. For doing this, as in \cite{mal2, malm, malm2}, one can
introduce suitable coordinates on $\O_\e$ near $K_\e$, expand
formally \eqref{eq:tpe} in powers of $\e$, and solve it term by term
using functions of the form
\begin{equation}\label{eq:ukei}
    u_{I,\e}(y,\z) = \left[ w_0 + \e w_1 + \dots + \e^I w_I \right]
    (\e y, \z' + \Phi_0(\e y) + \dots + \e^{I-2} \Phi_{I-2}(\e y),
    \zn); \quad \z = (\z', \zn).
\end{equation}
Here $\Phi_0, \dots, \Phi_{I-2}$ represent smooth sections of the
normal bundle $N K$, and the functions $(w_i)_i$ are determined
implicitly via equations of the type
\begin{equation}\label{eq:modlin2}
    \begin{cases}
    - \D w_i + w_i - p w_0(\z) w_i
    = F_i(\e y, w_0, \dots, w_{i-1}, \Phi_0, \dots, \Phi_{i-2}) & \hbox{ in }
    \R^{n+1}_+, \\ \frac{\partial w_i}{\partial \nu} = 0 & \hbox{
    on } \partial \R^{n+1}_+.
  \end{cases}
\end{equation}
Notice that the operator acting on $w_i$ is nothing but the
linearization of \eqref{eq:w0} at $w_0$ (shifted in $\z'$ by $\Phi_0
+ \dots + \e^{I-2} \Phi_{I-2}$), which has an $n$-dimensional kernel
due to the invariance by translation in $\z'$. The functions
$\Phi_i$ are chosen in order to obtain orthogonality of $F_i$ to the
kernel, and to guarantee solvability in $w_i$. In doing this, the
non-degeneracy condition on $K$ comes into play, since the
$\Phi_i$'s solve equations of the form $\mathfrak{J} \Phi_i =
G_i(\oy)$. $\mathfrak{J}$ denotes the {\em Jacobi operator} of $K$,
related to the second variation of the volume functional, which is
invertible by the non-degeneracy assumption on the minimal
submanifold. Notice also that we wrote the variable $y$ with a
factor $\e$ on the front. This is in order to emphasize the slow
dependence in $y$ of these functions. In fact, recalling that (in
the model problem described above) resonance occurs mostly when
dealing with highly oscillating eigenfunctions, if we require slow
dependence in $y$ then there is no obstruction in solving
\eqref{eq:tpe} up to an arbitrary order $\e^I$.

Next one linearizes \eqref{eq:tpe} near the approximate solutions
just found. Compared to the above model problem, the eigenvalues
will be perturbed by some amount, due to the presence of the
corrections $(w_i)_i$ and to the geometry of the problem. In fact
the amount will be in general of order $\e$, since this is the size
of the corrections (from the $w_i$'s and the expansions of the
metric coefficients, see Lemma \ref{l:expgeuz}). This prevents a
direct control of the small eigenvalues of the linearized operator
(at $u_{I,\e}$) since, as discussed above, the characteristic size
of the spectral gaps at resonance are of order $\e^2$ or $\e^k$.

To overcome this problem, we look at the eigenvalues as functions of
$\e$. The counterparts of the numbers $\s_{\e^2 \o_l}$ can be again
obtained via a Taylor's expansion in $\e$, and they turn out to be
constant multiples of $\e^2$ times the eigenvalues of $\mathfrak{J}$
(up to an error of order $o(\e^2)$), so they are never zero. On the
other hand, the counterparts of the $\eta_{\e^2 \rho_j}$'s could
vanish for some values of $\e$ but, recalling the expansion
$\eta_{\e^2 \rho_j} \simeq \eta_0 + \e^2 j^{\frac 2k}$, one can hope
that generically in $\e$ none of these eigenvalues will be zero.

This is indeed shown using a classical theorem due to T. Kato, see
\cite{ka}, pag. 445, which allows us to estimate the derivatives of
the eigenvalues with respect to $\e$. To apply this result one needs
some control not only on the initial eigenvalues but also on the
corresponding eigenfunctions, and this is what basically the last
sections are devoted to. There we prove that if $\l = o(\e^2)$ is an
eigenvalue of the linearized operator, the eigenfunctions (up to a
small error) are linear combinations of products like $\phi_j(\e y)
u_j(\z)$, for $j \simeq \e^{-k}$ and for suitable functions $u_j$
radial in $\z$. Then we deduce that $\frac{\pa \l}{\pa \e}$ is close
to a number depending on $\e, N, p$ and $K$ only. As a consequence,
the spectral gaps near zero will {\em shift}, as $\e$ varies, almost
without squeezing, yielding invertibility for suitable values of the
parameter. This method also provides estimates on the norm of the
inverse operator, which blows-up with rate $\max\{\e^{-k},
\e^{-2}\}$ when $\e$ tends to zero, see Remark \ref{r:final}.

Finally, a straightforward application of the implicit function
theorem gives the desired result. To fix the ideas, when $p \leq
\frac{N+2}{N-2}$, solutions of \eqref{eq:tpe} can be found as
critical points of the following functional
\begin{equation}\label{eq:ie}
  J_\e (u) = \frac{1}{2} \int_{\O_\e} \left( |\n u|^2 + u^2
  \right) - \frac{1}{p+1} \int_{\O_\e} |u|^{p+1}, \qquad u \in
  H^1(\O_\e).
\end{equation}
One proves that $\|J'_\e(u_{I,\e})\|_{H^1(\O_\e)} \leq C_{I,k}
\e^{I+1-\frac k2}$ for $\e$ small. Even if the norm of the inverse
linear operator blows-up when $\e$ tends to zero, choosing $I$
sufficiently large (depending only on $k$ and $p$), one can find a
solution using the contraction mapping theorem near $u_{I,\e}$.

\

The general strategy of this proof, and especially Kato's theorem,
has been used in \cite{mal2}, \cite{malm} and \cite{malm2}, so
throughout the paper we will be sketchy in the parts where simple
adaptations apply. However the present setting requires some new
ingredients: we are going to explain next what are the differences
with respect to these and to some other related papers.  First of
all, compared to \cite{malm}, \cite{malm2}, where the case $k = N -
1$ was treated, here we need to characterize the limit set among all
the possible ones, since the codimension is higher, and this
reflects in the fact that the limit problem \eqref{eq:w0} is
degenerate. This requires to introduce the normal sections $\Phi_0,
\dots, \Phi_{I-2}$ in \eqref{eq:ukei}, and to use the non-degeneracy
condition on $K$.

The localization of the limit set has been indeed also faced in
\cite{mal2}. Here, apart from including that result as a particular
case, allowing higher dimensions and codimensions, we need a more
geometric approach. The main issue, as we already remarked, is that
we cannot use parameterizations with product sets as in
\eqref{eq:pardo}, since the normal bundle of $K$ is not trivial in
general. At this point some interplay between the analytic and
geometric features of the problem is needed. In particular the first
and second eigenfunctions of the linearization of \eqref{eq:w0} (the
profile of $\tilde{u}_\e$ at every point $q$ of $K$) can be viewed
of {\em scalar} or {\em vectorial} nature. More precisely, the
eigenfunction corresponding to the first eigenvalue is radial and
unique up to a scalar multiple. On the other hand the eigenfunctions
corresponding to the second eigenvalue have the symmetry of the
first spherical harmonics in the unit sphere of $N_q K$, and they
are in one-to-one correspondence with the vectors of $N_q K$. The
same holds true for the eigenfunctions of problem \eqref{eq:lina}
when $\a \geq 0$. When $q$ varies over the limit set, these
eigenfunctions (which are the resonant ones), depending on their
symmetry determine respectively a scalar function on $K$ or a
section of the normal bundle $N K$, on which the Laplace-Beltrami
operator or the normal Laplacian act naturally, see in particular
Section \ref{s:mod}. Apart from these considerations some other
difficulties arise, more technical in nature, due to the more
general character of the present result compared to that in
\cite{mal2}. Heavier computations are involved, especially since the
curvature tensors have more components, and some extra terms appear.
Anyway, some of the arguments have been simplified.

Finally, we should point out the differences with respect to the
papers \cite{dkw}, \cite{mmp}, \cite{mp}, where also special
solutions of the Nonlinear Schr\"odinger equation or constant mean
curvature surfaces are found. In \cite{dkw} and \cite{mp} the {\em
spectral gaps} are relatively big, and the eigenvalues can be
located using direct comparison arguments, so there is no need to
invoke Kato's theorem. In \cite{mmp} arbitrarily small spectral gaps
are allowed, but while there one has to study a partial differential
equation on a surface only, here we need to analyze the equation on
the whole space, which takes some extra work. Also, the Riemannian
manifold we consider here, $\pa \O$, has an {\em extrinsic}
curvature as a subset of $\R^N$, and therefore some error terms turn
out to be of order $\e$, and not $\e^2$, see Remark \ref{r:l3} {\em
(a)}. Nevertheless, we take great advantage of the geometric
construction in \cite{mmp}, especially in their choice of
coordinates near the limit set. We believe that our method could
adapt to study concentration at general manifolds for the Nonlinear
Schr\"odinger equation as well, as conjectured in \cite{amn}.

\

The paper is organized in the following way. We first introduce some
notations and conventions. In Section \ref{s:geom} we collect some
notions in differential geometry, like the Fermi coordinates near a
minimal submanifold, the normal Laplacian, the Laplace-Beltrami and
the Jacobi operators as well as the asymptotics of their
eigenvalues. In Section \ref{s:as} we construct the approximate
solution $u_{I,\e}$. In Section \ref{s:mod} we study some spectral
properties for the limit problem \eqref{eq:w0} (with some extension)
and we then derive a model for the linearized operator at
$u_{I,\e}$. In Section \ref{s:real} we turn then to the real
linearized operator: we construct some approximate eigenfunctions
which allow us to split our functional space as direct sum of
subspaces for  which the linearized operator is almost diagonal. In
Section \ref{s:applic}, using this splitting we characterize the
eigenfunctions corresponding to resonant eigenvalues. From these
estimates we can obtain invertibility, via Kato's theorem, and prove
our main result Theorem \ref{t:m}.

\

\begin{center}
{\bf Acknowledgments}
\end{center}

\noindent The authors are supported by MURST, under the project {\em
Variational Methods and Nonlinear Differential Equations}. They are
grateful to F. Pacard for some discussions concerning Remark
\ref{r:pac}. F.M. is grateful to SISSA for the kind hospitality.

\

\begin{center}
{\bf Notation and conventions}
\end{center}

\no - Dealing with coordinates, Greek letters like $\a, \b, \dots$,
will denote indices varying between $1$ and $N-1$, while capital
letters like $A, B, \dots$ will vary between $1$ and $N$; Roman
letters like $a$ or $b$ will run from $1$ to $k$, while indices like
$i, j, \dots$ will run between $1$ and $n := N - k - 1$.

\ms \no - $\z_{1}, \dots, \z_{n}, \z_{n+1}$ will denote coordinates
in $\R^{n+1}=\R^{N-k}$, and they will also be written as
$\z'=(\z_{1}, \dots, \z_{n})$, $\z=(\z',\zn)$.

\ms \no - The manifold $K$ will be parameterized with coordinates
$\oy = (\oy_1, \dots, \oy_k)$. Its dilation $K_\e := \frac 1 \e K$
will be parameterized by coordinates $(y_1, \dots, y_k)$ related to
the $\oy$'s simply by $\ov{y} = \e y$.

\ms \no - Derivatives with respect to the variables $\oy$, $y$ or
$\z$ will be denoted by $\pa_{\oy}$, $\pa_y$, $\pa_\z$, and for
brevity sometimes we might use the symbols $\pa_{\oa}$ and $\pa_i$
for $\pa_{\oy_a}$ and $\pa_{\z_i}$ respectively.

\ms \no - In a local system of coordinates, $(\ov{g}_{\a \b})_{\a
\b}$ are the components of the metric on $\pa \O$ naturally induced
by $\R^N$. Similarly, $(\ov{g}_{AB})_{AB}$ are the entries of the
metric on $\O$ in a neighborhood of the boundary. $(H_{\a \b})_{\a
\b}$ will denote the components of the mean curvature operator of
$\partial \O$ into $\R^N$.

\

\

\noindent Below, for simplicity, the constant $C$ is allowed to vary
from one formula to another, also within the same line, and will
assume larger and lager values. It is always understood that $C$
depends on $\O$, the dimension $N$ and the exponent $p$. It will be
explicitly written $C_l$, $C_\d$, $\dots$, if the constant $C$
depends also on other quantities, like an integer $l$, a parameter
$\d$, etc. Similarly, the positive constant $\g$ will assume smaller
and smaller values.

For a real positive variable $r$ and an integer $m$, $O(r^m)$ (resp.
$o(r^m)$) will denote a function for which $\left|
\frac{O(r^m)}{r^m} \right|$ remains bounded (resp. $\left|
\frac{o(r^m)}{r^m} \right|$ tends to zero) when $r$ tends to zero.
We might also write $o_\e(1)$ for a quantity which tends to zero as
$\e$ tends to zero. With $\mathcal{O}(r^m)$ we denote functions
which depend on the above variables $(\oy, \z)$, which are of order
$r^m$, and whose partial derivatives of any order, with respect to
the vector fields $\partial_\alpha$, $r\,\partial_i$, are bounded by
a constant times $r^m$.

$L_i$ will stand in general for a differential operator of order at
most $i$ in both the variables $\oy$ and $\z$ (unless differently
specified), whose coefficients are assumed to be smooth in $\oy$.

For summations, we might use the notation $\sum_{c}^d$ to indicate
that the sum is taken over an integer index varying from $[c]$ to
$[d]$ (the integer parts of $c$ and $d$ respectively). We might use
the same convention when we make an integer index vary between $c$
and $d$. We also use the standard convention of summing terms where
repeated indices appear.

We will assume throughout the paper that the exponent $p$ is at most
critical, namely that $p \leq \frac{N+2}{N-2}$, so that problem
\eqref{eq:pe} is variational in $H^1(\O)$. We will indicate at the
end what are the arguments necessary to deal with the general case.

\

\section{Geometric background}\label{s:geom}

In this section we list some preliminary notions in differential
geometry. First of all we introduce Fermi coordinates near a
submanifold of $\pa \O$, recall the definition of minimal
submanifold, and introduce the Laplace-Beltrami and the Jacobi
operators, together with some of their spectral properties. We refer
for example to \cite{aul} and \cite{spi} as basic references in
differential geometry.

\subsection{Fermi coordinates on $\pa\O$ near $K$}\label{ss:fc}
Let $K$ be a $k$-dimensional submanifold of $(\partial\O,\ov g)$
($1\le k\le N-1$) and set $n=N-k-1$ (see our notation). We choose
along $K$ a local orthonormal frame field $((E_a)_{a=1,\cdots
k},(E_i)_{i=1,\cdots, n})$ which is oriented. At points of $K$, $T
\pa \O$ splits naturally as $T K \oplus N K$, where $T K$ is the
tangent space to $K$ and $N K$ represents the normal bundle, which
are spanned respectively by $(E_a)_a$ and $(E_j)_j$.

Denote by $\n$ the connection induced by the metric $\ov{g}$ and by
$\n^N$ the corresponding normal connection on the normal bundle.
Given $q \in K$, we use some geodesic coordinates $\oy$ centered at
$q$. We also assume that at $q$ the normal vectors $(E_i)_i$, $i =
1, \dots, n$, are transported parallely (with respect to $\n^N$)
through geodesics from $q$, so in particular
\begin{equation}\label{eq:parall}
    \ov g\left(\nabla_{E_a}E_j\,,E_i\right)=0  \quad \hbox{ at } q,
    \qquad \quad i,j = 1, \dots, n, a = 1, \dots, k.
\end{equation}
In a neighborhood of $q$, we choose {\em Fermi coordinates} $(\oy,
\z)$ on $\pa \O$ defined by~
\begin{equation}\label{eq:fc}
    (\ov y,\z)\longrightarrow \exp^{\pa \O}_{\ov y}(
\sum\limits_{i=1}^{n}\,\z_i\,E_i); \qquad \quad (\oy, \z) =
\left((\ov y_a)_a,(\z_i)_i\right),
\end{equation}
where $\exp^{\pa \O}_{\oy}$ is the exponential map at $\oy$ in $\pa
\O$.

By our choice of coordinates, on $K$ the metric $\ov{g}$ splits in
the following way
\begin{equation}\label{eq:splitovg}
    \ov g(q) = \ov g_{ab}(q)\,d\ov y_a\otimes d\ov y_b+\ov
g_{ij}(q)\,d\z_i\otimes d\z_j; \qquad \quad q \in K.
\end{equation}
We denote by $\Gamma_a^b(\cdot)$ the 1-forms defined on the normal
bundle of $K$ by~
\begin{equation}\label{eq:Gab}
    \Gamma_a^b(E_i)=\ov g(\nabla_{E_a}E_b,E_i).
\end{equation}
We will also denote by $R_{\a\b\g\d}$ the components of the
curvature tensor with lowered indices, which are obtained by means
of the usual ones  $R_{\b\g\d}^\s$ by~
\[R_{\a\b\g\d}=\ov g_{\a\s}\,R_{\b\g\d}^\s.\]
When we consider the metric coefficients in a neighborhood of $K$,
we obtain a deviation from formula \eqref{eq:splitovg}, which is
expressed by the next lemma, see Proposition 2.1 in \cite{mmp} for
the proof. Denote by $r$ the distance function from $K$.

\begin{lemma} In the above coordinates $(\oy, \z)$,
for any $a=1,...,k$ and any $i,j=1,...,n$, we
have
\[
\begin{array}{rllll}\ov g_{ij}(0,\z)&=\delta_{ij}+\frac{1}{3}\,R_{istj}\,\z_s\,\z_t\,
+\,{\mathcal O}(r^3);\\[3mm]
\ov g_{aj}(0,\z)&={\mathcal O}(r^2);\\[3mm]
\ov
g_{ab}(0,\z)&=\delta_{ab}-2\,\Gamma_{a}^b(E_i)\,\z_i+\left[R_{sabl}+\Gamma_{a}^c(E_s)\,
\Gamma_{c}^b(E_l) \right]\z_s\z_l+{\mathcal O}(r^3).
\end{array}
\]
Here $R_{istj}$ are  computed at the point $q$ of $K$ parameterized
by $(0,0)$. \label{lemovg}
\end{lemma}

\subsection{Normal Laplacian, Laplace-Beltrami and Jacobi
operators}\label{ss:op}

In this subsection we recall some basic definitions and spectral
properties of differential operators associated to minimal
submanifolds. We first  recall some notions about the {\em
Laplace-Beltrami} operator,  the normal
connection and the {\em normal Laplacian}. \\
\no If $(M,g)$ is an $m$-dimensional  Riemannian manifold, the
Laplace-Beltrami operator on $M$ is defined in local coordinates by~
\begin{equation}\label{eq:laplcoord}
  \Delta_{
g}=\frac{1}{\sqrt{\det g}}\,\partial_A(\,\sqrt{\det g}\,{
g}^{AB}\,\partial_B\,),
\end{equation}
where the indices $A$ and $B$ runs in $1,\dots,m$, and where
$g^{AB}$ denote the components of the inverse of the matrix
$g_{AB}$.

\

\no Let $K \subseteq M$ be a  $k$-dimensional submanifold, $k\le
m-1$. The normal connection $\n^N$ on a normal vector field $V$ is
defined as the projection of the connection $\n V$ onto $N K$.
Moreover, one has the following formula regarding the horizontal
derivative of the product $\langle \cdot, \cdot \rangle_N$ in the
normal bundle (see \cite{spi}, Volume 4, Chapter 7.C, for further
details)
$$
  X \langle V, W \rangle_N = \langle \n^N_X V, W \rangle + \langle
  V, \n_X^N W \rangle,
$$
for any smooth sections $V$ and $W$ in $N K$. If we  choose an
orthonormal frame $(E_i)_i$ for $N K$ along $K$, we can write
$$
  \n^N_{\pa_{\ov{a}}}E_j = \b^l_j \left( \pa_{\ov{a}}
  \right) E_l,
$$
for some differential forms $\b^l_j$ (we recall our notation
$\pa_{\ov{a}}=\frac{\pa}{\pa \oy_a}$). Since the normal fields
$(E_i)_i$ are chosen to be orthonormal, it follows that for any
horizontal vector field $X$ there holds $X \langle E_i, E_j
\rangle_N = 0$, and hence one has
\begin{equation}\label{eq:antib}
    \b^l_j \left( \pa_{\ov{a}} \right) = - \b^j_l \left(
   \pa_{\ov{a}} \right) \qquad \qquad \forall \; l, j =
   1, \dots, n := m-k.
\end{equation}
This holds true, in particular, if we choose Fermi coordinates.
Since indeed the normal fields are extended via (normal) parallel
transport from $q$ to some neighborhood through the exponential map,
it follows that $\b^i_j(\pa_{\ov{a}})(0,0,\dots, \oy_a,0,\dots, 0) =
0$, and hence
\begin{equation}\label{eq:bfermiq}
  \b^l_j \left( \pa_{\ov{a}} \right) = 0 \quad \hbox{ at } q \qquad \qquad \forall
  \; a = 1, \dots, k, \hbox{ and } \forall \; l, j = 1, \dots, n;
\end{equation}
\begin{equation}\label{eq:bfermiq2}
    \pa_{\ov{a}} \left( \b^l_j \left( \pa_{\ov{a}} \right)
    \right) = 0 \quad \hbox{ at } q \qquad \qquad \forall
  \; a = 1, \dots, k, \hbox{ and } \forall \; l, j = 1, \dots, n.
\end{equation}

\no Recalling these facts, we can derive the expression of the {\em
normal Laplacian } in Fermi coordinates in the following way:
 given a normal vector field $V = V^j E_j$, there holds
$$
  \n^N_{\pa_{\ov{a}}} V = \pa_{\ov{a}} V^j\, E_j +
  V^j \b^l_j \left( \pa_{\ov{a}} \right) E_l.
$$
For any two normal vector fields $V$ and $W$ we have, by the
definition of $\D^N_K$
$$
  \int_K \langle \n^N V, \n^N W \rangle_N\, dV_{\ov{g}} = - \int_K
  \langle \D^N_KV, W \rangle_N \,dV_{\ov{g}}.
$$
We compute now the expression of $\D^N_K$ evaluating the left-hand
side and integrating by parts
\begin{eqnarray*}
    \int_K \langle \n^N V, \n^N W \rangle_N dV_{\ov{g}} & = &
   \int_K \left\langle \pa_{\ov{a}} V^j E_j + V^j \b^l_j \left(
   \pa_{\ov{a}} \right) E_l, \pa_{\ov{b}} W^i E_i + W^i
   \b^h_i \left( \pa_{\ov{b}} \right) \right\rangle_N
   \ov{g}^{ab} \sqrt{\det \ov{g}}  \\
   &=&  \int_K \left[\pa_{\ov{a}} V^i\, \pa_{\ov{b}} W^i +
  \pa_{\ov{a}} V^j W^i \b^j_i \left( \pa_{\ov{b}} \right)
  + V^j \b^i_j \left( \pa_{\ov{a}} \right) \,\pa_{\ov{b}} W^i
  \right. \\ &  &\qquad + \left. V^j W^i \b^l_j \left( \pa_{\ov{a}} \right) \b^l_i
  \left( \pa_{\ov{b}} \right)  \right] \ov{g}^{ab} \sqrt{\det \ov{g}}
\end{eqnarray*}
This quantity, for any $V$ and $W$, has to coincide with $- \int_K
(\D^N_K V)^i W^i \sqrt{\det \ov{g}}$, so we deduce that
\begin{eqnarray}\label{eq:bbbb} \nonumber
  (\D^N_KV)^i & = & \D_K (V^i) + \frac{1}{\sqrt{\det \ov{g}}}\,
  \pa_{\ov{b}}\left( V^j \b^i_j \left( \pa_{\ov{a}} \right) \ov{g}^{ab}
  \sqrt{\det \ov{g}} \right) \\
  & - & \ov{g}^{ab} \left( \pa_{\ov{a}} V^j \b^j_i \left(
  \pa_{\ov{b}} \right) + W^j \b^l_j \left( \pa_{\ov{a}} \right) \b^l_i \left(
  \pa_{\ov{b}} \right) \right) \sqrt{\det \ov{g}}.
\end{eqnarray}
In Fermi coordinates at $q$, which is  parameterized by $(0,0)$, we
have that
\begin{equation}\label{eq:ffcc}
    \ov{g}_{ab} = \d_{ab}, \quad\qquad \pa_{\ov{c}} \ov{g}_{ab} =
0\qquad \hbox{and}\qquad \pa_{\ov{c}} \sqrt{\det \ov{g}} = 0,
\end{equation}
and we also have \eqref{eq:bfermiq}-\eqref{eq:bfermiq2}. Hence the
last formula simplifies in the following way
\begin{equation}\label{eq:nlaplfer}
  (\D^N_KV)^i = \D_K (V^i) \qquad \qquad \hbox{ at } q.
\end{equation}

\

\no Let $C^\infty (NK)$ be the space of smooth normal vector fields
on $K$. For $\Phi \in C^\infty (NK)$, we can define the
one-parameter family of submanifolds $t \mapsto K_{t,\Phi}$ by
\begin{equation}
\label{eqKtPhi} K_{t,\Phi}:=\{ \exp_{\oy}^{\pa \O} (t \Phi(\oy))
\,:\, \oy\in K \}.
\end{equation}
The  first variation formula of the volume is the equation~
\begin{equation}
\label{eqfvf} \frac{d}{dt}\bigg|_{t=0}\mbox{Vol}(K_{t,\Phi})=\int_K
\langle\Phi, {\bf h}\rangle_N\,dV_K,
\end{equation}
where   ${\bf h}$ is the {\em mean curvature} (vector) of $K$ in
$\pa \O$, $\langle \cdot, \cdot \rangle_N$ denotes the restriction
of $\ov{g}$ to $N K$, and $dV_K$ the volume element of $K$.

\

\no The submanifold $K$ is said to be  {\em minimal} if it is a
critical point for the volume functional, namely if~
\begin{equation}
\label{eqms} \frac{d}{dt}\bigg|_{t=0}\mbox{Vol}(K_{t,\Phi})=0 \qquad
\mbox{for any $\Phi\in C^\infty (NK)$}
\end{equation}
or, equivalently by \eqref{eqfvf}, if the mean curvature ${\bf h}$
is identically zero on $K$. It is possible to prove that, if
$\G^b_a(E_i)$ is as in \eqref{eq:Gab}, then
\begin{equation}\label{eq:min}
    K \hbox{ is minimal } \qquad \Leftrightarrow \qquad \G^a_a(E_i)
    = 0 \quad \hbox{ for any } i = 1, \dots n.
\end{equation}
We point out that in the last formula we are summing over the index
$a$, which is repeated.

\

\no The {\em Jacobi operator} $\mathfrak{J}$ appears in the
expression of the second variation of the volume functional for a
minimal submanifold $K$
\begin{equation}
\label{eqjac}
\frac{d^2}{dt^2}\bigg|_{t=0}\mbox{Vol}(K_{t,\Phi})=-\int_K
\langle{\mathfrak J}\Phi,\Phi\rangle_N\,dV_K; \qquad \quad \Phi\in
C^\infty (NK),
\end{equation}
and is given by~
\begin{equation}\label{eq:jacobi}
    {\mathfrak
J}\Phi:=-\Delta_K^N\Phi+\mathfrak{R}^N\Phi-\mathfrak{B}^N\Phi,
\end{equation}
where $\mathfrak{R}^N, \mathfrak{B}^N :NK\rightarrow NK$ are defined
as~
$$
\mathfrak{R}^N\Phi=(R(E_a,\Phi)E_a)^N; \qquad \qquad \ov
g(\mathfrak{B}^N\Phi,n_K):=\G_b^a(\Phi)\G_a^b(n_K),
$$
for any unit normal vector $n_K$ to $K$. The operator $\D^N_K$ is
the {\em normal Laplacian} on $K$  defined in \eqref{eq:nlaplfer}.

\

\no A submanifold $K$ is said to be {\em non-degenerate} if the
Jacobi operator ${\mathfrak J}$ is invertible, or equivalently if
the equation ${\mathfrak J}\Phi=0$ has only the trivial solution
among the sections in $N K$.

\

\no We recall now some Weyl asymptotic formulas, referring for
example to \cite{cha}, or to \cite{LSY} %\cite{lieb}, \cite{NS}
and \cite{Min-P} for further details. Let $(M,g)$ be a compact
closed Riemannian manifold of dimension $m$, and let $\D_g$ be the
Laplace-Beltrami operator. Letting $(\rho_i)_i$, $i = 0, 1, \dots$,
denote the eigenvalues of $- \D_g$ (ordered to be non-decreasing in
$i$ and counted with their multiplicity), we have that
\begin{equation}\label{eq:weyl1}
    \rho_i \sim C_m \,\left(\frac{i}{Vol(M)}\right)^\frac{2}{m}\qquad \hbox{as}\quad i
\rightarrow \infty,
\end{equation}
where  $Vol (M)$ is the volume of $(M,g)$ and $C_m$ is a constant
depending only on the dimension $m$ (the Weyl constant). A similar
estimate, which can be proved using \eqref{eq:bbbb} and
\eqref{eq:weyl1}, holds for the normal Laplacian $\D^N_K$ on a
$k$-dimensional submanifold $K \subseteq M$. In fact, letting
$(\o_j)_j$, $j = 0, 1, \dots$, denote the eigenvalues of $- \D^N_K$
(still chosen to be non-decreasing in $j$ and counted with
multiplicity), one has
\begin{equation}\label{eq:weyl3}
  \o_j \sim C_{m,k}\,\left(\frac{j}{Vol(K)}\right)^\frac{2}{k}\qquad \hbox{as}\quad
  j \rightarrow \infty,
\end{equation}
where $C_{m,k}$ depends on the dimensions $m$ and $k$ only.

Considering the Jacobi operator $\mathfrak{J}$ for a minimal
submanifold $K$, it is easy to see from \eqref{eq:jacobi} that,
since $\mathfrak{J}$ differs from $- \D^N_K$ only by a bounded
quantity, we have the same asymptotic formula for its eigenvalues
$(\mu_l)_l$, and thereby
\begin{equation}\label{eq:weyl4}
  \mu_l \sim C_{m,k}\,\left(\frac{l}{Vol(K)}\right)^\frac{2}{k}\qquad \hbox{as}\quad
  l \rightarrow \infty.
\end{equation}
In the following, we let $(\phi_i)_i$ (resp. $(\var_j)_j$,
$(\psi_l)_l$) denote a base of eigenfunctions of $- \D_K$ (resp. of
$- \D^N_K$, $\mathfrak{J}$), normalized in $L^2(K)$ (resp. in
$L^2(K;NK)$), namely the set of functions (resp. normal sections of
$K$) satisfying
$$
  - \D_K \phi_i = \rho_i \phi_i; \qquad - \D^N_K \var_j = \o_j \var_j;
  \qquad \mathfrak{J} \psi_l = \mu_l
  \psi_l, \qquad \qquad i,j,l = 1, 2, \dots.
$$

Finally, using the eigenvalues $(\rho_j)_j$ and $(\mu_l)_l$, one can
express the $L^2$ norms, or the Sobolev norms of linear combinations
of the $\phi_j$'s and the $\psi_l$'s. In particular, if $f = \sum_j
\a_j \phi_j$, and if $g = \sum_l \b_l \psi_l$ are an $L^2$ function
and an $L^2$ normal section of $K$, and if $L_1 = \sum_\a c_\a(\oy)
\pa^\a_{\oy}$, $L_2 = \sum_\a \tilde{c}_\a(\oy) (\n^N_{\oy})^\a$ are
differential operators of order $d$ with smooth coefficients acting
on functions and normal sections respectively, then one has
\begin{equation}\label{eq:l2l2}
    \|L_1 f\|_{L^2(K)}^2 \leq C_{L_1} \sum_j (1 + \rho_j^d) \a_j^2;
    \qquad \qquad \|L_2 g\|_{L^2(K;N K)}^2 \leq C_{L_2} \sum_l (1
    + |\mu_l|^d) \b_l^2.
\end{equation}
An estimate similar to the latter one in \eqref{eq:l2l2} holds by
replacing the $\mu_l$'s by the $\o_j$'s, namely if $g' = \sum_j
\b'_j \var_j$, then $\|L_2 g'\|_{L^2(K;N K)}^2 \leq C_{L_2} \sum_j
(1 + |\o_j|^d) (\b_j')^2$.

\section{Approximate solutions to \eqref{eq:tpe}}\label{s:as}

In this section, given any positive integer $I$, we construct
functions $u_{I,\e}$ which solve \eqref{eq:tpe} up to an error of
order $\e^I$. We will find approximate solutions of \eqref{eq:tpe}
in the following form
\begin{equation}\label{eq:ukf}
 \chi_\e(|\z|)\bigg( w_0 \left( \z' + \Phi(\e y), \zn \right)
 + \e w_1 (\e y, \z' + \Phi(\e y), \zn)
  + \dots + \e^I w_I (\e y, \z' + \Phi(\e y), \zn)\bigg),
\end{equation}
where $\Phi(\e y) = \Phi_0(\e y) + \dots + \e^{I-2} \Phi_{I-2} (\e
y)$ and where the cutoff function $\chi_\e$ satisfies the properties
\begin{equation}
\label{eq:chi} \left\{
  \begin{array}{ll}
    \chi_\e(t)=1 & \hbox{for } t\in[0,\frac12\e^{-\g}],  \\
    \chi_\e(t)=0 & \hbox{for } t\in[\frac34\e^{-\g},\e^{-\g}], \\
    |\chi_\e^{(l)}(t)|\le C_l\,\e^{l\g}, & \hbox{}l\in \N.
  \end{array}
\right.
\end{equation}
Here $\Phi_0, \dots, \Phi_{I-2}$ are smooth vector fields from $K$
into $N K$, while $w_1, \dots, w_I$ are suitable functions
determined recursively by an iteration procedure. For doing this we
choose a system of coordinates in a neighborhood of $\pa \O_\e$ for
which the new metric coefficients can be expanded in powers of $\e$,
see Lemma \ref{l:expgeuz} below. In this way we can also expand
\eqref{eq:tpe} formally in powers of $\e$ and solve it term by term.
The functions $(w_i)_i$ will be obtained as solutions of an equation
arising from the linearization of \eqref{eq:w0} at $w_0$, while the
normal sections $(\Phi_i)_i$ will be determined using the
invertibility of the Jacobi operator. Notice that, by the
translation invariance of \eqref{eq:w0}, the linearized operator
possesses a non-trivial kernel, which turns out to be spanned by $\{
\pa_{\z_1} w_0, \dots, \pa_{\z_n} w_0\}$. The role of $\Phi_0, \dots
\Phi_{I-2}$ is to obtain at every step orthogonality to this kernel
and to solve the equation using Fredholm's alternative.

The method here is similar in spirit to the one used in \cite{mal2}
except for the fact that, working in higher dimensions and
codimensions, more geometric tools are needed. Therefore, we will
mainly focus on the new and geometric aspects of the construction,
omitting some details about the rigorous estimates on the error
terms, which can be handled as in \cite{mal2}.

\subsection{Choice of coordinates near $\pa
\O_\e$ and properties of approximate solutions}\label{ss:coord}

Let $\U_0 : \mathcal{U} \to \partial \O$, where $\mathcal{U} =
\mathcal{U}_1 \times \mathcal{U}_2 \subseteq \R^{k} \times \R^n$ is
a neighborhood of $0$ in $\R^{N-1}$, be a parametrization of
$\partial \O$ near some point $q \in K$ through the Fermi
coordinates $(\oy,\z)$ described before.

Let $\g \in (0,1)$ be a small number which, we recall, is allowed to
assume smaller and smaller values throughout the paper. Then for $\e
> 0$ we set
$$
  B_{\e,\g} = \left\{ x \in \R^{n+1}_+ \; : \; |x| < \e^{-\g}
  \right\}.
$$
Next we introduce a parametrization of a neighborhood (in $\O_\e$)
of $\frac q \e \in \pa \O_\e$ though the map $\U_\e$ given by
\begin{equation}\label{eq:fe}
    \U_\e(y,\z',\zn) = \frac{1}{\e} \U_0(\e y, \e \z') +
\zn \nu(\e y, \e \z'), \qquad x = (y,\z',\zn) \in \frac{1}{\e}
\mathcal{U}_1 \times B_{\e,\g},
\end{equation}
where $\e y = \oy$ and where $\nu(\e y, \e \z')$ is the inner unit
normal to $\partial \O$ at $\U_0(\e y, \e \z')$. We have
$$
\frac{\partial \U_\e}{\partial y_a} = \frac{\partial \U_0}{\partial
y_a}(\e y, \e \z') + \e \zn \frac{\partial \nu}{\partial
\ov{y}_a}(\e y, \e \z'); \qquad \qquad \frac{\partial
\U_\e}{\partial \z_i} = \frac{\partial \U_0}{\partial \z_i}(\e y, \e
\z') + \e \zn \frac{\partial \nu}{\partial \z_i}(\e y, \e \z').
$$
Using the equation
\begin{equation}\label{eq:dn}
    d \nu_x [v] =  {\bf H}(x)[v],
\end{equation}
we find
\begin{equation}\label{eq:dfe1}
  \frac{\partial \U_\e}{\partial y_a} = \left[ Id + \e \zn
   {\bf H}(\e y, \e \z') \right] \frac{\partial
\U_0}{\partial \ov y_a}(\e y, \e \z'); \qquad \qquad \frac{\partial
\U_\e}{\partial \z_i} = \left[ Id + \e \zn
  {\bf H}(\e y, \e \z') \right] \frac{\partial
\U_0}{\partial \z_i}(\e y, \e \z').
\end{equation}
Differentiating $\U_\e$ with respect to $\zn$ we also get
\begin{equation}\label{eq:dfe3}
  \frac{\partial \U_\e}{\partial \zn} = \nu(\e y, \e \z').
\end{equation}
Hence, letting $g_{AB}$ be the coefficients of the flat metric $g =
g_\e$ (we are emphasizing the role of the parameter $\e$ in the
entries, which is due to the dependence in $\e$ of the map
$\Upsilon_\e$) of $\R^N$ in the coordinates $(y,\z',\zn)$, with easy
computations we deduce that
\begin{equation}\label{eq:geij}
    g_{\a\b} (\tilde{y},\zn) = \ov{g}_{\a\b} (\e \tilde{y})
  + \e \zn \left( H_{\a\d} \ov{g}_{\d \b} + H_{\b \d} \ov{g}_{\d \a}
  \right) (\e \tilde{y}) + \e^2 \zn^2 H_{\a \d} H_{\s \b}
  \ov{g}_{\d \s} (\e \tilde{y}), \qquad  \tilde{y} = (y, \z');
\end{equation}
\begin{eqnarray}\label{eq:ge2233}
  g_{\a N} \equiv 0; \qquad \qquad g_{NN} \equiv 1.
\end{eqnarray}
Using the parametrization in \eqref{eq:fe}, a solution $u$ of
\eqref{eq:tpe} satisfies the equation~
\begin{equation}\label{eq:i'0}
  - \frac{1}{\sqrt{\det g}} \left[ \partial_B \left( g^{AB} \sqrt{\det g}
  \right) \right] \partial_A u - g^{AB} \partial^2_{AB} u + u - u^p = 0
  \quad \hbox{ in } \frac{1}{\e} \mathcal{U}_1 \times B_{\e,\g}
\end{equation}
with Neumann boundary  conditions on $\{\zn = 0\}$. Looking at the
term of order $\e^i$ in this equation, we will determine recursively
the functions $(w_i)_i$ and $(\Phi_{i-2})_i$ (defined in
\eqref{eq:ukf}) for $i = 1, \dots, I$. The specific choice of the
integer $I$, which will be determined later, will depend on the
dimension $N$ of $\O$, the dimension $k$ of $K$, and the exponent
$p$. For the moment we let it denote just an arbitrary integer. The
main result of this section is the following one.~
\begin{pro}\label{p:uke}
Consider the Euler functional $J_\e$ defined in \eqref{eq:ie} and
associated to problem \eqref{eq:tpe} (for $p \leq
\frac{n+k+2}{n+k-2}$). Then for any $I \in \N$ there exists a
function $u_{I,\e} : \O_\e \to \R$ with the following properties
\begin{equation}\label{eq:uke1}
  \|J'_\e(u_{I,\e})\|_{H^1(\O_\e)} \leq C_I \e^{I + 1 - \frac{k}{2}};
  \qquad u_{I,\e} \geq 0 \quad \hbox{ in } \O_\e; \qquad
  \frac{\partial u_{I,\e}}{\partial \nu} = 0 \quad \hbox{
  on } \partial \O_\e,
\end{equation}
where $C_I$ depends only on $\O$, $K$, $p$ and $I$. Moreover in the
above coordinates there holds
\begin{equation}\label{eq:uke2}
  \begin{cases}
    \left| \n_y^{(m)} u_{I,\e} (y,\z) \right| \leq C_{m,I} \e^m
    e^{-|\z|} P_I(\z), & \\[3mm]
    \left\| \n_y^{(m)} \n_{\z}
    u_{I,\e} (y,\z) \right\| \leq C_{m,I} \e^m e^{-|\z|} P_I(\z),
    &  \\[3mm]
     \left\| \n_y^{(m)} \n_{\z}^2 u_{I,\e} (y,\z)
    \right\| \leq C_{m,I} \e^m e^{-|\z|} P_I(\z),  & \end{cases}
    \qquad y \in \frac{1}{\e} \mathcal{U}_1, \z \in B_{\g,\e}, m = 0, 1, \dots,
\end{equation}
where $\n^{(m)}_y$ (resp. $\n^{(i)}_\z$) is any derivative of order
$m$ with respect to the $y$ variables (resp. of order $i$ with
respect to the $\z$ variables), where $C_{m,I}$ is a constant
depending only on $\O$, $K$, $p$ and $m$, and where $P_I(\z)$ are
suitable polynomials in $\z$.
\end{pro}

\no In the next subsection we show how to construct the approximate
solution $u_{I,\e}$ and we give some general ideas for the
derivation of the estimates in \eqref{eq:uke2}. We refer to
\cite{mal2} for rigorous and detailed proofs.

\subsection{Proof of Proposition \ref{p:uke}}\label{ss:pfpr31}

This subsection is devoted to the explicit construction of
$u_{I,\e}$. First of all we expand the Laplace-Beltrami operator
(applied to an arbitrary function $u$) in Fermi coordinates, and
then by means of this expansion we define implicity and recursively
the functions $(w_i)_i$ and the normal sections $(\Phi_i)_i$.

\subsubsection{Expansion of $\D_{g_\e} u$ in Fermi coordinates}\label{ss:e2}

We first provide a Taylor expansion of the coefficients of the
metric $g = g_\e$. {From} Lemma \ref{lemovg} and formula
\eqref{eq:geij} we have immediately the following result.

\begin{lem}\label{l:expgeuz}
For the (Euclidean) metric $g_\e$ in the above coordinates we have
the expansions
\begin{eqnarray*} &&g_{ij}=\delta_{ij} + 2 \e \zn H_{ij} +
\frac{1}{3} \e^2\,R_{istj}\,\z_s\,\z_t + \e^2 \zn^2 (H^2)_{ij}
\,+\,{\mathcal O}(\e^3|\z|^3); \\[3mm]
 &&g_{aj}= 2 \e \zn H_{aj} + {\mathcal O}(\e^2|\z|^2);\\[3mm]
 &&g_{ab}=  \delta_{ab}-2 \e \Gamma_{a}^b(E_i)\z_i + 2 \e \zn
H_{ab} + \e^2 \left[R_{sabl} +\Gamma_{a}^c(E_s) \Gamma_{c}^b(E_l)
\right]\z_s\z_l + \e^2 \zn^2 (H^2)_{ab}+ {\mathcal O}(\e^3|\z|^3);
\\[3mm]
 && g_{\a N} \equiv 0; \qquad \qquad \qquad g_{NN} \equiv 1.
\end{eqnarray*}
\end{lem}

\

\no Using these formulas, we are interested in expanding $\D_{g_\e}
u$ in powers of $\e$ for a function $u$ of the form
$$
  u(\oy, \z) = u(\e y, \z).
$$
Such a function represents indeed an {\em ansatz} for each term of
the sum in \eqref{eq:ukf}.

We recall that, when differentiating functions with respect to the
variables $y, \z$, we will mean that $\partial_a =
\partial_{y_a}$ and $\partial_i = \partial_{\z_i}$. When dealing
with the scaled variables $\oy$ we will write explicitly
$\partial_{\oy_a}$, so that, if $u$ is as above, we have $\partial_a
u(\e y, \z) = \e \partial_{\oy_a} u(\oy,\z)$.

\begin{lem}\label{l:expDgeu}
Given any positive integer $I$ and a function $u : \frac{1}{\e}
\mathcal{U}_1 \times B_{\e,\g} \to \R$ of the form $u(\e y,\z)$, we
have
\begin{eqnarray}\label{eq:explaplu} \nonumber
  \D_{g_\e} u & = & \partial^2_{ii} u + \partial^2_{\zn \zn} u + \e \left[
H_\a^\a \partial_{\zn} u - 2 \zn H_{ij} \partial^2_{ij} u \right] \\
  & + & \e^2 \left[ L_{2,1} u + L_{2,2} u + L_{2,3} u \right] +
  \sum_{i=3}^I \e^i L_i u + \e^{I + 1} \tilde{L}_{I+1} u,
\end{eqnarray}
where
$$
  L_{2,1} u = \partial^2_{\oy_a \oy_a} u - 4 \zn H_{ia} \pa^2_{\z_i \oy_a}
  u;
$$
$$
  L_{2,2} u = 3 \zn^2 (H^2)_{ij} \pa^2_{\z_i \z_j} u + 2
  \zn H_{ab} \G^a_b(E_i) \pa_i u - 2 \zn tr(H^2) \pa_{\zn} u;
$$
\begin{eqnarray*}
  L_{2,3} u & = & \left( R_{iaal} + \frac 13 R_{ihhl} \right) \z_l \pa_i u
  - \frac 13 R_{mijl} \z_m \z_l \pa^2_{\z_i \z_j} u - \frac 13
  R_{miji} \z_m \pa_{\z_j} u \\
  & - & \z_j \G^b_a(E_i) \G^a_b(E_j)
  \partial_{\z_i} u + 2 \z_i H_{ab} \G^a_b(E_i) \pa_{\zn} u,
\end{eqnarray*}
and where the $L_i$'s are linear operators of order $1$ and $2$
acting on the variables $\oy$ and $\z$ whose coefficients are
polynomials (of order at most $i$) in $\z$ uniformly bounded (and
smooth) in $\oy$. The operator $\tilde{L}_{I+1}$ is still linear and
satisfying the same properties of the $L_i$'s, except that its
coefficients are not polynomials in $\z$, although they are bounded
by polynomials in $\z$.
\end{lem}

\begin{pf}
The proof is simply based on a Taylor expansion of the metric
coefficients in terms of the geometric properties of $\pa \O$ and
$K$, as in Lemma \ref{l:expgeuz}. Recall that the Laplace-Beltrami
operator is given by~
\[ \Delta_{
g_\e}=\frac{1}{\sqrt{\det g_\e}}\,\partial_A(\,\sqrt{\det g_\e}\,{
g_\e}^{AB}\,\partial_B\,)\,,\] where indices $A$ and $B$ run between
$1$ and $N$. We can write~
\[\Delta_{g_\e}={g_\e}^{AB}\,\partial^2_{AB}+\left(\partial_A\,{
g_\e}^{AB}\right)\,\partial_B+\frac{1}{2}\,\partial_A(\,\log{\det
g_\e}\,)\,{g_\e}^{AB}\,\partial_B.\] Using the expansions of Lemma
\ref{l:expDgeu}, we easily see that~
\[
\begin{array}{rllll}
{g_\e}^{AB}\,\partial^2_{AB}u&=\partial^2_{\z_i
\z_i}u+\partial^2_{\zn \zn}u-2\,\e \zn
H_{ij}\partial^2_{\z_i\z_j}u\\[3mm]
&+\e^2\left\{\partial^2_{\ov y_a\ov y_a}
+\left(3\zn^2(H^2)_{ij}-\frac13\,R_{mijl}\z_m\z_l\right)\partial^2_{\z_i
\z_i}u-4\,\zn H_{ia}\partial^2_{\z_i \ov y_a}u \right\} +
\mathcal{O}(\e^3|\z|^3).
\end{array}
\]
We can also prove~
\begin{eqnarray*}
 \sqrt{\det g_\e}& = & 1+\e\zn H_\a^\a+ \frac16 \e^2 R_{miil}\z_m \z_l+\frac12
 \e^2 \bigg( R_{maal}+ \G_a^c(E_m)
\G_c^a(E_l)  \bigg)\z_m \z_l \\
 & + & \e^2\left\{\frac12 \zn^2  {(H_\a^\a)}^2 -\zn tr(H^2)
 +2\zn\z_iH_{ab}\G_b^a(E_i)-\z_i \z_j \G_a^b(E_i)  \G_b^a(E_j)        \right
 \} \\ & + & \mathcal{O}(\e^3|\z|^3),
\end{eqnarray*}
which gives
\[
\begin{array}{rllll}
\log \sqrt{\det g_\e}&=\e\zn H_\a^\a +\e^2\bigg\{
2\zn\z_iH_{ab}\G_b^a(E_i)-\zn^2\,tr(H^2)
-\z_i \z_j \G_a^b(E_i)  \G_b^a(E_j) \bigg\}\\[3mm]
&+\frac16\e^2 R_{miil}\z_m \z_l+\frac12\e^2\bigg( R_{maal}+
\G_a^c(E_m) \G_c^a(E_l)  \bigg)\z_m \z_l + \mathcal{O}(\e^3|\z|^3).
\end{array}
\]
Hence, we obtain~
\[
\begin{array}{rllll}
\pa_A\left(\log \sqrt{\det
g_\e}\right)\,g^{AB}\pa_B&=\e^2\bigg\{2\zn H_{ab}\G_b^a(E_i) - \z_j
\G_a^b(E_i)  \G_b^a(E_j) +\frac13 R_{mhhl}
\z_l+R_{iaal}\z_l\bigg\}\pa_i u\\[3mm]
&+\e H_\a^\a\pa_{\zn}u+\e^2\bigg\{ 2\z_lH_{ab}\G_b^a(E_l)-2\zn
tr(H^2)\bigg\}\pa_\zn u + \mathcal{O}(\e^3|\z|^3).
\end{array}
\]
 Collecting these formulas together, we obtain the desired result.
\end{pf}

\begin{rems}\label{r:l3}
(a) The term of order $\e$ in the expansion of $\D_g u$ in
\eqref{eq:explaplu} depends on the fact that $\pa \O$ has an
extrinsic curvature in $\R^N$. Such a term does not appear in the
analogous expansion for the mean curvature of tubes condensing on
minimal subvarieties of an abstract manifold, see Proposition 4.1 in
\cite{mmp} (where the small parameter $\rho$ is the counterpart of
our parameter $\e$).

(b) For later purposes, see for example Lemma \ref{l:sehu2}, it is
convenient to analyze in further detail the operator $L_3$ in
\eqref{eq:explaplu}, and in particular the coefficients of the
second derivatives in the $\oy$ variables. It follows from the above
expansions that the coefficient of $\pa^2_{\oy_a \oy_b}$ in $L_3$ is
given by
$$
  2 \left( \z_i \G^b_a(E_i) - \zn H_{ab} \right).
$$
\end{rems}

\subsubsection{Construction of the approximate
solution}\label{sss:steps}

We show now how to construct the approximate solutions of
\eqref{eq:tpe} via an iterative method. Given $I-2$ smooth vector
fields $\Phi_0, \dots, \Phi_{I-2}$ we define first the following
function $\hat{u}_{I,\e}$ on $K \times \R^{n+1}$, see \eqref{eq:ukf}
$$
  \hat{u}_{I,\e}(\oy, \z) = w_0(\z' + \Phi(\ov{y}),\zn) + \e
  w_1(\oy,\z' + \Phi(\ov{y}),\zn) + \dots + \e^{I} w_I(\oy,\z' +
  \Phi(\ov{y}),\zn),
$$
where $\Phi = \Phi_0 + \e \Phi_1 + \dots + \e^{I-2} \Phi_{I-2}$. In
the following, with an abuse of notation, we will consider
$\hat{u}_{I,\e}$ (and $w_0, \dots, w_I$) as functions of the
variables $y$ and $\z$ through the change of coordinates $\oy = \e
y$.

To define the functions $(w_j)_j$ and $(\Phi_j)_j$ we expand
equation \eqref{eq:i'0} formally in powers of $\e$ for $u =
\hat{u}_{I,\e}$ (using mostly Lemma \ref{l:expDgeu}) and we analyze
each term separately. Looking at the coefficient of $\e$ in the
expansion we will determine $w_1$, while looking at the coefficient
of $\e^j$ we will determine $w_j$ and $\Phi_{j-2}$, for $j = 2,
\dots, I$. In this procedure we use crucially the invertibility of
the Jacobi operator (recall that we are assuming $K$ to be
non-degenerate) and the spectral properties of the linearization of
\eqref{eq:w0} at $w_0$.

\

\no $\bullet$ {\bf Step 1: Construction of $w_1$}

\ms

\no We begin by taking $I = 1$ and $\Phi = 0$. {From} Lemma
\ref{l:expDgeu} we get {\em formally}
\begin{eqnarray*}
  - \D_{g_\e} \hat{u}_{1,\e} + \hat{u}_{1,\e} - \hat{u}_{1,\e}^p & = &
  - \D_{\R^{n+1}} w_0 + w_0 - w_0^p + \e \left( - \D_{\R^{n+1}} w_1 + w_1
  - p w_0^{p-1} w_1 \right) \\
   & - & \e \left[ H_\a^\a \partial_{\zn} w_0 - 2 \zn H_{ij} \partial^2_{ij}
   w_0 \right] + O(\e^2).
\end{eqnarray*}
The term of order $1$ (in the power expansion in $\e$) vanishes
trivially since $w_0$ solves \eqref{eq:w0}, and in order to make the
coefficient of $\e$ vanish, $w_1$ must satisfy the following
equation
\begin{equation}\label{eq:eqw1}
    \mathcal{L}_0 w_1 = H_\a^\a \partial_{\zn} w_0 - 2 \zn H_{ij} \partial^2_{ij}
    w_0,
\end{equation}
where $\mathcal{L}_0$ is the linearization of \eqref{eq:w0} at
$w_0$, namely
$$
  \left\{
    \begin{array}{ll}
      - \D w_1 + w_1 - p w_0^{p-1} w_1 = H_\a^\a \partial_{\zn}
      w_0 - 2 \zn H_{ij} \partial^2_{ij} w_0, & \hbox{ in } \R^{n+1}_+, \\
      \frac{\pa w_1}{\pa \zn} = 0, & \hbox{ on } \{ \zn = 0\}.
    \end{array}
  \right.
$$
Since $\mathcal{L}_0$ is self-adjoint and Fredholm on
$H^1(\R^{n+1}_+)$, the equation is solvable if and only if the
right-hand side is orthogonal to the kernel of $\mathcal{L}_0$,
namely if and only if the $L^2$ product of the right-hand side with
$\frac{\pa w_0}{\pa \z_i}$ vanishes for $i = 1, \dots, n$, see
Proposition \ref{p:nd0} below. This is clearly satisfied in our case
since both $\partial_{\zn} w_0$ and $\partial^2_{ij} w_0$ are even
in $\z'$, while the $\frac{\pa w_0}{\pa \z_i}$'s are odd in $\z'$
for every $i$. Besides the existence of $w_1$, from elliptic
regularity estimates we can prove its exponential decay in $\z$ and
its smoothness in $\oy$ (see for example Lemma 3.4 in \cite{mal2}).
Precisely, there exists a positive constant $C_1$ (depending only on
$\O, K$ and $p$) such that for any integer $\ell$ there holds
\begin{equation}\label{eq:estw1}
    |\n^{(\ell)}_{\oy} w_1(\oy,\z)| \leq C_1 C_l (1 + |\z|)^{C_1}
    e^{-|\z|}; \qquad \quad (\oy, \z) \in K \times \R^{n+1},
\end{equation}
where $C_l$ depends only on $l$, $p$, $K$ and $\O$.

\

\no $\bullet$ {\bf Step 2: Expansion at an arbitrary order}

\ms

\no We consider next  the coefficient of $\e^{\tilde{I}}$ for an
integer $\tilde{I}$ between $2$ and $I$, and we assume that the
functions $w_1, \dots, w_{\tilde{I}-1}$ and the vector fields
$\Phi_0, \dots, \Phi_{\tilde{I}-3}$ have been determined by
induction in $\tilde{I}$. The couple $(w_{\tilde{I}},
\Phi_{\tilde{I}-2})$ will be found reasoning as for $w_1$: in
particular an equation for $\Phi_{\tilde{I}-2}$ (solvable by the
invertibility of $\mathfrak{J}$) is obtained by imposing
orthogonality of some expression to the kernel of $\mathcal{L}_0$,
and then $w_{\tilde{I}}$ is found again with Fredholm's alternative.

Expanding \eqref{eq:i'0} with $u = \hat{u}_{I,\e}$, we easily see
that (formally), in the coefficient of $\e^{\tilde{I}}$, the
function $w_{\tilde{I}}$ appears as solution of the equation
\begin{equation}\label{eq:eqwtigen}
  \left\{
    \begin{array}{ll}
      \mathcal{L}_{\Phi} w_{\tilde{I}} = F_{\tilde{I}}(\oy, \z,
w_0, w_1, \dots, w_{\tilde{I}-1}, \Phi_0, \dots,
\Phi_{\tilde{I}-2}) & \hbox{ in } \R^{n+1}_+; \\
      \frac{\pa w_{\tilde{I}}}{\pa \zn} = 0 & \hbox{ on } \{\zn =
0\},
    \end{array}
  \right.
\end{equation}
where $\mathcal{L}_{\Phi}$ is defined by
$$
\mathcal{L}_{\Phi} u = - \D u + u - p w_0^{p-1}(\z'+\Phi(\oy),\zn)
u,
$$
and where $F_{\tilde{I}}$ is some smooth function of its arguments
(which we are assuming determined by induction). Our next goal is to
understand the role of $\Phi_{\tilde{I}-2}$ in the orthogonality
condition on $F_{\tilde{I}}$ (to the kernel of
$\mathcal{L}_{\Phi}$). In order to do this, we notice that, using
Lemma \ref{l:expDgeu} for $u = \hat{u}_{I,\e}$, the function $\Phi$
(precisely its derivatives in $\oy$) appears through the chain rule
when we differentiate $u$ with respect to the $\oy$ variables.
Moreover, for testing the orthogonality of the right-hand side in
\eqref{eq:eqwtigen} to the kernel of $\mathcal{L}_\Phi$, we have to
multiply it by the functions $\frac{\pa w_0}{\pa \z_i} (\z' +
\Phi(\oy),\zn)$, $i = 1, \dots, n$, so this condition will yield an
equation for $\Phi$ (and in particular for $\Phi_{\tilde{I}-2}$)
through a change of variables of the form $\z' \mapsto \z' +
\Phi(\oy)$.

Therefore, in the expansion of $\D_g \hat{u}_{I,\e}$, we focus only
on the terms (of order $\e^{\tilde{I}}$) containing either
derivatives with respect to the $\oy$ variables, which we collected
in $L_{2,1}$, or containing explicitly the variables $\z'$, which
are listed in $L_{2,3}$. In particular, none of these terms appear
in the first line of \eqref{eq:explaplu}.

Denoting the components of $\Phi$ by $(\Phi^j)_j$ (in the basis
$(E_j)_j$ of $N K$),  there holds
$$
  \partial_{\oy_a} \left( u(\oy,\z'+\Phi(\oy),\zn) \right) =
  \partial_{\oy_a} u(\oy,\z'+\Phi,\zn) + \frac{\pa \Phi^j}{\pa
  \oy_a} \frac{\pa u}{\pa \z_j}(\oy,\z'+\Phi(\oy),\zn);
$$
\begin{eqnarray*}
  \partial^2_{\oy_a \oy_a} \left( u(\oy,\z'+\Phi(\oy),\zn) \right) & = &
  \pa^2_{\oy_a \oy_a} u(\oy,\z'+\Phi,\zn) + 2 \frac{\pa \Phi^j}{
  \pa \oy_a} \pa^2_{\oy_a \z_j} u(\oy,\z'+\Phi,\zn) \\
   & + & \frac{\pa^2 \Phi^j}{\pa_{\oy_a \oy_a}} \frac{\pa u}{\pa
   \z_j}(\oy,\z'+\Phi(\oy),\zn) + \frac{\pa \Phi^j}{\pa \oy_a} \frac{\pa \Phi^l}{\pa \oy_a}
   \frac{\pa^2 u}{\pa \z_j \pa \z_l}(\oy,\z'+\Phi(\oy),\zn);
\end{eqnarray*}
$$
  \frac{\pa^2}{\pa \z_l \pa \oy_a} \left( u(\oy,\z'+\Phi(\oy),\zn)\right)
  = \pa^2_{\pa \z_l \pa \oy_a} u(\oy,\z'+\Phi,\zn) + \frac{\pa
  \Phi^j}{\pa \oy_a} \frac{\pa^2 u}{\pa \z_j \pa \z_l}(\oy,\z'+\Phi(\oy),\zn).
$$
Therefore, recalling the definition of $\hat{u}_{I,\e}$, since
$\pa_{\oy_a} w_0 = 0$ we find that
\begin{eqnarray*}
  L_{2,1} \hat{u}_{I,\e} & = & \frac{\pa^2 \Phi^j}{\pa^2_{\oy_a \oy_a}}
  \frac{\pa w_0}{\pa \z_j} + \frac{\pa \Phi^j}{\pa \oy_a} \frac{\pa
  \Phi^l}{\pa \oy_a} \frac{\pa^2 w_0}{\pa \z_j \pa \z_l} -
  4 \zn H_{la} \frac{\pa
  \Phi^j}{\pa \oy_a} \frac{\pa^2 w_0}{\pa \z_j \pa \z_l} \\
  & + & \sum_{i=1}^{\tilde{I}} \e^i \left\{ \pa^2_{\oy_a \oy_a} w_i + 2 \frac{\pa \Phi^j}{
  \pa \oy_a} \pa^2_{\oy_a \z_j} w_i +
  \frac{\pa^2 \Phi^j}{\pa^2_{\oy_a \oy_a}} \frac{\pa w_i}{\pa
   \z_j} + \frac{\pa \Phi^j}{\pa \oy_a} \frac{\pa \Phi^l}{\pa \oy_a}
   \frac{\pa^2 w_i}{\pa \z_j \pa \z_l} \right. \\
   & - & \left. 4 \zn H_{la} \left( \pa^2_{ \z_l \oy_a} w_i + \frac{\pa
  \Phi^j}{\pa \oy_a} \frac{\pa^2 w_i}{\pa \z_j \pa \z_l} \right) \right\}.
\end{eqnarray*}

\

\no $\bullet$ {\bf Step 3: Determining $w_{\tilde{I}}$ and
$\Phi_{\tilde{I}-2}$ for $\tilde{I} \geq 2$}

\ms

\no When we look at the coefficient of $\e^{\tilde{I}}$ in $\e^2
L_{2,1} \hat{u}_{I,\e}$, the terms containing $\Phi_{\tilde{I}-2}$
are given by
$$
  \frac{\pa^2 \Phi^j}{\pa^2{\oy_a \oy_a}}
  \frac{\pa w_0}{\pa \z_j} - 4 \zn H_{la} \frac{\pa
  \Phi^j}{\pa \oy_a} \frac{\pa^2 w_0}{\pa \z_j \pa \z_l} \quad \left(
  + \frac{\pa \Phi^j}{\pa \oy_a} \frac{\pa
  \Phi^l}{\pa \oy_a} \frac{\pa^2 w_0}{\pa \z_j \pa \z_l} \hbox{ if }
  \tilde{I} = 2 \right).
$$

When we project $\D_{g_\e} \hat{u}_{I,\e} - \hat{u}_{I,\e} +
\hat{u}_{I,\e}^p$ onto the kernel of $\mathcal{L}_{\Phi}$, namely
when we multiply this expression by $\frac{\pa w_0}{\pa
\z_s}(\z'+\Phi(\oy),\zn)$, $s = 1, \dots, n$, considering the terms
of order $\e^{\tilde{I}}$ involving $\Phi_{\tilde{I}-2}$, we have no
contribution from the first line and from $L_{2,2}$ in
\eqref{eq:explaplu} (with $u = \hat{u}_{I,\e}$), as explained in
Step 2. Also, in \eqref{eq:explaplu}, the factors of $\e^i$ for $i
\geq 3$, multiplied by $\e^{\tilde{I}-2} \Phi_{\tilde{I}-2}$ will
give higher order terms. In conclusion, we only need to pay
attention to $L_{2,1}$ and $L_{2,3}$.

When we multiply $\e^2 L_{2,3} w_0(\z'+\Phi,\zn)$ by $\frac{\pa
w_0}{\pa \z_s}(\z'+\Phi,\zn)$, $s = 1, \dots, n$, we can obtain the
coefficient of $\e^{\tilde{I}} \Phi_{\tilde{I}-2}^h$ in the
following way.

Looking for example at the first term in $\e^2 L_{2,3}$ we get
\begin{eqnarray*}
  & & \e^2 \int_{\R^{n+1}_+} \left( R_{iaal} + \frac 13 R_{ihhl} \right) \z_l
\pa_i w_0(\z' + \Phi, \zn) \pa_s w_0(\z'+\Phi,\zn) d \z \\ & = &
 \e^2 \int_{\R^{n+1}_+} \left( R_{iaal} + \frac 13 R_{ihhl} \right) (\z_l -
\Phi^l) \pa_i w_0(\z', \zn) \pa_s w_0(\z',\zn) d \z \\
  & = & \e^2 \int_{\R^{n+1}_+} \left( R_{iaal} + \frac 13 R_{ihhl}
\right) \z_l \pa_i w_0(\z', \zn) \pa_s w_0(\z',\zn) d \z \\ & - &
\e^2 \sum_{j=0}^{I-2} \e^j \Phi_j^l \int_{\R^{n+1}_+} \left(
R_{iaal} + \frac 13 R_{ihhl} \right) \pa_i w_0(\z', \zn) \pa_s
w_0(\z',\zn) d \z.
\end{eqnarray*}
Since $w_0$ is even in $\z'$, it follows by symmetry that the term
of order $\e^{\tilde{I}}$ containing $\Phi_{\tilde{I}-2}$ in the
last expression is given by
\begin{equation}\label{eq:rc1}
    - C_0 \left( R_{saal} + \frac 13 R_{shhl} \right) \Pti^l,
\end{equation}
where we have set
\begin{equation}\label{eq:C0}
    C_0 = \int_{\R^{n+1}_+} (\partial_{1} w_0)^2.
\end{equation}
{From} similar arguments, the third and the fourth terms in $L_{2,3}
w_0$ give respectively
\begin{equation}\label{eq:rc2}
    \frac 13 R_{lisi} C_0 \Pti^l,
\end{equation}
and
$$
  C_0 \G^b_a(E_s) \G^a_b(E_l) \Pti^l.
$$
The last term in $L_{2,3} w_0$ gives no contribution since the
coefficient of $\Phi_{\tilde{I}-2}$ vanishes by oddness, so it
remains to consider the second term. Integrating by parts we find
$$
  \frac 23 R_{mijl} \Pti^l \int_{\R^{n+1}_+} \z_m \partial_{\z_s}
  w_0 \partial^2_{\z_i \z_j} w_0 d \z \quad \left( + \pa^2_{\z_i\ov
y_a}\Pti^m \pa^2_{\z_j\ov y_a}\Pti^l
  \int_{\R^{n+1}_+} \partial^2_{\z_i \z_j} w_0 \partial_{\z_s} w_0
  d \z \quad \hbox{ if } \tilde{I} = 2 \right).
$$
In case $\tilde{I} = 2$ the quantity within round brackets cancels
by oddness, therefore in any case we only need to estimate the first
one. Still by oddness in $\z'$, the first integral is non-zero only
if, either $i = j$ and $m = s$, or $i = s$ and $j = m$, or $i = m$
and $j = s$.

In the latter case we have vanishing by the antisymmetry of the
curvature tensor in the first two indices. Therefore the only terms
left to consider are
$$
  \sum_{i} \frac 23 R_{siil} \Pti^l \int_{\R^{n+1}_+} \z_s \partial_{\z_s}
  w_0 \partial^2_{\z_i \z_i} w_0 d \z + \sum_{i} \frac 23 R_{isil} \Pti^l
  \int_{\R^{n+1}_+} \z_i \partial_{\z_s}
  w_0 \partial^2_{\z_s \z_i} w_0 d \z.
$$
Observe that, integrating by parts, when $s \neq i$ there holds
$$
  \int_{\R^{n+1}_+} \z_s \partial_{\z_s} w_0 \partial^2_{\z_i \z_i}
  w_0 d \z = - \int_{\R^{n+1}_+} \z_i \partial_{\z_s}
  w_0 \partial^2_{\z_s \z_i} w_0 d \z.
$$
Hence, still by the antisymmetry of the curvature tensor we are left
with
$$
  - \sum_{i} \frac 43 R_{siil} \Pti^l
  \int_{\R^{n+1}_+} \z_i \partial_{\z_s}
  w_0 \partial^2_{\z_s \z_i} w_0 d \z.
$$
The last integral can be computed with a further integration by
parts and is equal to $- \frac 12 C_0$, so we get
$$
  \frac 23 R_{siil} C_0 \Pti^l.
$$
This quantity cancels exactly with the second term in \eqref{eq:rc1}
and with \eqref{eq:rc2}.

\no When we multiply $\e^2 L_{2,1} w_0(\z'+\Phi,\zn)$ by $\frac{\pa
w_0}{\pa \z_s}(\z'+\Phi,\zn)$, $s = 1, \dots, n$, the terms
containing $\e^{\tilde{I}} \Phi_{\tilde{I}-2}^h$ are given by
\begin{eqnarray*}
  &&\int_{\R^{n+1}_+}\frac{\pa^2 \Phi_{\tilde{I}-2}^j}{\pa^2{\oy_a \oy_a}}
\pa_{\z_j} w_0\pa_{\z_s}
  w_0  d \z - 4 \int_{\R^{n+1}_+} \zn H_{la} \frac{\pa
  \Phi_{\tilde{I}-2}^j}{\pa \oy_a} \pa^2_{\z_j\z_l} w_0\partial_{\z_s}
  w_0  d \z \quad \\
  && \left(
  +  \int_{\R^{n+1}_+}\frac{\pa \Phi_{\tilde{I}-2}^j}{\pa \oy_a} \frac{\pa
  \Phi_{\tilde{I}-2}^l}{\pa \oy_a}
  \pa^2_{\z_j\z_l}w_0\,\partial_{\z_s}w_0
   d \z \hbox{ if }
  \tilde{I} = 2 \right),
\end{eqnarray*}
which give  by oddness
$$
 C_0\,\frac{\pa^2 \Phi_{\tilde{I}-2}^j}{\pa^2{\oy_a \oy_a}}.
$$
Collecting the above computations, we conclude that
$F_{\tilde{I}}(\oy, \z, w_0, w_1, \dots, w_{\tilde{I}-1}, \Phi_0,
\dots, \Phi_{\tilde{I}-2})$, the right-hand side of
\eqref{eq:eqwtigen}, is $L^2$-orthogonal to the kernel of
$\mathcal{L}_\Phi$ if and only if $\Phi_{\tilde{I}-2}$ satisfies an
equation of the form
$$
  C_0 \left( \frac{\pa^2 \Pti^s}{\pa_{\oy_a \oy_a}} - R_{saal}
  \Pti^l + \G^b_a(E_s) \G^a_b(E_l) \Pti^l \right)=G_{\tilde{I}-2}(\oy, \z,
w_0, w_1, \dots, w_{\tilde{I}-1}, \Phi_0, \dots,
\Phi_{\tilde{I}-3}),
$$
for some expression $G_{\tilde{I}-2}$. This equation can indeed be
solved in $\Pti$. In fact, observe that the operator acting on
$\Pti$ in the left hand side is nothing but the Jacobi operator,
which is invertible by the non-degeneracy condition on $K$.

\no Having defined $\Pti$ in this way, we turn to the construction
of $w_{\tilde I}$ which, we recall, satisfies equation
\eqref{eq:eqwtigen}. Having imposed the orthogonality condition, we
get again solvability and, as for $w_1$, one can prove the following
estimates
\begin{equation}\label{eq:estwtI}
    |\n^{(\ell)}_{\oy} w_{\tilde{I}}(\oy,\z)| \leq C_{\tilde{I}}
    C_l (1 + |\z|)^{C_{\tilde{I}}} e^{-|\z|};
    \qquad \quad (\oy, \z) \in K \times \R^{n+1},
\end{equation}
where $C_l$ depends only on $l$, $p$, $K$ and $\O$.

\

As already mentioned, we limit ourselves to the formal construction
of the functions $u_{I,\e}$, omitting the details about the rigorous
estimates of the error terms, which can be obtained reasoning as in
\cite{mal2}. We only mention that the number $\g$ has to be chosen
sufficiently small to obtain the positivity of $u_{I,\e}$, after we
multiply $\hat{u}_{I,\e}$ by the cutoff function $\chi_\e$, see
\eqref{eq:ukf} and \eqref{eq:chi}.

\section{A model linear problem}\label{s:mod}

In this section we consider a model for the linearized equation at
approximate solutions which, for $p \leq \frac{N+2}{N-2}$ (as we are
assuming until the last subsection), corresponds to
$J''_\e(u_{I,\e})$. We first study a one-parameter family of
eigenvalue problems, which include the linearization at $w_0$ of
\eqref{eq:w0}. Then we turn to the model for $J''_\e(u_{I,\e})$,
which can be studied, roughly, using separation of variables.

\subsection{Some spectral analysis in $\R^{n+1}_+$}\label{ss:spa}

In this subsection we consider a class of eigenvalue problems, being
mainly interested in the symmetries of the corresponding
eigenfunctions. We denote points of $\R^{n+1}$ by $(n+1)$-tuples
$\zeta_1, \zeta_2, \dots, \zeta_n, \z_{n+1} = (\z',\zn)$, and we let
$$
  \R^{n+1}_+ = \left\{ (\zeta_1, \zeta_2, \dots,
\zeta_n, \zn) \in \R^{n+1} \; : \; \zn > 0 \right\}.
$$
For $p \in \left( 1, \frac{n+3}{n-1} \right)$ ($\frac{n+3}{n-1}$ is
the critical exponent in $\R^{n+1}$) we consider problem
\eqref{eq:w0} which, we recall, is
$$
  \begin{cases}
    - \D u + u = u^p & \text{ in } \R^{n+1}_+, \\
    \frac{\partial u}{\partial \nu} = 0 & \text{ on } \partial
    \R^{n+1}_+, \\
    u > 0, u \in H^1(\R^{n+1}_+).
  \end{cases}
$$
It is well-known, see e.g. \cite{kwo}, that this problem possesses a
radial solution $w_0(r)$, $r^2 = \sum_{i=1}^{n+1} \zeta_i^2$, which
satisfies the properties
\begin{equation}\label{eq:dew0}
  \begin{cases}
    w'_0(r) < 0, & \text{ for every } r > 0, \\
    \lim_{r \to \infty} e^r r^{\frac{n}{2}} w_0(r) = \a_{n,p} > 0, &
    \lim_{r \to \infty} \frac{w'_0(r)}{w_0(r)} = - 1,
  \end{cases}
\end{equation}
where $\a_{n,p}$ is a positive constant depending only on $n$ and
$p$. Moreover, it turns out that all the solutions of \eqref{eq:w0}
coincide with $w_0$ up to a translation in the $\z'$ variables, see
\cite{gnn}, \cite{gnn2}.

Solutions of \eqref{eq:w0} can be found as critical points of the
functional $\ov{J}$ defined by
\begin{equation}\label{eq:I0}
  \ov{J}(u) = \frac{1}{2} \int_{\R^{n+1}_+} \left( |\n u|^2 + u^2 \right) -
  \frac{1}{p+1} \int_{\R^{n+1}_+} |u|^{p+1}; \qquad u \in H^1(\R^{n+1}_+).
\end{equation}
We have the following non-degeneracy result, see e.g. \cite{o}.

\begin{pro}\label{p:nd0}
The kernel of $\ov{J}''(w_0)$ is generated by the functions
$\frac{\partial w_0}{\partial \zeta_1}, \dots, \frac{\partial
w_0}{\partial \zeta_n}$. More precisely, there holds
$$
\ov{J}''(w_0) [w_0,w_0] = - (p-1) \|w_0\|_{H^1(\R^{n+1}_+)}^2,
$$
and
$$
\ov{J}''(w_0)[v,v] \geq C^{-1} \|v\|_{H^1(\R^{n+1}_+)}^2, \qquad
\forall v \in H^1(\R^{n+1}_+), v \perp w_0, \partial_{\zeta_1} w_0,
\dots,
\partial_{\zeta_n} w_0
$$
for some positive constant $C$. In particular, we have $\eta < 0$,
$\s = 0$ and $\t > 0$, where $\eta$, $\s$ and $\t$ are respectively
the first, the second and the third eigenvalue of $\ov{J}''(w_0)$.
Furthermore the eigenvalue $\eta$ is simple while $\s$ has
multiplicity $n$.
\end{pro}

Notice that, writing the eigenvalue equation $\ov{J}''(w_0) [u] = \l
u$ in $H^1(\R^{n+1}_+)$, taking the scalar product with an arbitrary
test function and integrating by parts one finds that $u$ satisfies
$$
\left\{
  \begin{array}{ll}
    - \D u + u - p w_0^{p-1} u = \l (- \D u + u) & \hbox{ in }
    \R^{n+1}_+,
    \\ \frac{\pa u}{\pa \zn} = 0 & \hbox{ on } \pa \R^{n+1}_+.
  \end{array}
\right.
$$
The goal of this subsection (the motivation will become clear in the
next one) is to study a more general version of this eigenvalue
problem, namely
\begin{eqnarray}\label{eq:eua}
  \begin{cases}
    - \D u + (1 + \a) u - p w_0^{p-1} u = \l \left(
    - \D u + (1 + \a) u \right) & \text{ in } \R^{n+1}_+, \\
    \frac{\partial u}{\partial \nu} = 0 & \text{ on } \partial
    \R^{n+1}_+,
  \end{cases}
\end{eqnarray}
where $\a \geq 0$. It is convenient to introduce the Hilbert space
(which coincides $H^1(\R^{n+1}_+)$, but endowed with an equivalent
norm)
$$
H_{\a} = \left\{ u \in H^1(\R^{n+1}_+) \, : \, \|u\|^2_{\a} =
\int_{\R^{n+1}_+} (|\n u|^2 + (1 + \a) u^2) \right\},
$$
with corresponding scalar product $( \cdot , \cdot )_\a$. We also
let $T_{\a} : H_\a \to H_\a$ be defined by duality in the following
way
\begin{equation}\label{eq:Tie}
  (T_\a u, v)_{H_\a} = \int_{\R^{n+1}_+} ((\n u \cdot \n v) + (1 + \a) u v) -
  p \int_{\R^{n+1}_+} w_0^{p-1} u v; \qquad u, v \in H_\a.
\end{equation}
When $\a = 0$, the operator $T_0$ is nothing but $\ov{J}''(w_0)$.
For $\a \geq 0$, the eigenfunctions of $T_\a$ satisfy
\eqref{eq:eua}. We want to study the first three eigenvalues of
$T_\a$ depending on the parameter $\a$.

\begin{pro}\label{p:ta}
Let $\eta_\a, \s_\a$ and $\t_\a$ denote the first three eigenvalues
of $T_\a$. Then $\eta_\a, \s_\a$ and $\t_\a$ are non-decreasing in
$\a$. For every value of $\a$, $\eta_\a$ is simple and there holds
$$
   \frac{\partial \eta_\a}{\partial \a} > 0; \qquad \qquad
   \qquad \lim_{\a \to + \infty} \eta_\a = 1.
$$
The eigenvalue $\s_\a$ has multiplicity $n$ and for $\a$ small it
satisfies $\frac{\partial \s_\a}{\partial \a} > 0$. The
eigenfunction $u_\a$ corresponding to $\eta_\a$ is radial in $\z$
and radially decreasing, while the eigenfunctions corresponding to
$\s_\a$ are spanned by functions $v_{\a,i}$ of the form
$v_{\a,i}(\z) = \hat{v}_\a(|\z|) \frac{\zeta_i}{|\z|}$, $i = 1,
\dots, n$, for some radial function $\hat{v}_\a(|\z|)$. If $u_\a$
and $v_\a$ are normalized so that $\|u_\a\|_\a = \|v_{\a,i}\|_\a =
1$, then they depend smoothly on $\a$. Moreover we have
$$
|\n^{(l)} u_\a(x)| + |\n^{(l)} (v_{\a,i})(x)| \leq C_l e^{-
\frac{|x|}{C_l}},
$$
provided $\a$ stays in a fixed bounded set of $\R$.
\end{pro}

\noindent Before proving the proposition we state a preliminary
lemma.

\begin{lem}\label{l:sphhar}
Let $\t$ denote the third eigenvalue of $\ov{J}''(w_0)$. Then, for
$\a \geq 0$, every eigenfunction corresponding to an eigenvalue $\l
\leq \frac \t 2$ of \eqref{eq:eua} is either radial and corresponds
to the least eigenvalue, or is a radial function times a first-order
spherical harmonic (in the angular variable $\th = \frac{\z}{|\z|}$)
with zero coefficient in $\z'$, and correspond to the second
eigenvalue.
\end{lem}

\begin{pf}
First of all we notice that, extending evenly across $\pa
\R^{n+1}_+$ any function $u \in H^1(\R^{n+1}_+)$ which is a solution
of \eqref{eq:eua}, we obtain a smooth entire solution of $- \D u +
(1 + \a) u - p w_0^{p-1} u = \l \left( - \D u + (1 + \a) u \right)$.
Next, we decompose $u$ in spherical harmonics in the angular
variable $\th$ (we are using only spherical harmonics which are even
in $\zn$)
$$
  u = \sum_{i=0}^\infty u_i(|\z|) Y_{i,e}(\th); \qquad \quad \z \in
  \R^{n+1}, \th = \frac{\z}{|\z|} \in S^n.
$$
Here $Y_{i,e}$ is the $j-th$ eigenfunction of $- \D_{S^n}$ (which is
even in $\zn$), namely it satisfies $\D_{S^n} Y_{i,e} =
\l_{i,e}^{S^n} Y_{i,e}$, where we have denoted by $\l_{i,e}^{S^n}$
the $i$-th eigenvalue of $- \D_{S^n}$ on the space of even functions
in $\zn$. In particular, the function $Y_{0,e}$ is constant on $S^n$
and correspond to $\l_{1,e}^{S^n} = 0$, while $\l_{2,e}^{S^n} = n$
has multiplicity $n$. The eigenfunctions corresponding to
$\l_{2,e}^{S^n}$ are (up to a constant multiple) the restrictions,
from $\R^{n+1}$ to $S^n$, of the linear functions in $\z'$.

The laplace equation in polar coordinates writes as
$$
  \D_{\R^{n+1}} u = \D_r u + \frac{1}{r^2} \D_{S^n} u,
$$
where $\D_r = \frac{d^2}{d r^2} + \frac nr \frac{d}{d r}$.
Therefore, if $u = \sum_{i=0}^\infty u_i(|\z|) Y_{i,e}(\th)$ is a
solution of \eqref{eq:eua}, then every radial component $u_i$
satisfies the equation
\begin{equation}\label{eq:laj}
    \left\{
      \begin{array}{ll}
        (1-\l) \left( - v'' - \frac n r v ' + \left( 1 + \a +
        \frac{\l_{i,e}^{S^n}}{r^2} \right)
  v  \right) - p w_0^{p-1} v = 0   & \hbox{ in } \R_+; \\
        v'(0) = 0. &
      \end{array}
    \right.
\end{equation}
We also notice that, since the space of functions $\{ v(r)
Y_{i,e}(\th) \}$ (for a fixed $i$) is sent into itself by the
Laplace operator, every Fourier component (in the angular variables)
of an eigenfunction  of \eqref{eq:eua} is still an eigenfunction.

We call $\l_{\a,i,j}$ the $j$-th eigenvalue of \eqref{eq:laj}. From
Proposition \ref{p:nd0} it follows that $\l_{0,1,1} = -(p-1) < 0$
and that $\l_{0,1,j} > \t$ for $j \geq 2$. In fact, a radial
eigenfunction of $\ov{J}''(w_0)$ which is not (a multiple of) $w_0$
itself must correspond to an eigenvalue greater or equal than $\t$,
which is positive. On the other hand, it follows from Proposition
\ref{p:nd0} that $\l_{0,2,1} = 0$, and also that $\l_{0,2,j} \geq \t
> 0$ for $j \geq 2$. Finally, since $\l_{0,i,1} \geq \t > 0$
for $i \geq 3$, we have in addition $\l_{0,i,j} \geq \t$ for every
$i \geq 3$ and for every $j \geq 1$.

After these considerations, we turn to the case $\a > 0$, for which
similar arguments will apply. Solutions of \eqref{eq:laj} can be
found as extrema (minima, for example) of the Rayleigh quotient
\begin{equation}\label{eq:rq}
  \frac{\int_{\R_+} r^n \left[ (v')^2 + \left( 1 + \a +
  \frac{\l_{i,e}^{S^n}}{r^2} \right) v^2 \right] - p
  \int_{\R_+} r^n w_0^{p-1} v^2}{\int_{\R_+} r^n \left[ (v')^2
  + \left( 1 + \a + \frac{\l_{i,e}^{S^n}}{r^2} \right)
  v^2 \right]}
\end{equation}
from a standard min-max procedure. Using elementary inequalities it
is easy to see that the above quotient is non-decreasing in $\a$.
Therefore it follows that $\l_{\a,1,j} > 0$ for $j \geq 2$, that
$\l_{\a,2,j} \geq \t > 0$ for $j \geq 2$ and that $\l_{\a,i,j} \geq
\t$ for every $i \geq 3$ and for every $j \geq 1$. This concludes
the proof.
\end{pf}

\

\begin{pfn} {\sc of Proposition \ref{p:ta}}
The simplicity of $\eta_\a$ can be proved as in \cite{malm}, Section
3, using spherical rearrangements and the maximum principle. The
weak monotonicity in $\a$ of the eigenvalues can be easily shown
using the Rayleigh quotient in the space $H_\a$, as for
\eqref{eq:rq}.

The smoothness of $\a \mapsto \eta_\a$ and of $\a \mapsto u_\a$ can
be deduced in the following way. Since the two spaces
$H^1(\R^{n+1}_+)$ and $H_\a$ coincide, and since the eigenvalues of
an operator do not depend on the choice of the (equivalent) norms,
we can consider $T_\a$ acting on $H^1(\R^{n+1}_+)$ endowed with its
standard norm (independent of $\a$). Having fixed the space, we
notice that the explicit expression of $T_\a$ is given by
\begin{equation}\label{eq:Taue}
    T_\a u = \left[ - \D + 1 \right]^{-1} \left( - \D u + (1 + \a) u -
  p w_0^{p-1} u \right).
\end{equation}
In fact, letting $T_\a u = q \in H^1(\R^{n+1}_+)$, taking the scalar
product with any $v \in H^1(\R^{n+1}_+)$ and using \eqref{eq:Tie} we
find
$$
  \int_{\R^{n+1}_+} \left[ (\n q \cdot \n v) + q v \right] =
  \int_{\R^{n+1}_+} \left[(\n u \cdot \n v) + (1 + \a) u v\right] - p
  \int_{\R^{n+1}_+} w_0^{p-1} u v,
$$
which leads to \eqref{eq:Taue} by the arbitrarity of $v$. It is
clear that the operator in \eqref{eq:Taue} depends smoothly on $\a$
and therefore, being $\eta_\a$ simple, the smooth dependence on $\a$
of $\eta_\a$ and $u_\a$ follows.

We now compute the derivative of $\eta_\a$ with respect to $\a$. The
function $u_\a$ satisfies
\begin{equation}\label{eq:eva}
  \begin{cases}
    (1 - \eta_\a) \left( - \D u_\a + (1 + \a) u_\a \right) =
    p w_0^{p-1} u_\a & \text{ in } \R^{n+1}_+, \\ \frac{\partial
    u_\a}{\partial \nu} = 0 & \text{ on } \partial \R^{n+1}_+.
  \end{cases}
\end{equation}
Differentiating with respect to $\a$ the equation $\|u_\a\|_\a^2 =
1$,  we find
\begin{equation}\label{eq:dnva}
  \frac{d}{d \a} \|u_\a\|^2_\a = 0 \qquad \qquad \Rightarrow
  \qquad \qquad \left( \frac{d u_\a}{d \a}, u_\a \right)_\a = -
  \int_{\R^{n+1}_+} u_\a^2.
\end{equation}
On the other hand, differentiating \eqref{eq:eva}, we obtain
\begin{equation}\label{eq:edva}
  \begin{cases}
    - \frac{d \eta_\a}{d \a} \left( - \D u_\a + (1 + \a)
    u_\a \right) + (1 - \eta_\a) \left( - \D \left( \frac{d u_\a}{d \a}
    \right) + (1 + \a) \frac{d u_\a}{d \a} + u_\a \right) = p
    w_0^{p-1} \frac{d u_\a}{d \a} & \text{ in } \R^{n+1}_+,  \\
    \frac{\partial}{\partial \nu} \left( \frac{d u_\a}{d \a}
    \right) = 0  & \text{ on } \partial \R^{n+1}_+.
  \end{cases}
\end{equation}
Multiplying \eqref{eq:edva} by $u_\a$, integrating by parts and
using \eqref{eq:dnva}, one gets
\begin{equation}\label{eq:dma}
  \frac{d \eta_\a}{d \a} = (1 - \eta_\a) \int_{\R^{n+1}_+} u_\a^2 > 0.
\end{equation}
Indeed, since $T_\a \leq Id_{H^1(\R^{n+1}_+)}$, every eigenvalue of
$T_\a$ is strictly less than $1$, and in particular $(1 - \eta_\a)
> 0$.
We now consider the second eigenvalue $\s_\a$. For any $\a \geq 0$
it is possible to make a separation of variables, finding
eigenfunctions of \eqref{eq:eua} of the form $Y_{i,e}
\hat{v}_{\a,i}$, where $Y_{i,e} = \frac{\z_i}{|\z|}$, $i = 1, \dots,
n$, correspond to $\l_{2,e}^{S^n}$. Also, from Lemma \ref{l:sphhar}
we know that for $\a$ close to $0$ (indeed, as long as $\s_\a < \t$)
every eigenfunction corresponding to $\s_\a$ is of this form, for
some $i \in \{1, \dots, n\}$. Therefore, if we restrict ourselves to
the space of functions of the form $\hat{v}(|\z|) \frac{\z_i}{|z|}$
for a fixed $i \in \{1, \dots, n\}$, the first eigenvalue for
\eqref{eq:eua} becomes simple, so we can reason as before, obtaining
smoothness in $\a$ and the strict monotonicity of $\s_\a$.

We prove next that the eigenvalue $\eta_\a$ converges to $1$ as $\a
\to + \infty$. There holds
$$
\eta_\a = \inf_{u \in H_\a} \frac{\int_{\R^{n+1}_+} \left[ |\n u|^2
+ (1 + \a) u^2 - p w_0^{p-1} u^2 \right]}{\int_{\R^{n+1}_+} \left[
|\n u|^2 + (1 + \a) u^2 \right]}.
$$
Fixing any $\d > 0$, it is sufficient to notice that
$$
  |\n u|^2 + \left( (1 + \a) - p w_0^{p-1} \right) u^2 \geq (1 - \d)
  \left[ |\n u|^2 + (1 + \a) u^2 \right] \qquad \hbox{ for every } u,
$$
provided $\a$ is sufficiently large. This concludes the proof of the
claim.

The decay on $u_\a$, $v_{\a,i}$ and their derivatives is standard
and can be shown as in \cite{malm}, so we do not give details here.
\end{pfn}

%\
%
%\noindent In view of formula \eqref{eq:dma}, for later applications
%it is convenient for us to define, given $\a \geq 0$, the function
%$\tilde{F}(\a)$ as
%\begin{equation}\label{eq:fa}
%  \tilde{F}(\a) = 2 \a (1 - \eta_\a) \int_{\R^{n+1}_+} u_\a^2 > 0.
%\end{equation}

\begin{rem}\label{r:a0}
Proposition \ref{p:ta} implies in particular that there is a unique
$\ov{\a} > 0$ such that $\eta_{\ov{\a}} = 0$. Moreover, we have also
$$
u_0 = \tilde{C}_0 w_0; \qquad \qquad v_0^h = \ov{C}_0 \partial_h
w_0,
$$
for some positive constants $\tilde{C}_0$ and $\ov{C}_0$.
\end{rem}

\

\noindent We also need to introduce a variant of the eigenvalue
problem \eqref{eq:eua}, for which we impose vanishing of the
eigenfunctions outside a certain set. For $\e > 0$ and for $\g \in
(0,1)$ we define
\begin{equation}\label{eq:je}
    B_{\e,\g} = \left\{ x \in \R^{n+1}_+ \; : \; |x| < \e^{-\g} \right\},
\end{equation}
and let
$$
H^1_{\e} = \left\{ u \in H^1(B_{\e,\g}) \; : \; u (x) = 0 \text{ for
} |x| = \e^{-\g} \right\}.
$$
We let $H_{\a,\e}$ denote the space $H^1_\e$ endowed with the norm
$$
  \|u\|_{\a,\e}^2 = \int_{B_{\e,\g}} \left[ |\n u|^2 + (1 + \a) u^2
  \right]; \qquad u \in H^1_\e,
$$
and the corresponding scalar product $( \cdot , \cdot )_{\a,\e}$.
Similarly, we define $T_{\a,\e}$ by
$$
(T_{\a,\e} u, v)_{\a,\e} = \int_{B_{\e,\g}} \left[ (\n u \cdot \n v)
+ (1 + \a) u v - p w_0^{p-1} u v \right]; \qquad \qquad u, v \in
H_{\a,\e}.
$$
The operator $T_{\a,\e}$ satisfies properties analogous to $T_\a$.
We list them in the next Proposition, which also gives a comparison
between the first eigenvalues and eigenfunctions of $T_\a$ and
$T_{\a,\e}$.

\begin{pro}\label{p:eicompe}
There exists $\e_0 > 0$ such that for $\e \in (0, \e_0)$ the
following properties hold true. Let $\eta_{\a, \e}$, $\s_{\a, \e}$
and $\t_{\a,\e}$ denote the first three eigenvalues of $T_{\a,\e}$.
Then $\eta_{\a, \e}$, $\s_{\a, \e}$ and $\t_{\a,\e}$ are
non-decreasing in $\a$. For every value of $\a$, $\eta_{\a, \e}$ is
simple and $\frac{\partial \eta_{\a,\e}}{\partial \a} > 0$. For $\a$
sufficiently small, $\s_{\a,\e}$ has multiplicity $n$ and
$\frac{\partial \s_{\a,\e}}{\partial \a} > 0$. The eigenfunction
$u_{\a,\e}$ corresponding to $\eta_{\a,\e}$ is radial in $\z$ and
radially decreasing, while the eigenfunctions corresponding to
$\s_{\a,\e}$ are spanned by functions $v_{\a,\e,i}$ of the form
$v_{\a,\e,i}(\z) = \hat{v}_{\a,\e}(|\z|) \frac{\zeta_i}{|\z|}$, $i =
1, \dots, n$, for some radial function $\hat{v}_{\a,\e}(|\z|)$. The
eigenvector $u_{\a,\e}$ (resp. $v_{\a,\e,i}$), normalized with
$\|u_{\a,\e}\|_{H_{\a,\e}} = 1$ (resp.
$\|v_{\a,\e,i}\|_{H_{\a,\e,i}} = 1$) corresponding to $\eta_{\a,\e}$
(resp. $\s_{\a,\e}$ for $\a$ small) depend smoothly on $\a$.
Moreover for some fixed $C > 0$ there holds
\begin{eqnarray}\label{eq:niuae}
  |\n^{(l)} u_{\a,\e}(\z)| + |\n^{(l)} v_{\a,\e,i}(\z)| \leq C_l e^{-
\frac{|\z|}{C_l}}, \qquad \hbox{ for } i = 0, \dots, n;
\end{eqnarray}
\begin{equation}\label{eq:muaae}
  |\eta_\a - \eta_{\a,\e}| + \|u_\a - u_{\a,\e}\|_{H^1(\R^{n+1}_+)}
  + |\s_\a - \s_{\a,\e}| + \|v_{\a,i} - v_{\a,\e,i}\|_{H^1(\R^{n+1}_+)}
  \leq C e^{- \frac{\e^{-\g}}{C}},
\end{equation}
provided $\a$ stays in a fixed bounded set of $\R$. The functions
$u_{a,\e}$ and  $v_{\a,\e,i}$ in this formula have been set
identically $0$ outside $B_{\e,\g}$. Furthermore, $\t_{\a,\e} \geq
\t_\a \geq \t$ for every value of $\a$ and $\e$.
\end{pro}

\noindent The proof follows that of Proposition 2.3 in \cite{malm2},
and hence we omit it here. It is still based on some elementary
inequalities and on the Rayleigh quotient. The quantitative
estimates in \eqref{eq:muaae} can be deduced using cutoff functions
and the Green's representation formula for the operator $- \D + (1 +
\a)$ in $\R^{n+1}_+$.

As a consequence of this proposition (taking $\a = 0$) we obtain
that, if (for $\e$ small) $u \in H^1_\e$ has no Fourier components
(in $\th$) with indices less or equal to $n$, then $(T_{0,\e} u,
u)_{0,\e} \geq \frac{\t}{2} (u,u)_{0,\e}$. Equivalently, there holds
\begin{equation}\label{eq:varo}
    p \int_{B_{\e,\g}} w_0^{p-1}(|\z|) u^2 \leq \left( 1 - \frac \t
    2 \right)
    \int_{B_{\e,\g}} \left( - \D u + u \right) u d\z \qquad \quad
    \hbox{ for any } u = \sum_{i=n+1}^\infty u_j(|\z|) Y_{i,e}(\th),
    u \in H^1_\e.
\end{equation}

\subsection{A model for $J''_\e(u_{I,\e})$}\label{ss:model}

In this subsection, using the analysis of the previous one, we
construct a model operator which, up to some extent, mimics the
properties of $J''_\e(u_{I,\e})$, and for which we can give an
explicit description of the spectrum. Although the related
construction in \cite{mal2} is a particular case of the one made
here, the general spirit is quite different, and is more geometric
in nature.

\

\no First of all, we choose an orthonormal frame $(E_i)_i$ as
before, and we define a metric $\hat{g}$ on $N K$ as follows. For $v
\in N K$, a tangent vector $V \in T_v N K$ can be identified with
the velocity of a curve $v(t)$ in $N K$ which is equal to $v$ at
time $t = 0$. The same holds true for another tangent vector $W \in
T_v N K$. Then the metric $\hat{g}$ on $N K$  is defined on the
couple $(V, W)$ in the following way (see \cite{do2}, pag. 79)
$$
  \hat{g}(V,W) = \ov{g} \left( \pi_* V, \pi_* W \right) + \left\langle
  \frac{D^N v}{d t}|_{t=0}, \frac{D^N w}{d t}|_{t=0} \right\rangle_N.
$$
In this formula $\pi$ denotes the natural projection from $N K$ onto
$K$, and $\frac{D^N v}{d t}$ denotes the (normal) covariant
derivative of the vector field $v(t)$ along the curve $\pi \; v(t)$.
In the notation of Subsection \ref{ss:op} we have that, if $v(t) =
v^j(t) E_j(t)$, then
$$
  \frac{D^N v}{d t} = \frac{d v^j(t)}{dt} E_j(t) + v^j(t) \b_j^l
  \left( \pi_* \frac{d v(t)}{dt} \right) E_l.
$$
Therefore, if we choose a system of coordinates $\oy$ on $K$ and
then a system of coordinates on $N K$ defined by
$$
  (\oy, \ov{\z}) \in \R^k \times \R^n \qquad \mapsto \qquad \ov{\z}^j E_j(\oy),
$$
we get that
$$
  \hat{g}_{\oa \ob}(\oy,\ov{\z}) = \ov{g}_{\oa \ob}(\oy) + \ov{\z}^i
  \ov{\z}^j \left\langle \n^N_{\pa_{\oa}} E_i, \n^N_{\pa_{\ob}} E_j
  \right\rangle_N =
  \ov{g}_{\ov{a} \ov{b}}(\oy) + \ov{\z}^i \ov{\z}^j \b^l_i \left(
  \pa_{\oa} \right) \b^l_j \left( \pa_{\ob} \right),
$$
and
$$
  \hat{g}_{\oa \ov{i}}(\oy,\ov{\z}) = \ov{\z}^j \b^i_j \left(
  \pa_{\ov{a}} \right); \qquad \qquad \hat{g}_{\ov{i} \ov{j}}(\oy,\ov{\z})
  = \d_{\ov{i} \ov{j}},
$$
where we have set $\pa_{\ov{i}} = \frac{\pa}{\pa \ov{\z}_i}$. We
notice also that the following co-area type formula holds, for any
smooth compactly supported function $f : N K  \to \R$
\begin{equation}\label{eq:coarea}
    \int_{N K} f dV_{\hat{g}} = \int_K \left( \int_{N_{\oy} K}
    f(\ov{\z}) d \ov{\z} \right) dV_{\ov{g}}(\oy).
\end{equation}
This follows immediately from the fact that $\det \hat{g} = \det
\ov{g}$, which in turn can be verified by expressing $\hat{g}$ as a
product of three matrices like
$$
  \left(
    \begin{array}{cc}
      Id & \z \b \\
      0 & Id \\
    \end{array}
  \right) \left(
            \begin{array}{cc}
              \ov{g} & 0 \\
              0 & Id \\
            \end{array}
          \right) \left(
                    \begin{array}{cc}
                      Id & 0 \\
                      \z \b & Id \\
                    \end{array}
                  \right),
$$
the first and the third having determinant equal to $1$.

Having defined the metric $\hat{g}$, we express  the Laplacian of a
function $u$ defined on $N K$ with respect to this metric. In Fermi
coordinates centered at some point $q \in K$, using
\eqref{eq:bfermiq}, \eqref{eq:bfermiq2} and \eqref{eq:ffcc}, it
turns out that (for $\oy = 0$)
\begin{equation}\label{eq:lapluhg}
    \D_{\hat{g}} u = \pa^2_{\oa \oa} u + \pa^2_{\ov{i} \ov{i}} u.
\end{equation}

\

Next we define the set $S_\e$ as
$$
  S_\e = \left\{ (v, \zn) \in N K_\e \times \R_+ \; : \; \left(|v|^2 + \zn^2
  \right)^{\frac 12} \leq \e^{-\g} \right\}, \qquad \quad \R_+ = \left\{
  \zn \; : \; \zn > 0 \right\},
$$
where $N K_\e$ stands for the normal bundle of $K_\e$ (in $\O_\e$).
We next endow $S_\e$ with a natural metric, inherited by $\hat{g}$
through a scaling. If $R_\e$ denotes the dilation $x \mapsto \e x$
in $\R^N$ (extended naturally to its subsets), we define a metric
$\tilde{g}_\e$ on $S_\e$ by
$$
  \tilde{g}_\e = \frac{1}{\e^2} [(R_\e)_* \hat{g}] \otimes  d \zn^2.
$$
In particular, choosing coordinates $(y,\z')$ on $N K_\e$ via the
scaling $(\oy,\ov{\z}) = \e (y,\z')$, one easily checks that the
components of $\tilde{g}_\e$ are given by
$$
  (\tilde{g}_\e)_{ab}(y,v) = (\ov{g})_{\oa \ob}(\e y) + \e^2
  v^i v^j \b^l_i \left( \pa_{\oa} \right) (\e y)
  \b^l_j \left( \pa_{\ob} \right) (\e y),
$$
$$
  (\tilde{g}_\e)_{ai}(y,v) = \e v^j \b^i_j \left( \pa_{\oa}
  \right)(\e y); \qquad \qquad (\tilde{g}_\e)_{ij}(\oy,v) = \d_{ij},
$$
and also
$$
  (\tilde{g}_\e)_{N N} \equiv 1; \qquad \qquad (\tilde{g}_\e)_{N \a} \equiv
  0.
$$
Therefore, if $u$ is a smooth function in $S_\e$, it follows that in
the above coordinates $(y, \z', \zn)$ (at $y = 0$)
\begin{equation}\label{eq:lapltg}
    \D_{\tilde{g}_\e} u = \pa^2_{a a} u + \pa^2_{i i} u +
  \pa^2_{\zn \zn} u.
\end{equation}

In the following, to emphasize a slow dependence of a function $u$
in the variables $y$, we will often write $u(y, \z) = u(\e y, \z)$
(where, we recall, $\z = (\z',\zn)$), identifying with an abuse of
notation the variable $y$ parameterizing $K_\e$ with $\oy$,
parameterizing $K$. In this case we have that (at the origin of the
Fermi coordinates)
\begin{equation}\label{eq:lapltge}
    \D_{\tilde{g}_\e} u = \e^2 \pa^2_{\oa \oa} u +
    \pa^2_{i i} u + \pa^2_{\zn \zn} u.
\end{equation}

For later purposes, we evaluate $\D_{\tilde{g}_\e}$ on functions
with a special structure. In particular, if we deal with a function
$u$ of the form $u (\oy,\z) = \phi(\oy) v(|\z|)$, we have that
\begin{equation}\label{eq:lapradu}
  \D_{\tilde{g}_\e} u = \e^2 (\D_K \phi (\oy)) v(|\z|) + \phi(\oy) \D_{\z}
  v,
\end{equation}
and if instead $u(\oy,\z) = v(|\z|) \psi^h(\oy) \frac{\z_h}{|\z|}$
for some smooth normal section $\psi = \psi^h E_h$, then we find
\begin{equation}\label{eq:laplinu}
  \D_{\tilde{g}_\e} u = \e^2 (\D_K^N \psi)^h (\oy) \frac{\z_h}{|\z|}
  v(|\z|) + \psi^h (\oy) \D_{\z} \left(v(|\z|) \frac{\z_h}{|\z|} \right).
\end{equation}

\

\no Now we introduce the function space $H_{S_\e}$ defined as the
family of functions in $H^1(S_\e)$ which vanish on $\{ |v|^2 + \zn^2
= \e^{-2\g} \}$, endowed with the scalar product
\begin{equation}\label{eq:sphse}
    (u,v)_{H_{S_\e}} = \int_{S_\e} \left( \n_{\tilde{g}_\e}
    u \cdot \n_{\tilde{g}_\e}v + u v \right) dV_{\tilde{g}_\e}.
\end{equation}
We consider next the operator $T_{S_\e} : H_{S_\e} \to H_{S_\e}$
defined by duality as
\begin{equation}\label{eq:TSe}
    (T_{S_\e} u, v)_{H_{S_\e}} = \int_{S_\e} \left( \n_{\tilde{g}_\e}
    u \cdot \n_{\tilde{g}_\e}v + u v
    - p w_0^{p-1}(|\z|) u v \right) dV_{\tilde{g}_\e},
\end{equation}
for arbitrary $u, v \in H_{S_\e}$. Our goal is to characterize some
of the eigenvalues of $T_{S_\e}$, with the corresponding
eigenfunctions.

For simplicity, if $u_{\a,\e}$, $v_{\a,\e,i}$, $\eta_{\a,\e}$ and
$\s_{\a,\e}$ are given by Proposition \ref{p:eicompe}, recalling our
notation from Subsection \ref{ss:op}, we also set
\begin{equation}\label{eq:ukvk}
  u_{j,\e} = u_{\e^2 \rho_j, \e}; \qquad v_{l,\e,i} = v_{\e^2
  \o_l, \e, i}; \qquad \qquad \eta_{j,\e} = \eta_{\e^2 \rho_j,\e};
\qquad \s_{l,\e} = \s_{\e^2 \o_l,\e}.
\end{equation}
We also assume that these functions are normalized so that
\begin{equation}\label{eq:normeigenf}
\left\{
  \begin{array}{ll}
    \|u_{j,\e}\|^2_{\e^2 \rho_j,\e} = \int_{B_{\g,\e}} \left( |\n u_{j,\e}|^2
    + (1 + \e^2 \rho_j) u_{j,\e}^2 \right) = 1;  &  \\
    \|v_{l,\e,i}\|^2_{\e^2 \o_l,\e} = \int_{B_{\g,\e}} \left( |\n v_{l,\e,i}|^2
    + (1 + \e^2 \o_l) v_{l,\e,i}^2 \right) = 1. &
  \end{array}
\right.
\end{equation}
After these preliminaries, we can state our result.

\begin{pro}\label{p:dec}
Let $\e_0, \e$ be as in Proposition \ref{p:eicompe}. Let $\l <
\frac{\t}{4}$ be an eigenvalue of $T_{S_\e}$. Then either $\l =
\eta_{j,\e}$ for some $j$, or $\l = \s_{l,\e}$ for some index $l$.
The corresponding eigenfunctions $u$ are of the form
\begin{equation}\label{eq:etse}
  u(y,\z) = \sum_{\{j \; : \; \eta_{j,\e} = \l\}}
  \a_j \phi_j(\e y) u_{j,\e}(\z) + \sum_{\{l \; : \; \s_{l,\e} = \l\}}
  \b_l \varphi^i_l(\e y) v_{l,\e,i}(\z),
\end{equation}
where $(y,\z)$ denote the above coordinates on $S_\e$, and where
$(\a_j)_j$, $(\b_l)_l$ are arbitrary constants. Viceversa, every
function of the form \eqref{eq:etse} is an eigenfunction of
$T_{S_\e}$ with eigenvalue $\l$. In particular the eigenvalues of
$T_{S_\e}$ which are smaller than $\frac{\t}{4}$ coincide with the
numbers $(\eta_{j,\e})_j$ or $(\s_{l,\e})_l$ which are smaller than
$\frac{\t}{4}$.
\end{pro}

\begin{pf}
The proof is based on separation of variables and the spectral
analysis of Proposition \ref{p:eicompe}. Integrating by parts, one
can check that the eigenfunction $u$ of $T_{S_\e}$ satisfies the
following equation
\begin{equation}\label{eq:eetse}
    \left\{
      \begin{array}{ll}
        (1 - \l) \left( - \D_{\tilde{g}_\e} u + u \right) -
         p w_0^{p-1}(\z) u = 0 & \hbox{ in } S_\e, \\
        \frac{\pa u}{\pa \zn} = 0 & \hbox{ on } \{\zn = 0\}.
      \end{array}
    \right.
\end{equation}
As before, we can extend $u$ evenly in $\zn$, to obtain a smooth
solution of the differential equation in \eqref{eq:eetse} in the set
$\{ (v,\zn) \in N K_\e \times \R \; : \; (|v|^2 + \zn^2)^{\frac 12}
\leq \e^{-\g} \}$. Hence, fixing $y \in K_\e$, we can use Fourier
decomposition in the angular variable of $\z$, and we can write
$$
  u(y,\z) = \sum_{l=0}^\infty u_l(y,|\z|) Y_{l,e}(\th),
$$
where $\th = \frac{\z}{|\z|} \in S^n$, and where $Y_{l,e}$ is the
$l$-th spherical harmonic function which is even in $\zn$. We now
decompose $u$ further in a convenient way as
\begin{equation}\label{eq:dechatu}
    u = \un{u}_0 + \un{u}_1 + \un{u}_2,
\end{equation}
where
$$
  \un{u}_0 = \frac{1}{\sqrt{|S^n|}}u_0(y,|\z|); \qquad \quad \un{u}_1
  = \sum_{l=1, \dots, n} u_l(y,|\z|) Y_{l,e}(\th); \qquad \quad
  \un{u}_2 = \sum_{l \geq n+1} u_l(y,|\z|) Y_{l,e}(\th).
$$
Integrating by parts, the last formula, together with
\eqref{eq:coarea}, \eqref{eq:lapradu} and \eqref{eq:laplinu} (recall
that $Y_{l,e}$ for $l = 1, \dots, n$ are linear combinations of
$\frac{\z_h}{|\z|}$ on $S^n$, $h = 1, \dots, n$) easily imply that
$(\un{u}_i, \un{u}_j)_{H_{S_\e}} = 0$ for $i \neq j$ and that
$(T_{S_\e} \un{u}_i, \un{u}_j)_{H_{S_\e}} = 0$ for $i \neq j$,
namely that $T_{S_\e}$ diagonalizes with respect to the above
decomposition \eqref{eq:dechatu}.

We begin by considering the action of $T_{S_\e}$ on $\un{u}_0$.
Using a Fourier decomposition of $\un{u}_0(y,|\z|)$ through the
eigenfunctions $(\phi_j)_j$ of the Laplace-Beltrami operator on
$K_\e$ we set
$$
  \un{u}_0(y,|\z|) = \sum_{j=0}^\infty \phi_j(\e y) \tilde{u}_j(|\z|).
$$
By \eqref{eq:lapradu} we get immediately that for any $j$
$$
  \D_{\tilde{g}_\e} (\phi_j(\e y) \tilde{u}_j(|\z|)) =
(\e^2 \D_{\ov{g}} + \D_{\z}) (\phi_j(\e y) \tilde{u}_j(|\z|)) =
(\D_{\z} - \e^2 \rho_j) \phi_j(\e y) \tilde{u}_j(|\z|).
$$
As a consequence we find that $\un{u}_0 \in H^1_\e$ satisfies the
following partial differential equation in $B_{\e,\g}$, with Neumann
boundary conditions on $\{\zn = 0\}$
$$
  - \D_{\tilde{g}_\e} \un{u}_0 + \un{u}_0 - p w_0^{p-1}(|\z|) \un{u}_0
  = \sum_{j=0}^\infty
  \phi_j(\e y) \left( - \D_{\z} \tilde{u}_j(|\z|) + (1 + \e^2 \rho_j)
  \tilde{u}_j(|\z|) - p w_0^{p-1}(|\z|) \tilde{u}_j(|\z|) \right).
$$
{From} this formula it follows that if $T_{S_\e} u = \l u$ for some
$\l$, then by the orthogonality to $\un{u}_1$, $\un{u}_2$ we have
also $T_{S_\e} \un{u}_0 = \l \un{u}_0$, and each of the components
$\tilde{u}_j$ (which are radial in $\z$) satisfies the eigenvalue
equation $T_{\e^2 \rho_j, \e} \tilde{u}_j = \l \tilde{u}_j$ in
$H_{\e^2 \rho_j,\e}$ with the same value of $\l$, where we are using
the notation of Subsection \ref{ss:spa}. Using the same terminology,
we can further decompose $\tilde{u}_j$ as
$$
  \tilde{u}_j(|\z|) = \a_j u_{j,\e} + \ov{u}_{j,\e} \qquad
  \hbox{ with } \a_j \in \R \hbox{ and with } (u_{j,\e},
  \ov{u}_{j,\e})_{\e^2 \rho_j,\e} = 0.
$$
{From} the spectral analysis carried out in the previous subsection
it follows that if $\l < \frac{\t}{4}$ (and $\e$ is sufficiently
small), then $\ov{u}_{j,\e} = 0$ for every $j$, and $\l =
\eta_{j,\e}$ for some set of indices $j$.

\

\no We now turn to the evaluation of $T_{S_\e}$ on $\un{u}_1$.
Similarly as before, expanding with respect to the eigenfunctions of
the normal Laplacian we can decompose $\un{u}_1$ in the following
way
$$
  \un{u}_1(y,\z) = \sum_{l \geq 0} \sum_{i=1}^n \tilde{v}_{l}(|\z|)
\varphi_{l,i}(\e y) \frac{\z_i}{|\z|},
$$
and from \eqref{eq:laplinu} we deduce that
\begin{eqnarray*}
  \D_{\tilde{g}_\e} \left( \sum_{i=1}^n \tilde{v}_{l}(|\z|)
\varphi_{l,i}(\e y) \frac{\z_i}{|\z|} \right) & = & \sum_{i=1}^n
(\e^2 \D_{K}^N \varphi_l)^i \frac{\z_i}{|\z|} v_l(|\z|) +
\sum_{i=1}^n \varphi_{l,i}(\e y)  \D_{\z} \left( \tilde{v}_{l}(|\z|)
\frac{\z_i}{|\z|}
\right) \\
  & = & (\D_{\z} - \e^2 \o_l) \left( \sum_{i=1}^n
\tilde{v}_{l}(|\z|) \varphi_{l,i}(\e y) \frac{\z_i}{|\z|} \right).
\end{eqnarray*}
As a consequence we find that also
\begin{eqnarray*}
   & & - \D_{\tilde{g}_\e} \un{u}_1 + \un{u}_1 - p w_0^{p-1}(|\z|)
  \un{u}_1 \\
   & = & \sum_{l \geq 0} \sum_{i=1}^n \varphi_{l,i}(\e y) \left[
- \D_{\z} \left( \tilde{v}_{l}(|\z|) \frac{\z_i}{|\z|} \right) + (1
+ \e^2 \o_l)
 \tilde{v}_{l}(|\z|) \frac{\z_i}{|\z|} - p w_0^{p-1}(|\z|) \tilde{v}_{l}(|\z|)
 \frac{\z_i}{|\z|} \right].
\end{eqnarray*}
Hence, by the spectral analysis of the previous subsection,
reasoning as for $\un{u}_0$ we deduce that if $\un{u}_1$ satisfies
$T_{S_e} \un{u}_1 = \l \un{u}_1$ with $\l < \frac \t 4$, then
$\tilde{v}_{l}(|\z|) \frac{\z_i}{|\z|} = v_{l,\e,i}$, and hence it
follows that $\l = \s_{l,\e}$ for some set of indices $l$.

Finally, we turn to $\un{u}_2$. Proceeding as for the definition of
the metric $\hat{g}$ (and using the same notation), we can introduce
a bilinear form $\mathfrak{g}$ (semi-positive definite) on $T \, N
K$ defined by
$$
  \mathfrak{g}(V,W) = \left \langle
  \frac{D^N v}{d t}|_{t=0}, \frac{D^N w}{d t}|_{t=0} \right\rangle_N.
$$
Using again a scaling in $\e$, we can also introduce the following
bilinear form on $S_\e$
$$
  \mathfrak{g}_\e = \frac{1}{\e^2} (R_\e)_* \mathfrak{g} \otimes  d \zn^2.
$$
The components of this form in the above coordinates $(y,\z)$ are
given by
$$
(\mathfrak{g}_\e)_{ab}(y,v) = \e^2 v^i v^j \b^l_i \left( \pa_{\oa}
\right) (\e y) \b^l_j \left( \pa_{\ob} \right) (\e y); \qquad \quad
(\mathfrak{g}_\e)_{ai}(y,v) = \e v^j \b^i_j \left( \pa_{\oa}
  \right)(\e y);
$$
$$
(\mathfrak{g}_\e)_{ij}(\oy,v) = \d_{ij}; \qquad (\mathfrak{g}_\e)_{N
N} \equiv 1; \qquad (\mathfrak{g}_\e)_{N \a} \equiv 0.
$$
We then define by duality the operator $\mathfrak{T}_\e$ through the
formula
$$
  ( \mathfrak{T }_\e u , u)_{H_{S_\e}} := \int_{S_\e}
  \left[ \mathfrak{g}_\e(\n_{\tilde{g}_\e}u, \n_{\tilde{g}_\e}u) +
  u^2 - p w_0^{p-1}(|\z|) u^2 \right] dV_{\tilde{g}_\e}.
$$
%Since $\mathfrak{g} \leq \un{g}$ on $N K$, it follows easily that
%\begin{equation}\label{eq:compforms}
%  ( \mathfrak{T }_\e u , u)_{H_{S_\e}} \leq ( T_{S_\e} u ,
%  u)_{H_{S_\e}} \qquad \qquad \hbox{ for every } u \in H_{S_\e}.
%\end{equation}
Moreover, computing the pointwise action of $\mathfrak{T}_\e$
integrating by parts, reasoning as for the derivation of
\eqref{eq:lapltge}, and using \eqref{eq:coarea}, one finds that
\begin{equation}\label{eq:split2}
    ( \mathfrak{T }_\e u , u)_{H_{S_\e}} =
\int_{K_\e}\left[\int_{S_{y,\e}}\left( -u \,\D_\z u +u^2 - p
w_0(|\z|)^{p-1} \right) d \z \right]dV_{\ov{g}_\e}(y), \qquad \qquad
u \in H_{S_\e},
\end{equation}
where we have set $\ov{g}_\e = \frac{1}{\e^2} (R_\e)_* \ov{g}$ and
$S_{y,\e} = \left\{ (v, \zn) \in N_y K_\e \times \R_+ \; : \;
\left(|v|^2 + \zn^2 \right)^{\frac 12} \leq \e^{-\g} \right\}$.

Hence, using \eqref{eq:coarea} (with the scaled metric
$\tilde{g}_\e$), \eqref{eq:varo} with $u = \un{u}_2$ and
\eqref{eq:split2} we find
$$
  p \int_{S_\e} w_0^{p-1} \un{u}_2^2 dV_{\tilde{g}_\e} = p \int_{K_\e}
  \left( \int_{S_{y,\e}} w_0^{p-1} \un{u}_2^2 \right) dV_{\ov{g}_\e}(y) \leq
  \left( 1 - \frac \t
    2 \right) \int_{K_\e} \left[\int_{S_{y,\e}}\left( -\un{u}_2 \,\D_\z \un{u}_2
+\un{u}_2^2\right)\right]dV_{\ov{g}_\e}(y).
%\\  & = & \left( 1 - \frac \t
%    2 \right) ( \mathfrak{T }_\e \un{u}_2, \un{u}_2)_{H_{S_\e}} \leq \left( 1 - \frac \t
%    2 \right)
%  ( T_{S_\e} \un{u}_2, \un{u}_2)_{H_{S_\e}}.
$$
Since $\t < 1$ (being an eigenvalue of $\ov{J}''(w_0) \leq
Id_{H^1(\R^{n+1}_+)}$), we deduce that
\begin{eqnarray*}
  (T_{S_\e} u, u)_{H_{S_\e}} & = & ( \mathfrak{T }_\e u , u)_{H_{S_\e}} +
  \int_{S_\e}
  \left[ (\hat{g}_\e - \mathfrak{g}_\e ) (\n_{\tilde{g}_\e}u, \n_{\tilde{g}_\e}u) +
  u^2 \right] dV_{\tilde{g}_\e} \\
  & \geq & \frac{\t}{2} \int_{S_\e}
  \left[ \mathfrak{g}_\e(\n_{\tilde{g}_\e}u, \n_{\tilde{g}_\e}u) +
  u^2 \right] dV_{\tilde{g}_\e} + \int_{S_\e}
  \left[ (\hat{g}_\e - \mathfrak{g}_\e ) (\n_{\tilde{g}_\e}u, \n_{\tilde{g}_\e}u) +
  u^2 \right] dV_{\tilde{g}_\e} \\ & \geq & \frac{\t}{2}
  \|u\|^2_{H_{S_\e}}.
\end{eqnarray*}
If follows that there are no eigenvectors of the form $\un{u}_2$
corresponding to eigenvalues smaller than $\frac \t2$. This
concludes the proof.
\end{pf}

\begin{rem}\label{r:split2}
For later purposes, it is convenient to consider a splitting of the
functions in $H_{S_\e}$ which is slightly different from the one in
\eqref{eq:dechatu}. If $\un{u}_0$, $\un{u}_1$ and $\un{u}_2$ are as
above, with
$$
  \un{u}_0 = \sum_{j \geq 0} \phi_j(\e y)
  \tilde{u}_j(|\z|); \qquad \qquad \un{u}_1 = \sum_{l \geq 0}
  \sum_{i=1}^n \tilde{v}_l(|\z|) \varphi_{l,i}(\e y)
  \frac{\z_i}{|\z|},
$$
for some real sequences $(\a_j)_j$, $(\b_l)_l$, we can write
$$
  \tilde{u}_j(|\z|) = \a_j u_{j,\e}(|\z|) + \ov{u}_{j,\e}(|\z|), \qquad
\hbox{ with } \qquad ( u_{j,\e}, \ov{u}_{j,\e})_{\e^2 \rho_j, \e} =
0;
$$
$$
  \tilde{v}_l(|\z|) \frac{\z_i}{|\z|} = \b_l v_{l,\e,i}(\z) +
  \ov{v}_{l,\e}(|\z|) \frac{\z_i}{|\z|} := \b_l v_{l,\e,i}(\z) +
  \ov{v}_{l,\e,i}(\z), \qquad
\hbox{ with } \qquad ( v_{l,\e,i} , \ov{v}_{l,\e,i})_{\e^2 \o_l, \e}
= 0.
$$
Now we set $u = \mathfrak{u}_0 + \mathfrak{u}_1 + \mathfrak{u}_2$,
where
$$
  \mathfrak{u}_0 = \sum_{j=0}^\infty \a_j u_{j,\e}(|\z|) \phi_j(\e y); \qquad
\qquad
  \mathfrak{u}_1 = \sum_{l=0}^\infty \b_l v_{l,\e,i}(\z) \varphi_l^i(\e y);
$$
$$
  \mathfrak{u}_2 = \sum_{j=0}^\infty \ov{u}_{j,\e}(|\z|) \phi_j(\e y) +
  \sum_{l=0}^\infty \ov{v}_{l,\e,i}(\z) \varphi_l^i(\e y) + \un{u}_2.
$$
Then by \eqref{eq:normeigenf} one can check that $(\mathfrak{u}_i,
\mathfrak{u}_j)_{H_{S_\e}} = 0$ for $i \neq j$, and that
\begin{equation}\label{eq:rsplnu2}
    \|u\|_{H_{S_\e}}^2 = \|\mathfrak{u}_0\|_{H_{S_\e}}^2 + \|\mathfrak{u}_1\|_{H_{S_\e}}^2 +
\|\mathfrak{u}_2\|_{H_{S_\e}}^2 =  \frac{1}{\e^k} \sum_{j=0}^\infty
\a_j^2 + \frac{1}{\e^k} \sum_{l=0}^\infty \b_l^2 +
\|\mathfrak{u}_2\|_{H_{S_\e}}^2;
\end{equation}
\begin{equation}\label{eq:rspltseu}
    (T_{S_\e} u, u)_{H_{S_\e}} = \sum_{j=0}^\infty \eta_{j,\e} \a_j^2 +
    \sum_{l=0}^\infty
\s_{l,\e} \b_l^2 + (T_{S_\e} \mathfrak{u}_2,
\mathfrak{u}_2)_{H_{S_\e}}; \qquad \quad (T_{S_\e} \mathfrak{u}_2,
\mathfrak{u}_2)_{H_{S_\e}} \geq C \|\mathfrak{u}_2\|_{H_{S_\e}}^2,
\end{equation}
for some fixed positive constant $C$.
\end{rem}

\

\noindent From the last proposition we deduce the following
corollary, regarding the Morse index of $T_{S_\e}$.

\begin{cor}\label{c:eigenv}
Let $\g \in (0,1)$, and let $T_{S_\e} : H_{S_\e} \to H_{S_\e}$ be
defined as before. Then, as $\e$ tends to zero, the Morse index of
$T_{S_\e}$ satisfies the estimate
$$
  M.I.(T_{S_\e}) \simeq \left( \frac{\ov{\a}}{C_k} \right)^{\frac
  k2}Vol(K) \e^{-k},
$$
where $\ov{\a}$ is the unique real number for which $\eta_{\ov{\a}}
= 0$ (see Remark \ref{r:a0}).
\end{cor}

\begin{pf} From Proposition \ref{p:dec} we have that the Morse index
of $T_{S_\e}$ is equal to the number of negative $\eta_{j,\e}$'s. By
the estimate in \eqref{eq:muaae}, this number is asymptotic to the
number of $j$'s for which $\eta_{\e^2 \rho_j}$ is negative.
Therefore it is sufficient to count the number of eigenvalues
$\rho_j$ for which $\e^2 \rho_j$ is less than $\ov{\a}$. By the
Weyl's asymptotic formula, see \cite{LSY}, we have that $\rho_j
\simeq C_k \left(\frac{j}{Vol(K)}\right)^{\frac k2}$ so the
conclusion follows immediately.
\end{pf}

\section{Accurate analysis of the linearized operator}\label{s:real}

In this section we first compare $J''_\e(u_{I,\e})$ to the model
operator introduced in the previous one. A naive direct comparison
will give errors of order $\e$, see Lemma \ref{l:comp} and Corollary
\ref{c:comp}, but sometimes we will need estimates of order $\e^2$.
Therefore we will expand at a higher order the eigenvalues (of the
linearized operator at $u_{I,\e}$) close to zero with the
corresponding eigenfunctions, to get sufficient control on the
errors. Finally, using these expansions, we will define a suitable
decomposition of the functional space for which the linearized
operator is almost {\em diagonal}.

\subsection{Comparison of $J''_\e(u_{I,\e})$ and $T_{S_\e}$}\label{ss:comp}

We define first a bijection $\tilde{\Upsilon}_\e$ from $S_\e$ into a
neighborhood of $K_\e$ in $\O_\e$ in the following way. Given the
section $\Phi = \Phi_0 + \e \Phi_1 + \dots + \e^{I-2} \Phi_{I-2}$ in
$N K$ constructed in Section \ref{s:as}, for any $(v,\zn) \in S_\e$,
$v \in N_y K_\e$, $\zn \in \R_+$, we set
$$
  \tilde{\Upsilon}_\e(v,\zn) = \exp_y^{\pa \O_\e} (v + \Phi(\e y))
  + \zn \nu \left( \exp_y^{\pa \O_\e} (v + \Phi(\e y)) \right).
$$
Then we define the set $\Sig_\e \subseteq \O_\e$ to be
$$
  \Sig_\e = \tilde{\Upsilon}_\e(S_\e),
$$
endowed with the standard Euclidean metric induced from $\R^N$. For
$u \in H_{S_\e}$, we define the function $\tilde{u} : \Sig_\e \to
\R$ by
$$
  \tilde{u}(z) = u \left( \tilde{\Upsilon}_\e^{-1}(z) \right),
  \qquad \qquad z \in \Sig_\e,
$$
and letting $\L_\e$ to be the map $u \mapsto \tilde{u}$, we define
$$
  H_{\Sig_\e} = \L_\e(H_{S_\e}).
$$
$H_{\Sig_\e}$ has a natural structure of Hilbert (Sobolev) space
inherited by $H^1(\O_\e)$, and we denote by $( \cdot,
\cdot)_{H_{\Sig_\e}}$, $\| \cdot \|_{H_{\Sig_\e}}$ the corresponding
scalar product and norm. More precisely, we can identify the space
$H_{\Sig_\e}$ with the family of functions in $H^1(\O_\e)$ which
vanish identically in $\O_\e \setminus \Sig_\e$.

We introduce next the operator $T_{\Sig_\e} : H_{\Sig_\e} \to
H_{\Sig_\e}$ defined as the restriction to $H_{\Sig_\e}$ of
$J''_\e(u_{I,\e})$ which, using the duality in $H_{\Sig_\e}$, has
the following expression
\begin{equation}\label{eq:TSige}
    (T_{\Sig_\e} u, v)_{H_{\Sig_\e}} = \int_{\Sig_\e} \left( \n u \cdot
  \n v + u v \right) - p \int_{\Sig_\e} u_{I,\e}^{p-1} u v =
  ( u, v)_{H_{\Sig_\e}} - p \int_{\Sig_\e} u_{I,\e}^{p-1} u
  v.
\end{equation}

Fixing these notations and definitions, following the arguments at
the beginning of Section 4 in \cite{mal2} one can easily prove the
following result.

\begin{lem}\label{l:comp}
Identifying the functions in $H_{S_\e}$ with the corresponding ones
in $H_{\Sig_\e}$ via the map $\L_\e$, for $\e$ sufficiently small
one has
$$
  ( u, v)_{H_{\Sig_\e}}=( u, v)_{H_{S_\e}}+O(\e^{1-\g})
  \|u\|_{H_{S_\e}} \|v\|_{H_{S_\e}};
$$
$$
  ( T_{\Sig_\e} u, v)_{H_{\Sig_\e}}=( T_{S_\e} u, v)_{H_{S_\e}}
  +O(\e^{1-\g}) \|u\|_{H_{S_\e}} \|v\|_{H_{S_\e}}.
$$
with error $O(\e^{1-\g})$ independent of $u, v \in H_{\Sig_\e}$.
\end{lem}

We introduced the operator $T_{\Sig_\e}$ because it represents an
accurate model for $J''_\e(u_{I,\e})$. In fact, since most of the
functions we consider have an exponential decay away from $K_\e$, it
is reasonable to expect that the spectrum of $J''_\e(u_{I,\e})$ will
be affected only by negligible quantities if we work in
$H_{\Sig_\e}$ instead of $H^1(\O_\e)$. More precisely, one has the
following result (we recall the definition of $\t$ from the previous
section).

\begin{lem}\label{l:expcomp}
There exists a fixed constant $C$, depending on $\O$, $K$ and $p$
such that the eigenvalues of $J''_\e(u_{I,\e})$ and $T_{\Sig_\e}$
satisfy
$$
  \left| \l_j(J''_\e(u_{I,\e})) - \l_j(T_{\Sig_\e}) \right| \leq C
  e^{-\frac{1}{C \e^{-\g}}}, \qquad \hbox{ provided }
  \l_j(J''_\e(u_{I,\e})) \leq \frac \t 2.
$$
Here we are indexing the eigenvalues in non-decreasing order,
counted with multiplicity.
\end{lem}

We omit the proof of this result because it is very similar in
spirit to that of Lemma 5.5 in \cite{malm}. This is based on the
fact that the number of the eigenvalues of $T_{S_\e}$ which are less
or equal than $\frac 34 \t$ is bounded by $\e^{-D}$ for some $D > 0$
(see Proposition \ref{p:dec} and the Weyl's asymptotic formulas in
Subsection \ref{ss:op}), together with the exponential decay of the
eigenfunctions of $J''_\e(u_{I,\e})$, which can be shown as in
\cite{malm}, Lemma 5.1.

As a consequence of Lemmas \ref{l:comp} and \ref{l:expcomp} we
obtain the following result.

\begin{cor}\label{c:comp}
In the above notation, for $\e$ small one has that
\begin{equation}\label{eq:estccomp}
    \left| \l_j(J''_\e(u_{I,\e})) - \l_j(T_{S_\e}) \right| \leq C
  \e^{1-\g}, \qquad \hbox{ provided } \l_j(J''_\e(u_{I,\e})) \leq
  \frac \t 2.
\end{equation}
\end{cor}

Using Proposition \ref{p:dec} and Corollary \ref{c:comp}, it is
possible to obtain some qualitative information about the spectrum
of the linearized operator $J''_\e(u_{I,\e})$. However, this kind of
estimate is not sufficiently precise by the following
considerations. First of all, since the eigenvalues of $T_{S_\e}$
can approach zero at a rate $\min\{ \e^2, \e^k \}$, the estimate
\eqref{eq:estccomp} need to be improved if we want to guarantee the
invertibility of $J''_\e(u_{I,\e})$. Furthermore, it would be
natural to expect that the Jacobi operator (and its invertibility)
plays some role in the expansion of the eigenvalues, and this is not
apparent here.

On the other hand, Lemma \ref{l:expcomp} gives an accurate estimate
on the eigenvalues of $J''_\e(u_{I,\e})$ in terms of those of
$T_{\Sig_\e}$, so it will be convenient to analyze $T_{\Sig_\e}$
directly.

\subsection{Approximate eigenfunctions of $T_{\Sig_\e}$}\label{ss:ae}

In this subsection we construct approximate eigenfunctions to the
linearized operator at the approximate solutions $u_{I,\e}$. By the
reasons explained at the end of the previous subsection, we need a
refined expansion of the small eigenvalues of $T_{\Sig_\e}$, and in
particular here we want to understand how the $\s_{l,\e}$'s change
when we pass from $T_{S_\e}$ to $T_{\Sig_\e}$.

It is sufficient here to take $I = 2$, because the terms of order
higher than $\e^2$ do not affect the expansions below. As for the
construction of the approximate solutions $u_{I,\e}$, we proceed by
expanding the eigenvalue equation formally in powers of $\e$. By the
construction of $u_{2,\e}$, {\em formally} the following equation
holds
$$
  - \D_{{g}_\e} u_{2,\e} + u_{2,\e} - u_{2,\e}^p = O(\e^3).
$$
Using Fermi coordinates as in Section \ref{s:as} and differentiating
with respect to $\z_h$, we get
\begin{equation}\label{eq:canc}
    - \pa_{h} (\D_{{g}_\e} u_{2,\e}) + \pa_{h} u_{2,\e} - p
  u_{2,\e}^{p-1} \pa_{h} u_{2,\e} = O(\e^3).
\end{equation}
{From} the general expression of the Laplace-Beltrami operator, see
formula \eqref{eq:laplcoord}, we can easily see that
\begin{eqnarray}\label{eq:dhLu} \nonumber
  \pa_h (\D_{{g}_\e} u) & = & \D_{{g}_\e}(\pa_h u) + \pa_h {g}_\e^{AB} \pa_{AB} u +
  \pa_h (\pa_A {g}_\e^{AB}) \pa_B u \\
   & + & \frac 12 {g}_\e^{AB} \pa^2_{hA} \left( \log (\det { g}_\e) \right)
   \pa_B u + \frac 12 \pa_A \left( \log (\det { g}_\e) \right) (\pa_h
   { g}_\e^{AB}) \pa_B u.
\end{eqnarray}
Let us now consider the second term on the right-hand side of
\eqref{eq:dhLu}: dividing the indices this is equivalent to
$$
  \pa_h { g}_\e^{ij} \pa^2_{ij} u + 2 \pa_h { g}_\e^{ib} \pa^2_{ib} u + \pa_h { g}_\e^{ab}
  \pa_{ab} u + 2 \pa_h { g}_\e^{AN} \pa_A \pa_{\zn} u.
$$
{From} Lemma \ref{l:expgeuz}, and using the fact that we get an $\e$
factor each time we differentiate $u$ with respect to $y_a, y_b,
\dots$, we find that
$$
  \pa_h { g}_\e^{AB} \pa^2_{AB} u = - \frac 23 \e^2 R_{ihtj} \z_t \pa^2_{ij} u
  + O(\e^3).
$$
Similarly we get
$$
  \pa_h \pa_A { g}_\e^{AB} \pa_B u = \frac 13 \e^2 R_{hiij} \pa_j u + O(\e^3);
$$
\begin{eqnarray*}
  \frac 12 { g}_\e^{AB} \pa^2_{hA} \left( \log (\det { g}_\e) \right) \pa_B u & = &
  \e^2 \left( \frac 13 R_{illh} + R_{iaah} - \G^b_a(E_i)
  \G^a_b(E_h) \right) \pa_i u \\
  & + & 2 H_{ab} \G^b_a(E_h) \pa_{\zn} u + O(\e^3),
\end{eqnarray*}
and
$$
  \frac 12 \pa_A \left( \log (\det { g}_\e) \right) (\pa_h
   { g}_\e^{AB}) \pa_B u = O(\e^3).
$$
Putting together all these terms we deduce that
\begin{equation}\label{eq:dhDu}
  \pa_h (\D_{{ g}_\e} u) = \D_{{ g}_\e} (\pa_h u) - \frac 23 \e^2
  R_{ihtj} \z_t \pa_{ij} u + \e^2 \left( \frac 23 R_{illh}
  + R_{iaah} - \G^b_a(E_i) \G^a_b(E_h) \right) \pa_i u + O(\e^3).
\end{equation}

\no To construct the approximate eigenfunctions $\mathfrak{v}_\e$
and the approximate eigenvalues $\mu$, we make an {\em ansatz} of
the type
$$
\mathfrak{v}_\e = \left(\psi^h(\oy) \pa_h
u_{2,\e}(\oy,\z'+\Phi(\oy),\zn) + \e^2 z_2(\oy,\z)\right) + O(\e^3);
\qquad \quad \mu = \e^2 \ov{\mu} + O(\e^3),
$$
where the normal section $\psi = (\psi^h)_h$, the function $z_2$ and
the real number $\ov{\mu}$ have to be determined.

\ms

\no We notice that the eigenvalue equation $J''_\e(u_{2,\e}) v = \l
v$ in $H^1(\O_\e)$, with an integration by parts becomes
$$
- \D_{{ g}_\e}v
  + v - p \left( u_{2,\e}
  \right)^{p-1} v
  =    \l \left( - \D_{{ g}_\e} v + v\right),
$$
see also the derivation of \eqref{eq:eva}.

For $v=\mathfrak{v}_\e$ and $\l=\mu$, we have the following
expansion
\begin{eqnarray*}
  & - & \D_{{ g}_\e} \left(\psi^h(\oy) \pa_h u_{2,\e} + \e^2 z_2(\oy,\z)\right)
  \\ &
  + & \psi^h(\oy) \pa_h u_{2,\e} + \e^2 z_2(\oy,\z) - p \left( u_{2,\e}
  \right)^{p-1} \left( \psi^h(\oy) \pa_h
u_{2,\e} + \e^2 z_2(\oy,\z) \right) \\[3mm]
  & = & \e^2 \ov \mu \left[ - \D_{{ g}_\e} \left( \psi^h(\oy) \pa_h u_{2,\e}
  + \e^2 z_2(\oy,\z)
  \right) + \left( \psi^h(\oy) \pa_h u_{2,\e} + \e^2 z_2(\oy,\z)
  \right) \right] \\[3mm] & = & \e^2 \ov \mu \left[ \psi^h(\oy)
  \left( - \D_{{ g}_\e} \pa_h w_0 + \pa_h w_0 \right)
  \right] + O(\e^3) \\[3mm] & = & \e^2 \ov \mu p \psi^h(\oy) w_0^{p-1}
  \pa_h w_0 + O(\e^3).
\end{eqnarray*}
{From} \eqref{eq:dhDu} we can expand the Laplacian in the last
formula as
\begin{eqnarray*}
  - \D_{{ g}_\e} \left(\psi^h(\oy) \pa_h u_{2,\e} \right) & = &
  - \e^2
  \pa^2_{\oy_a \oy_a} \psi^h \pa_h w_0 - 2 \e^2 \pa_a
  \psi^h \pa^2_{jh} w_0 \pa_{\oy_a} \Phi^j_0 - \psi^h \D_{{ g}_\e}
  (\pa_h u_{2,\e}) \\ & + & 4 \e^2\zn H_{aj}
   \pa_{\oy_a} \psi^h \pa^2_{jh} w_0 + O(\e^3) \\
   & = & - \e^2
  \pa^2_{\oy_a \oy_a} \psi^h \pa_h w_0 - 2 \e^2 \pa_a
  \psi^h \pa^2_{jh} w_0 \pa_{\oy_a} \Phi^j_0
   - \psi^h \pa_h (\D_{{ g}_\e} u_{2,\e}) \\ & + & 4 \e^2\zn H_{aj}
   \pa_{\oy_a} \psi^h \pa^2_{jh} w_0
   + \frac 23 \e^2 \psi^h R_{ihtj} \z_t \pa_{ij} w_0 \\
   & - & \e^2 \psi^h \left( \frac 23 R_{illh}
  + R_{iaah} - \G^b_a(E_i) \G^a_b(E_h) \right) \pa_i w_0 + O(\e^3).
\end{eqnarray*}
Using \eqref{eq:canc} jointly with the last equality, and recalling
our previous notation (from Section \ref{s:as})
$$
\mathcal{L}_{\Phi} u = - \D u + u - p w_0^{p-1}(\z'+\Phi(\e y),\zn),
$$
we obtain the following condition on $z_2$
\begin{eqnarray}\label{eq:eqz2} \nonumber
  \mathcal{L}_\Phi z_2 & = & \pa^2_{\oy_a \oy_a} \psi^h \pa_h w_0 + 2 \pa_{\oy_a}
  \psi^h \pa^2_{j h} w_0 \pa_{\oy_a} \Phi^j_0 - \frac 23 \psi^h R_{ihtj}
  \z_t \pa^2_{ij} w_0 \\
  & + & \psi^h \left( \frac 23 R_{illh} + R_{iaah} - \G^b_a(E_i)
  \G^a_b(E_h) \right) \pa_i w_0 + p \ov \mu \psi^h w_0^{p-1}
  \pa_h w_0 \\ & - & 2 H_{ab} \G^b_a(E_h) \pa_\zn w_0 - 4 \zn
  H_{aj} \pa_{\oy_a} \psi^m \pa^2_{jm} w_0 + O(\e). \nonumber
\end{eqnarray}
In order to get solvability of this equation (in $z_2$), we need to
impose that the right-hand side is orthogonal to the kernel of
$\mathcal{L}_\Phi$ namely that, multiplying it by $\pa_s w_0$ and
integrating in $\z$, $s = 1, \dots, n$, we must get zero. If we do
this, reasoning as at the end of Subsection \ref{ss:e2}, we obtain
the following condition on $\psi$
$$
  C_0 \mathfrak{J} \psi = C_1 \ov \mu \psi, \qquad \quad \hbox{ where
  } \quad \qquad C_1 = p \int_{\R^{n+1}_+} w_0^{p-1} (\pa_1 w_0)^2 d
  \z,
$$
and where $C_0$ is given in \eqref{eq:C0}. With the choices
$$
  \ov \mu = \frac{C_0}{C_1} \mu_l; \qquad \qquad  \psi = \psi_l,
$$
where $\mu_l$ is an eigenvalues of $\mathfrak{J}$ with eigenfunction
$\psi_l$, the right-hand side of \eqref{eq:eqz2} is perpendicular to
the kernel of $\mathcal{L}_\Phi$, and we get solvability in $z_2$.
Using the eigenvalue equation for $\psi_l$, \eqref{eq:eqz2} can be
simplified as
\begin{eqnarray*}
  \mathcal{L}_\Phi z_2 & = & \mu_l \psi_l^h \pa_h w_0 \left( p \frac{C_0}{C_1} w_0^{p-1}
  - 1 \right) + 2 \pa_{\oy_a} \psi^h_l \left( \pa_{\oy_a} \Phi^j_0 -
  2 \zn H_{aj} \right) \pa^2_{jh} w_0 \\
  & + & \frac 23 \psi_l^h \bigg( R_{ijjh} \pa_i w_0 - R_{ihtj} \z_t
  \pa^2_{ij} w_0 - 3 H_{ab} \G^b_a(E_h) \pa_\zn w_0 \bigg).
\end{eqnarray*}

Next, we set
$$
  g_0^h(\oy,\z) = \mathcal{L}_\Phi^{-1} \left[ \pa_h w_0 \left( p \frac{C_0}{C_1}
  w_0^{p-1} - 1 \right) \right]; \qquad
  g_1^h(\oy,\z) = 2 \mathcal{L}_\Phi^{-1} \left[ \left( \pa_{\oy_a} \Phi^j_0 -
  2 \zn H_{aj} \right) \pa^2_{jh} w_0 \right];
$$
$$
  g_2^h(\oy,\z) = \frac 23 \mathcal{L}_\Phi^{-1} \left[ \left( R_{illh} \pa_i w_0
  - R_{ihtj} \z_t \pa^2_{ij} w_0 - 3 H_{ab} \G^b_a(E_h) \pa_\zn w_0
  \right) \right] + \pa_h w_2(\oy,\z'+\Phi(\oy),\zn),
$$
and $$  g_3^h(\oy,\z) = \pa_h w_1(\oy,\z'+\Phi(\oy),\zn). $$

\no We notice that, by the definitions of $C_0, C_1$, the
computations in Subsection \ref{sss:steps} and by oddness, the
arguments of $\mathcal{L}_\Phi^{-1}$ in the definitions of $g_0^h$,
$g_1^h$ and $g_2^h$ are all perpendicular to the kernel of
$\mathcal{L}_\Phi$, and therefore $g_0, g_1$ and $g_2$ are well
defined.

\medskip

\no Finally, with this notation,  we define the approximate
eigenfunction $\Psi_l$ as $\mathfrak{v}_\e$ times a suitable cut-off
function of $\z$, namely
\begin{equation}\label{eq:Psil}
  \Psi_l(\oy,\z) = \chi_\e(|\z|) \bigg\{ \psi_l^h(\oy) \left[
  \pa_h w_0 + \e g_3 ^h(\oy,\z) + \e^2 g_2^h(\oy,\z) \right]
  + \e^2\mu_l \psi_l^h(\oy) g_0^h(\oy,\z) +
  \e^2\pa_{\oy_a} \psi^h_l(\oy) g_1^h(\oy,\z) \bigg\},
\end{equation}
where $\chi_\e$ is as in \eqref{eq:chi}, and, as usual, $\oy = \e
y$.

A more accurate analysis, which we omit, shows that the above error
terms not only are of order $\e^3$, but they decay exponentially to
zero as $|\z|$ tends to infinity. Moreover, as we already remarked,
in the above estimates one can replace $u_{2,\e}$ with $u_{I,\e}$.
Precisely, one can prove the following result.

\begin{lem}\label{l:decae}
If $\Psi_l$ is given in \eqref{eq:Psil}, then there exist a
polynomial $P(\z)$ and a sequence of positive constants $(C_l)_l$,
depending on $\O$, $K$, $p$ and $I$ such that
$$
  \left| - \D_{g_\e} \Psi_l + \Psi_l - p u^{p-1}_{I,\e} \Psi_l - \e^2
  \frac{C_0}{C_1} \mu_l (- \D_{g_\e} \Psi_l + \Psi_l ) \right| \leq
  C_l \e^3 P(\z) e^{- |\z|}.
$$
\end{lem}

\subsection{A splitting of the functional space}\label{ss:dec}

In the previous subsection we expanded in $\e$ some of the
eigenvalues of $T_{\Sig_\e}$, precisely those which are the
counterparts of the $\s_{l,\e}$'s for $T_{S_\e}$. Actually,
$T_{S_\e}$ possesses another type of resonant eigenvalues, namely
the $\eta_{j,\e}$'s for suitable values of $j$, which in principle
could approach zero even faster. One of the differences between
these two families of eigenvalues is that the eigenfunctions
corresponding to the resonant $\s_{l,\e}$'s oscillate {\em slowly}
along $\pa \O_\e$, and this allowed us to perform the above
expansion. On the contrary, the eigenfunctions related to the
$\eta_{j,\e}$'s possess only high Fourier modes, and therefore such
an expansion is not possible anymore. Nevertheless, we can deal with
the counterparts of these eigenvalues applying Kato's theorem, which
on the other hand requires to characterize the corresponding
eigenfunctions up to some extent.

The purpose of the present subsection is to identify appropriate
subspaces of $H_{\Sig_\e}$ with respect to which $T_{\Sig_\e}$ is
approximately in block form. Recalling the definitions in
Proposition \ref{p:eicompe}, in \eqref{eq:ukvk} and in
\eqref{eq:Psil} (and also our convention about the range of an
integer index), for $\d \in \left(0,k\right)$, $\ov{C} \in (0,1)$,
we define the following subspaces
\begin{equation}\label{eq:h1}
  H_1  = \hbox{span} \left\{ \phi_i (\e y) u_{i,\e}(\z), i = 0, \dots,
\infty \right\};
\end{equation}
\begin{equation}\label{eq:h2}
    \hat{H}_2 = \hbox{span} \left\{ \Psi_l, l = 0, \dots, \e^{-\d}
\right\}; \quad \tilde{H}_2 = \hbox{span} \left\{ \psi_j^m (\e y)
\hat v_{j,\e}(|\z|)\,\frac{\z_m}{|\z|}, j = \e^{-\d} + 1, \dots,
\ov{C} \e^{-k} \right\};
\end{equation}
\begin{equation}\label{eq:h23}
    H_2 = \hat{H}_2 \oplus \tilde{H}_2; \qquad \qquad H_3 = \left( H_1
\oplus H_2 \right)^\perp,
\end{equation}
where $X^\perp$ denotes the orthogonal complement to the subspace
$X$ with respect to the scalar product in $H_{\Sig_\e}$. We have the
following result, which is the counterpart of Proposition 4.2 in
\cite{mal2}. The proof follows the same arguments, but for the
reader's convenience we prefer to give details since the notation
and the estimates are affected by the different dimensions and
codimensions we are dealing with.

\begin{pro}\label{p:ud}
There exists a small value of the constant $\ov{C} > 0$ in
\eqref{eq:h2}, depending on $\O$, $K$ and $p$, such that the
following property holds. For $\e$ sufficiently small and choosing
$\d \in \left( \frac{k}{2}, k \right)$ in \eqref{eq:h2}, every
function $u \in H_{\Sig_\e}$ decomposes uniquely as
$$
u = u_1 + u_2 + u_3, \qquad \hbox{ with } \qquad u_1 \in H_1, u_2
\in H_2, u_3 \in H_3.
$$
Moreover there exists a positive constant $C$, also depending on
$\O$, $K$ and $p$ such that
$$
(T_{\Sig_\e} u_3, u_3) \geq \frac{1}{C \ov{C}^{\frac{2}{k}}}
\|u_3\|^2_{H_{\Sig_\e}}.
$$
\end{pro}

\noindent The proof requires some  preliminary Lemmas. Before
stating them, we recall our convention about the symbol $\sum_c^d$,
for two positive real values $c$ and $d$.

\begin{lem}\label{l:ntu2}
Let $\tilde{u}_2 = \sum_{j=\e^{-\d}+1}^{\ov{C} \e^{-k}} \b_j
\psi_j^m(\e y) \hat v_{j,\e}(|\z|)\frac{\z_m}{|\z|} \in
\tilde{H}_2$. Then
\begin{equation}\label{eq:claim1}
    \|\tilde{u}_2\|_{H_{\Sig_\e}}^2 = (1 + O(\e^{1-\g}))
    \frac{1}{\e^k} \sum_{j=\e^{-\d}}^{\ov{C} \e^{-k}} \b_j^2.
\end{equation}
\end{lem}

\begin{pf}
By Lemma \ref{l:comp}, it is sufficient to estimate
$\|\tilde{u}_2\|_{H_{S_\e}}^2$. We notice that by \eqref{eq:jacobi}
there holds
$$
  - \D^N_K\psi_j = \mathfrak{J} \psi_j + (\mathfrak{B} -
  \mathfrak{R}) \psi_j = \mu_l \psi_j + ((\mathfrak{B} -
  \mathfrak{R}) \psi)_j.
$$
Integrating by parts, using \eqref{eq:laplinu} and the last formula
one finds that $\|\tilde{u}_2\|_{H_{S_\e}}^2$ becomes
\begin{eqnarray}\label{eq:aaaa} \nonumber
 &- &\int_{S_\e}\sum\limits_{j,l=\e^{-\d}+1}^{\ov{C} \e^{-k}}\D_{\tilde
 g_\e}\bigg(   \sum\limits_{m=1}^n\b_j
\psi_j^m (\e y) \hat v_{j,\e}(|\z|)\,\frac{\z_m}{|\z|}
\bigg)\cdot\bigg( \sum\limits_{h=1}^n\b_l \psi_l^h (\e y) \hat
v_{l,\e}(|\z|)\,\frac{\z_h}{|\z|}   \bigg)
\\
&+&\int_{S_\e}\sum\limits_{j,l=\e^{-\d}+1}^{\ov{C} \e^{-k}}\bigg(
\sum\limits_{m=1}^n\b_j \psi_j^m (\e y) \hat
v_{j,\e}(|\z|)\,\frac{\z_m}{|\z|}   \bigg)\cdot\bigg(
\sum\limits_{h=1}^n\b_l \psi_l^h (\e y) \hat
v_{j,\e}(|\z|)\,\frac{\z_h}{|\z|}   \bigg) = A_1 + A_2,
\end{eqnarray}
where
$$
  A_1 = \int_{S_\e}\sum\limits_{j,l=\e^{-\d}+1}^{\ov{C} \e^{-k}}
\left[ \left( - \D_{\z} + (1+\e^2\mu_j) \right)\bigg(
\sum\limits_{m=1}^n\b_j \psi_j^m (\e y) \hat
v_{j,\e}(|\z|)\,\frac{\z_m}{|\z|} \bigg) \right] \cdot \bigg(
\sum\limits_{h=1}^n\b_l \psi_l^h (\e y) \hat
v_{j,\e}(|\z|)\,\frac{\z_h}{|\z|}   \bigg);
$$
$$
 A_2 = \e^2 \int_{S_\e}\sum\limits_{j,l=\e^{-\d}+1}^{\ov{C}
\e^{-k}}\bigg( \sum\limits_{m=1}^n\b_j \left((\mathfrak{B} -
  \mathfrak{R})\psi_j  \right)^m (\e y) \hat v_{j,\e}(|\z|)\,\frac{\z_m}{|\z|}   \bigg)\cdot\bigg(
\sum\limits_{h=1}^n\b_l \psi_l^h (\e y) \hat
v_{j,\e}(|\z|)\,\frac{\z_h}{|\z|}   \bigg).
$$
Looking at $A_1$, the integral over any fiber $N_y K_\e$ is non zero
if and only if $m = h$ (and by symmetry, when computing the integral
we can assume both the indices to be $1$). Then, from
\eqref{eq:coarea} and from the orthogonality among different
$\psi_l$'s (which now are scaled in $\e$), recalling that
$\hat{v}_{j,\e}(|\z|) \frac{\z_m}{|\z|} = v_{j,\e,m}$, $A_1$ becomes
$$
  \frac{1}{\e^k} \sum_{j=\e^{-\d}+1}^{\ov{C} \e^{-k}}
  \b_j^2 \|v_{j,\e,1}\|^2_{\e^2 \eta_j, \e} = \frac{1}{\e^k}
  \sum_{j=\e^{-\d}+1}^{\ov{C} \e^{-k}} \b_j^2 \left[ \int_{\R^{n+1}_+} \left( |\n v_{j,\e,1}|^2 +
  (1 + \e^2 \mu_j) v_{j,\e,1}^2 \right) \right].
$$
Recalling the normalization \eqref{eq:normeigenf} and the fact that
$\eta_j = \o_j + O(1)$ (independently of $j$), see Subsection
\ref{ss:op}, we obtain that
\begin{equation}\label{eq:estA1}
    A_1 = \frac{1}{\e^k} \sum_{j=\e^{-\d}+1}^{\ov{C} \e^{-k}} (1 + O(\e^2)) \b_j^2.
\end{equation}
We turn now to the estimate of $A_2$. By the orthogonality of the
$\psi_l$'s, using again \eqref{eq:coarea} and \eqref{eq:normeigenf}
one finds
$$
  \int_{S_\e} \tilde{u}_2^2 dV_{\tilde{g}_\e} = \frac{1}{\e^k} \sum_{j=\e^{-\d}+1}^{\ov{C} \e^{-k}}
  \b_j^2 \|v_{j,\e,1}\|^2_{L^2(\R^{n+1}_+)} \leq \frac{1}{\e^k} \sum_{j=\e^{-\d}+1}^{\ov{C} \e^{-k}}
  \b_j^2.
$$
Working in a local system of coordinates $(y,z)$ as in Subsection
\ref{ss:model}, it is also convenient to write $\tilde{u}_2$ as
$$
  \tilde{u}_2(y,\z) = \sum_{m=1}^n f_m(y,|\z|) \z_m, \qquad \hbox{ where }
  \qquad f_m(y,|\z|) = \sum_{j=\e^{-\d}+1}^{\ov{C} \e^{-k}} \b_j \psi_j^m(\e y)
  \frac{\hat{v}_{j,\e}(|\z|)}{|\z|}.
$$
If $\mathcal{U}$ is a neighborhood of some point $q$ in $K$, where
the coordinates $y$ are defined, letting $\mathcal{U}_\e = \frac 1
\e \mathcal{U}$, one has
$$
  \int_{N \mathcal{U}_\e} \tilde{u}_2^2 dV_{\tilde{g}_\e} =
  \sum_{m=1}^n
  \int_{\mathcal{U}_\e} \left( \int_{\R^{n+1}_+} f_m^2 (y,|\z|) \z_1^2
  d \z \right) dV_{\ov{g}_\e}(y),
$$
so it follows that
\begin{equation}\label{eq:comparison}
    \sum_{m=1}^m
  \int_{\mathcal{U}_\e} \left( \int_{\R^{n+1}_+} f_m^2 (y,|\z|) \z_1^2
  d \z \right) dV_{\ov{g}_\e}(y) \leq \int_{S_\e} \tilde{u}_2^2
  dV_{\tilde{g}_\e} \leq \frac{1}{\e^k} \sum_{j=\e^{-\d}+1}^{\ov{C} \e^{-k}} \b_j^2.
\end{equation}

Now, we can write
$$
  A_2 = \e^2 \int_{S_\e} \tilde{\mathfrak{u}}_2 \tilde{u}_2
  dV_{\tilde{g}_\e}, \qquad \hbox{ where } \qquad \tilde{\mathfrak{u}}_2 =
  \sum_{j=\e^{-\d}+1}^{\ov{C} \e^{-k}}
  \b_j \left((\mathfrak{B} - \mathfrak{R})\psi_j  \right)^m (\e y)
   \hat v_{j,\e}(|\z|)\frac{\z_m}{|\z|}.
$$
As for $\tilde{u}_2$, we can write $\tilde{\mathfrak{u}}_2 =
\sum_{m=1}^n \mathfrak{f}_m(y,|\z|) \z_m$, where $\mathfrak{f}_m =
\sum_{j=\e^{-\d}+1}^{\ov{C} \e^{-k}} (\mathfrak{B} -
\mathfrak{R})_{mj} f_j(y,|\z|)$, and compute
$$
  \int_{N \mathcal{U}_\e} \tilde{\mathfrak{u}}_2^2 dV_{\tilde{g}_\e} =
  \sum_{m=1}^n
  \int_{\mathcal{U}_\e} \left( \int_{\R^{n+1}_+} \mathfrak{f}_m^2 (y,|\z|) \z_1^2
  d \z \right) dV_{\ov{g}_\e}(y).
$$
In conclusion, from the H\"older inequality, from
\eqref{eq:comparison}, covering $K_\e$ with finitely-many
$\mathcal{U}_\e$'s we derive
\begin{equation}\label{eq:estA2}
    |A_2| \leq \e^2 \left( \int_{S_\e} \tilde{u}_2^2 dV_{\tilde{g}_\e}
    \right)^{\frac 12} \left( \int_{S_\e} \tilde{\mathfrak{u}}_2^2
    dV_{\tilde{g}_\e} \right)^{\frac 12} \leq C \e^2 \frac{1}{\e^k} \|
    \mathfrak{B} - \mathfrak{R} \|_{L^\infty} \sum_{j=\e^{-\d}+1}^{\ov{C} \e^{-k}} \b_j^2.
\end{equation}
Then the conclusion follows from \eqref{eq:estA1} and
\eqref{eq:estA2}.
\end{pf}

In order to estimate the norm $\|\hat{u}_2\|_{H_{\Sig_\e}}$, it is
convenient to introduce an abstract result.

\begin{lem}\label{l:abstr}
For $j \in \{ 0, \dots, \e^{-\d}\}$, and for a sequence $(\b_j)_j$,
let us consider a function $u : S_\e \to \R$ of the form
$$
  u(y,\z) = \sum_{j=0}^{\e^{-\d}} \sum_{m=1}^n \b_j
  (L_{d,\oy} \psi_j^m)(\oy) g_m(\z),
$$
where $\oy = \e y$, where $L_{d,\oy}$ is a linear differential
operator of order $d$ with smooth coefficients in $\oy$, and where
the functions $g_m(\z)$ are also smooth and have an exponential
decay at infinity.

Then there exists a positive constant $C$, independent of $\e$, $\d$
and $(\b_j)_j$ such that
$$
  \|u\|_{L^2(S_\e)}^2 \leq C \frac{1}{\e^k} \sum_{j=0}^{\e^{-\d}} \left( 1 + \e^{2d}
  |\mu_j|^d \right) \b_j^2.
$$
\end{lem}

\begin{pf}
The proof is similar in spirit to that of Lemma \ref{l:ntu2}, but
here we take advantage of the fact that the {\em profile} $g_m(\z)$
is independent of the index $j$ (this lemma applies in particular to
each of the summands in the definition of $\Psi_l$, see
\eqref{eq:Psil}).

Using local coordinates, \eqref{eq:coarea} and the exponential decay
of the $g_m$'s, after integration in $\z$ we find
$$
  \|u\|_{L^2(S_\e)}^2 = \sum_{j,l=0}^{\e^{-\d}} \sum_{m,h=1}^n
  \b_j \b_l c_{mh} \int_{\mathcal{U}_\e}
  (L_{d,\oy} \psi_j^m)(\oy) (L_{d,\oy} \psi_l^h)(\oy) dV_{\ov{g}},
$$
for some bounded coefficients $(c_{mh})$. As for \eqref{eq:estA2}
then we find $\|u\|_{L^2(S_\e)} \leq C \|\psi\|_{H^d(K_\e,N K_\e)}$
and the last quantity, with a change of variables and by
\eqref{eq:l2l2}, can be estimated with $\frac{C}{\e^k}
\sum_{j=0}^{\e^{-\d}} (1 + \e^{2d} |\mu_j|^d) \b_j^2$. This
concludes the proof.
\end{pf}

\begin{lem}\label{l:nu2}
Let $u_2 = \hat{u}_2 + \tilde{u}_2 = \sum\limits_{j=0}^{\e^{-\d}}
\b_j \Psi_j(\e y, \z) + \sum\limits_{j=\e^{-\d}+1}^{\ov{C} \e^{-k}}
\b_j \psi_j^m(\e y) \hat v_{j,\e}(|\z|)\frac{\z_m}{|\z|} \in H_2$.
Then, choosing $\d \in \left( \frac{k}{2}, k \right)$ in
\eqref{eq:h2}, one has
\begin{equation}\label{eq:claim3}
    \|u_2\|_{H_{\Sig_\e}}^2 = \frac{1}{\e^k} (1 + O(\e^{1-\g}+\e^{2-\frac{2\d}{k}}))
    \left[ \sum_{j=0}^{\e^{-\d}} \b_j^2 \|\partial_1 w_0\|_{H^1(\R^{n+1}_+)}
    + \sum_{j=\e^{-\d}+1}^{\ov{C} \e^{-k}} \b_j^2 \right].
\end{equation}
\end{lem}

\begin{pf}
We first claim that the following formula holds
\begin{equation}\label{eq:claim2}
    \|\hat{u}_2\|_{H_{S_\e}}^2 = \frac{ 1}{\e^k} \sum_{j=0}^{\e^{-\d}}
    \b_j^2 \left(1 + O(\e^{2-2\g}+\e^{2-\frac{2\d}{k}})\right) \|\partial_1
    w_0\|_{H^1(\R^{n+1}_+)}^2.
\end{equation}

\noindent  {\bf Proof of \eqref{eq:claim2}.} We write
$$
\hat{u}_2 = \hat{u}_{2,1}+ \hat{u}_{2,2} := \sum_{j=0}^{\e^{-\d}}
\b_j \psi^m_j (\e y) \partial_m w_0 (\z) \chi_\e(|\z|) +
\sum_{j=0}^{\e^{-\d}} \b_j \ov{\Psi}_j(\e y, \z).
$$
where $\ov{\Psi}_j$ is the term of order $\e$ (and higher) in
$\Psi_j$. Reasoning as in the proof of Lemma \ref{l:ntu2} we get
\begin{eqnarray}\label{eq:esthu21} \nonumber
    \|\hat{u}_{2,1}\|_{H_{S_\e}}^2 & = & \frac{1}{\e^k} \sum_{j=0}^{\e^{-\d}} \b_j^2 (1 +
\e^2 \mu_j + O(\e^2)) \|\partial_m w_0 \chi_\e\|_{H^1(\R^{n+1}_+)}^2 \\
& = & \frac{1}{\e^k} \sum_{j=0}^{\e^{-\d}} \b_j^2 \left(1 +
O(\e^{2-\frac{2\d}{k}})\right) \|\partial_1
w_0\|_{H^1(\R^{n+1}_+)}^2,
\end{eqnarray}
where the last equality follows from the Weyl's asymptotic formula
\eqref{eq:weyl4}.

On the other hand, using Lemma \ref{l:abstr}, the Weyl's formula and
some computations, one also finds
\begin{eqnarray*}
  \e^k \|\hat{u}_{2,2}\|_{H_{S_\e}}^2 & \leq & C \e^2 \sum_{j=0}^{\e^{-\d}} \b_j^2
  \left( 1 + \e^2 |\mu_j|  \right) + C \e^4 \sum_{j=0}^{\e^{-\d}} \b_j^2
  \mu_j^2 \left( 1 + \e^2 |\mu_j| \right) \nonumber \\
  & + & C \e^4 \sum_{j=0}^{\e^{-\d}} \b_j^2 \left( 1 + |\mu_j| + \e^2
  |\mu_j|^3 \right) \leq C \left(\e^2+\e^{4-\frac{4\d}{k}} +
  \e^{6-\frac{6\d}{k}} \right) \sum_{j=0}^{\e^{-\d}} \b_j^2.
\end{eqnarray*}
By our choice of $\d$, the last formula reads
\begin{equation}\label{eq:nhu22}
  \|\hat{u}_{2,2}\|_{H_{S_\e}}^2 \leq \frac{C}{\e^k}
  \e^{4-\frac{4\d}{k}} \sum_{j=0}^{\e^{-\d}} \b_j^2.
\end{equation}
Finally, from \eqref{eq:esthu21} and \eqref{eq:nhu22} we also obtain
$$
(\hat{u}_{2,1},\hat{u}_{2,2})_{H_{S_\e}}\le \frac{C}{\e^k}
\sum_{j=0}^{\e^{-\d}} \b_j^2 \left(\e +
O(\e^{2-\frac{2\d}{k}})\right),
$$
which concludes the proof of \eqref{eq:claim2}.

\

\noindent  {\bf Proof of \eqref{eq:claim3}.} We write again
$\hat{u}_2 = \hat{u}_{2,1} + \hat{u}_{2,2}$. Then, by the
orthogonality relations among the $\psi_j$'s, reasoning as in the
proof of Lemma \ref{l:ntu2}, we get that $(\tilde{u}_2,
\hat{u}_{2,1})_{H_{S_\e}}$ becomes
$$
  \e^2 \sum\limits_{j=\e^{-\d}+1}^{\ov{C}
\e^{-k}} \sum\limits_{l=0}^{\e^{-\d}} \int_{S_\e} \bigg(
\sum\limits_{m=1}^n\b_j \left((\mathfrak{B} -
  \mathfrak{R})\psi_j  \right)^m (\e y) \hat v_{j,\e}(|\z|)\,\frac{\z_m}{|\z|}
    \bigg)\cdot\bigg( \chi_\e(|\z|)
\sum\limits_{h=1}^n\b_l \psi_l^h (\e y) \pa_h w_0   \bigg).
$$
As above, with some computations we find
\begin{eqnarray*}
    (\tilde{u}_2, \hat{u}_{2,1})_{H_{S_\e}} = O(\e^2)
    \|\tilde{u}_2\|_{H_{S_\e}} \|\hat{u}_{2,1}\|_{H_{S_\e}} =
    O(\e^2) \frac{1}{\e^k} \sum_{j=0}^{\ov{C} \e^{-k}} \b_j^2.
\end{eqnarray*}
From Lemma \ref{l:ntu2} and \eqref{eq:nhu22} we also find
$$
  (\tilde{u}_2, \hat{u}_{2,1})_{H_{S_\e}} \leq C \frac{1}{\e^k}
  \left( \sum_{j=0}^{\ov{C} \e^{-k}} (1 + O(\e^{1-\g}) \b_j^2)
  \right)^{\frac{1}{2}}
  \e^{2 - \frac{2\d}{k}} (\sum_{j=0}^{\ov{C} \e^{-k}} \b_j^2)^{\frac 12}.
$$
The result follows from the last two formulas.
\end{pf}

\begin{rem}\label{r:dec2}
From the proof of \eqref{eq:claim3} it also follows that every
function $u_2 \in H_2$ can be written uniquely as $u_2 = \hat{u}_2 +
\tilde{u}_2$, with $\hat{u}_2 \in \hat{H}_2$ and $\tilde{u}_2 \in
\tilde{H}_2$.
\end{rem}

\

\begin{pfn} {\sc of Proposition \ref{p:ud}.} In order to prove the
uniqueness of the decomposition it is sufficient to show that, for
$\e$ small
\begin{equation}\label{eq:this}
    (u_1, u_2)_{H_{\Sig_\e}} = o_\e(1) \|u_1\|_{H_{\Sig_\e}}
\|u_2\|_{H_{\Sig_\e}}, \qquad u_1 \in H_1, u_2 \in H_2,
\end{equation}
where $o_\e(1) \to 0$ as $\e \to 0$. Indeed, by Lemma \ref{l:comp}
we have
$$
(u_1,u_2)_{H_{\Sig_\e}} = (u_1, u_2)_{H_{S_\e}} + O(\e^{1-\g})
\|u_1\|_{H_{\Sig_\e}} \|u_2\|_{H_{\Sig_\e}},
$$
and since the functions $\partial_h w_0$, $g_0^h$, $g_3^h$ and
$v_{l, \e,i}$ are odd in $\z'$ (and so also $\tilde{u}_2$ and
$\hat{u}_{2,1}$), we get
$$
(u_1, u_2)_{H_{S_\e}} = (u_1, \hat{u}_{2,2})_{H_{S_\e}},
$$
where we have used the notation in the proof of Lemma \ref{l:nu2}.
Hence from the last three formulas, \eqref{eq:nhu22} and form
\eqref{eq:claim3} we deduce
\begin{equation}\label{eq:qperp}
    (u_1, u_2)_{H_{\Sig_\e}} \leq C (\e^{1-\g} + \e^{2-2\frac{\d}{k}})
    \|u_1\|_{H_{\Sig_\e}} \|u_2\|_{H_{\Sig_\e}},
\end{equation}
which implies \eqref{eq:this}, since $\d \in (\frac{k}{2},k)$.

\ms

\no To prove the second statement, it is sufficient to show that
\begin{equation}\label{eq:smu2u3}
    (u_3,v)_{H_{S_\e}} \leq \frac 12 \|u_3\|_{H_{S_\e}}
    \|v\|_{H_{S_\e}}; \qquad \hbox{ as } \e \to 0,
\end{equation}
for all $u_3 \in H_3$ and for all the functions $v$ of the form
$$
v = \sum_{l=0}^{\frac{1}{2} \ov{C} \e^{-k}} \tilde{\b}_l \varphi_l^m
(\e y) v_{l,\e,m}(\z).
$$
In fact, if we write $u_3 = \mathfrak{u}_{3,0} + \mathfrak{u}_{3,1}
+ \mathfrak{u}_{3,2}$ as in Remark \ref{r:split2} (with an obvious
change of notation),
$$
  \mathfrak{u}_{3,0} = \sum_{j=0}^\infty \a_j u_{j,\e}(|\z|) \phi_j(\e y); \qquad
\qquad
  \mathfrak{u}_{3,1} = \sum_{l=0}^\infty \b_l v_{l,\e,i}(\z) \varphi_l^i(\e
y),
$$
from \eqref{eq:rsplnu2} we find
\begin{equation}\label{eq:ffff}
    \|u_3\|^2_{S_\e} = \frac{1}{\e^k} \sum_{l=0}^\infty \left( \a_l^2 + \b_l^2 \right)  +
  \|\mathfrak{u}_{3,2}\|_{H_{S_\e}}^2.
\end{equation}
{From} \eqref{eq:rsplnu2}, from Lemma \ref{l:comp} and from the fact
that $u_3$ is perpendicular in $H_{\Sig_\e}$ to $\mathfrak{u}_{3,0}
\in H_1$, we deduce
$$
\frac{1}{ \e^k} \sum_{l=0}^\infty \a_l^2 = (\mathfrak{u}_{3,0},
\mathfrak{u}_{3,0})_{H_{S_\e}} = (\mathfrak{u}_{3,0},
u_3)_{H_{S_\e}} = O(\e^{1-\g})
\|u_3\|_{H_{S_\e}}\|\mathfrak{u}_{3,0}\|_{H_{S_\e}} \leq C \e^{1-\g}
\|u_3\|_{H_{S_\e}}^2.
$$
Moreover from \eqref{eq:smu2u3}, choosing $v =
\sum_{l=0}^{\frac{1}{2} \ov{C} \e^{-k}} \b_l \varphi_l^m (\e y)
v_{l,\e,m}(\z)$, and using \eqref{eq:ffff} we get
$$
 \frac{1}{ \e^k} \sum_{l \leq
\frac{1}{2} \ov{C} \e^{-k}} \b_l^2 = (u_3,v)_{H_{S_\e}} \leq \frac
12 \|u_3\|_{S_\e}^2.
$$
The last two formulas and \eqref{eq:ffff} then imply
\begin{equation}\label{eq:label}
    \|u_3\|_{H_{S_\e}}^2 \leq C \left( \sum_{l >
\frac{1}{2} \ov{C} \e^{-k}} \b_l^2 +
\|\mathfrak{u}_{3,2}\|_{H_{S_\e}}^2 \right),
\end{equation}
for some fixed constant $C$.

On the other hand, by \eqref{eq:rspltseu} we also have
\begin{eqnarray*}
  (T_{S_\e} u_3, u_3)_{S_\e} \geq \frac {1}{ \e^k} \sum_{l > \frac{1}{2} \ov{C}
  \e^{-k} + 1} \s_{l,\e} \b_l^2 + \frac {1}{C}
\|\mathfrak{u}_{3,2}\|_{H_{S_\e}}^2.
\end{eqnarray*}
Using the fact that $\s_{i,\e}\sim \s_{\e^2 w_i,\e}\sim \e^2
i^{\frac{2}{k}}$ by Proposition \ref{p:eicompe},  from
\eqref{eq:label} and the last formula it follows that
$$
  (T_{S_\e} u_3, u_3)_{S_\e} \geq \frac{1}{\e^k} \frac{1}{C \ov{C}^{\frac{2}{k}}}
  \sum_{i > \frac{1}{2} \ov{C} \e^{-k} + 1}  \b_l^2 + \frac {1}{C}
\|\mathfrak{u}_{3,2}\|_{H_{S_\e}}^2
  \geq \frac{1}{C \ov{C}^{\frac{2}{k}}}
  \|u_3\|_{H_{S_\e}}^2.
$$
This yields our conclusion, hence we are reduced to prove
\eqref{eq:smu2u3}.

\

\no {\bf Proof of \eqref{eq:smu2u3}.} By the form of $v$ and by
\eqref{eq:rsplnu2}, we have
\begin{equation}\label{eq:nnv}
    \|v\|^2_{H_{S_\e}} = \frac{1}{\e^k} \sum_{l=0}^{\frac{1}{2}
   \ov C \e^{-k}} \tilde{\b}_l^2.
\end{equation}
Using the $L^2$ basis $(\psi_l)_l$ of eigenfunctions of
$\mathfrak{J}$, we define the function $\varphi$ and the
coefficients $\{\b_l\}_{l=1,\dots, \infty}$ as
$$
\varphi (\oy) = \sum_{l=0}^{\frac 12 \ov{C} \e^{-k}} \tilde{\b}_l
\varphi_l(\oy) = \sum_{l=0}^\infty \b_l \psi_l(\oy) :=
\sum_{l=0}^\infty \b_l \psi_l^h(\oy) E_h(\oy),
$$
so we have
\begin{equation}\label{eq:nnnvp}
    \|\varphi\|_{L^2(K;NK)}^2 = \sum_{l=0}^{\frac 12 \ov{C} \e^{-k}}
\tilde{\b}_l^2 = \sum_{l=0}^\infty \b_l^2.
\end{equation}
Using these new coefficients $\b_j$, we set (see \eqref{eq:ukvk})
$$
  \tilde{v} (y,\z) = \ov{C}_0 \sum_{j=0}^{\e^{-\d}} \b_j \Psi_j(\e y, \z)
  + \sum_{j=\e^{-\d}+1}^{\ov{C} \e^{-k}} \b_j \psi_j^h(\e y)
  \hat v_{j,\e}(|\z|)\frac{\z_h}{|\z|} \in H_2.
$$
where $\ov{C}_0$ is given in Remark \ref{r:a0}. Hence we can write
\begin{eqnarray*}
    v - \tilde{v} = A_1 + A_2 + A_3 + A_4 + A_5,
\end{eqnarray*}
with
$$
A_1 = \sum_{l=0}^{\frac 12 \ov{C} \e^{-k}} \tilde{\b}_l
\varphi_l^m(\e y) \left[ v_{l,\e,m}(\z) - v_{0,\e,m}(\z) \right];
\qquad A_2 = \sum_{l=\ov{C} \e^{-k}+1}^\infty \b_l \psi_l^h (\e y)
v_{0,\e,h}(\z);
$$
$$
A_3 = - \ov{C}_0 \sum_{j=0}^{\e^{-\d}} \b_j \ov{\Psi}_j(\e y, \z);
\qquad \qquad A_4 = \sum_{l=\e^{-\d}+1}^{\ov{C} \e^{-k}} \b_l
\psi_l^h(\e y) \left( v_{0,\e,h} - v_{l,\e,h} \right);
$$
$$
A_5 = \sum_{l=0}^{\e^{-\d}} \b_l \psi_l^h \left( v_{0,\e,h} -
\ov{C}_0 \chi_\e(|\z|) \pa_h w_0 \right),
$$
and where $\ov{\Psi}_j$ is defined in the proof of Lemma
\ref{l:nu2}. Since $u_3$ is orthogonal to $H_2$, we get $(u_2,
\tilde{v})_{H_{\Sigma_\e}} = 0$, and so
\begin{equation}\label{eq:this2}
    (u_3, v)_{H_{\Sig_\e}} =
(u_3, A_1)_{H_{\Sig_\e}} + (u_3, A_2)_{H_{\Sig_\e}} + (u_3,
A_3)_{H_{\Sig_\e}} + (u_3, A_4)_{H_{\Sig_\e}} + (u_3,
A_5)_{H_{\Sig_\e}}.
\end{equation}
We prove now that $\|A_i\|_{H_{S_\e}}$ is small for every $i = 1,
\dots, 5$. From \eqref{eq:coarea}, the proof of Proposition
\ref{p:dec}, Proposition \ref{p:eicompe} and \eqref{eq:nnv} there
holds
$$
\|A_1\|^2_{H_{S_\e}} = \frac{1}{\e^k} \sum_{l=0}^{\frac 12 \ov{C}
\e^{-k}} \tilde{\b}_l^2 \|v_{l,\e,1} - v_{0,\e,1} \|_{l,\e}^2 \leq C
\ov{C}^2 (1 + \ov{C}^2) \|v\|^2_{H_{S_\e}} < \frac{1}{16}
\|v\|^2_{H_{S_\e}},
$$
provided $\ov{C}$ is sufficiently small.

To estimate $A_2$ we can use Lemma \ref{l:abstr} and some
computations to find
\begin{equation}\label{eq:A2}
    \|A_2\|^2_{H_{S_\e}} \leq C \frac{1}{\e^k} \sum_{l=\ov{C}
    \e^{-k} + 1}^\infty \b_l^2 (1 + \e^2 |\mu_l|).
\end{equation}
We now set $\tilde{\varphi} = \sum_{l=\ov{C} \e^{-k}+1}^\infty \b_l
\psi_l$. Since $\mathfrak{J} = - \D_K^N + O(1)$, for any integer $m$
one finds
\begin{eqnarray*}
    (\mathfrak{J}^m \tilde{\varphi}, \tilde{\varphi})_{L^2(K;NK)} &
\leq & (\mathfrak{J}^m {\varphi}, {\varphi})_{L^2(K)} \\ & \leq &
((- \D_K^N)^m {\varphi}, {\varphi})_{L^2(K;NK)} + C_m \left[ ((-
\D_K^N)^{m-1} {\varphi}, {\varphi})_{L^2(K;NK)} + ({\varphi},
{\varphi})_{L^2(K;NK)} \right].
\end{eqnarray*}
Since $\varphi = \sum_{l=0}^{\frac 12 \ov{C} \e^{-k}} \tilde{\b}_l
\varphi_l$, from \eqref{eq:nnnvp} we deduce that
\begin{eqnarray}\label{eq:aaa1} \nonumber
    (\mathfrak{J}^m \tilde{\varphi}, \tilde{\varphi})_{L^2(K;NK)}
& \leq & \left( \frac{\ov{C}}{2} \right)^{\frac{2m}{k}} \e^{-2m}
\|\varphi\|_{L^2(K;NK)}^2 + O(\e^{-2(m-1)})
\|\varphi\|_{L^2(K;NK)}^2 \\ & \leq & \left[ \left( \frac{\ov{C}}{2}
\right)^{\frac{2m}{k}} \e^{-2m} + O(\e^{-2(m-1)}) \right] \left(
\sum_{l=0}^{\frac 12 \ov{C} \e^{-k}} \tilde{\b}_l^2 \right).
\end{eqnarray}
On the other hand, since in the basis $(\psi_l)_l$, the function
$\tilde{\varphi}$ has non zero components only when $l \geq \ov{C}
\e^{-k}$, by the Weyl's asymptotic formula we have also that
\begin{equation}\label{eq:aaa2}
  (\mathfrak{J}^m \tilde{\varphi}, \tilde{\varphi})_{L^2(K;NK)} \geq
\left\{
  \begin{array}{ll}
    \sum_{l=\ov{C} \e^{-k} + 1}^\infty \mu_l^m \b_l^2; & \\ & \\
    C \ov{C}^{\frac{2m}{k}} \e^{-2m} \sum_{l=\ov{C} \e^{-k} + 1}^\infty
\b_l^2. &
  \end{array}
\right.
\end{equation}
Using \eqref{eq:aaa1} and the first inequality in \eqref{eq:aaa2}
with $m = 1$ we get
$$
  \e^2 \sum_{l=\ov{C} \e^{-k} + 1}^\infty \mu_l \b_l^2 \leq \left(
C \ov{C}^{\frac{2}{k}} + o_\e(1) \right) \sum_{l=0}^{\frac 12 \ov{C}
\e^{-k}} \tilde{\b}_l^2.
$$
Moreover, using  \eqref{eq:aaa1} and the second inequality in
\eqref{eq:aaa2} with $m$ arbitrary one also finds
$$
  \sum_{l=\ov{C} \e^{-k} + 1}^\infty \b_l^2 \leq \left(
 \left( \frac 12 \right)^{\frac{2m}{k}} + o_\e(1) \right)
\sum_{l=0}^{\frac 12 \ov{C} \e^{-k}} \tilde{\b}_l^2.
$$

Using \eqref{eq:nnv}, \eqref{eq:A2} and the last two inequalities
(for the second one we take $m$ large enough), for sufficiently
small $\ov{C}$ we find $\|A_2\|_{H_{S_\e}} < \frac{1}{16}
\|v\|_{H_{S_\e}}$.

Now we estimate  $\|A_3\|_{H_{S_\e}}$. Reasoning as for
\eqref{eq:nhu22}, from \eqref{eq:nnv} and \eqref{eq:nnnvp} we get
$$
\|A_3\|_{H_{S_\e}}^2 \leq C \frac{ 1}{ \e^k} \e^{4-4\frac{\d}{k}}
\sum_0^{\e^{-\d}} \b_j^2 \leq C
\e^{4-4\frac{\d}{k}}\|v\|_{H_{S_\e}}^2.
$$
Next, similarly to the estimate of $A_1$, for small $\ov{C}$ we find
$$
\|A_4\|^2_{H_{S_\e}} \leq \frac{1}{\e^k} C
\sum_{l=\e^{-\d}+1}^{\ov{C} \e^{-k}} \b_l^2 \left\|\hat v_{0,\e,1} -
v_{l,\e,1} \right\|_{l,\e}^2 \leq C \ov{C}^2 (1 + \ov{C}^2)
\|v\|^2_{H_{S_\e}} < \frac{1}{16} \|v\|^2_{H_{S_\e}}.
$$
Finally, from Proposition \ref{p:eicompe} and reasoning as for
$A_2$, we obtain also
$$
\|A_5\|^2_{H_{S_\e}} \leq \frac {1}{ \e^k} C e^{-C^{-1} \e^{-\g}}
\sum_{l=0}^{\e^{-\d}} \b_l^2(1+\e^2\o_l)\,\ov C
\e^{-k}l^{-\frac{2l}{k}}\leq C\e^{-k} e^{-C^{-1} \e^{-\g}}
\|v\|_{H_{S_\e}}^2.
$$
Taking \eqref{eq:this2} into account, this concludes the proof of
\eqref{eq:smu2u3}, provided we choose $\ov{C}$ and $\e$ sufficiently
small.
\end{pfn}

\section{Diagonalization of $T_{\Sig_\e}$ and applications}\label{s:applic}

In this section we study how the operator $T_{\Sig_\e}$ behaves with
respect to the above splitting of $H_{\Sigma_\e}$ in the three
subspaces $H_1, H_2$ and $H_3$. We prove that its form is almost
diagonal and we apply this analysis to study its invertibility for
suitable values of $\e$.

\subsection{Diagonalization}\label{ss:diag}

\noindent  Integrating by parts, we can evaluate the operator
$T_{\Sig_\e}$ multiplying a test function by the following quantity
\begin{equation}\label{eq:Seu}
    \mathfrak{S}_\e(u) = \sqrt{\det g} \left( - \D_g u + u - p u_{I,\e}^{p-1} u
  \right)
\end{equation}
and integrating in the variables $y$ and $\z$ (using
\eqref{eq:coarea}). In Lemma \ref{l:decae} we studied
$\mathfrak{S}_\e$ acting on the functions $\Psi_l$, for any $l$
fixed. In that lemma, our estimates depend on the value of the index
$l$, and in general one can expect that they become worse and worse
as $l$ increases. The goal of this subsection is to derive estimates
in terms of both $\e$ and $l$ and, evaluating $\mathfrak{S}_\e(u)$
on the functions $\hat{u}_2 \in \hat{H}_2$, we will keep track also
of the terms of order $\e^3$ and higher.

In the following, we will sometimes omit the factor $\chi_\e$
appearing in \eqref{eq:Psil} since this will only produce error
terms exponentially small in $\e$, which are negligible for our
purposes.

\begin{lem}\label{l:sehu2}
There exist linear differential operators $L_1, L_2, L_3$ (acting on
the variables $\oy$) of order $1$, $2$ and $3$ respectively, whose
coefficients (independent of $l$) are smooth and satisfy the bounds
\begin{equation}\label{eq:estci}
    c_\a(L_i) \leq C (1 + |\z|^C) e^{- \frac{|\z|}{C}},
\end{equation}
and such that in local coordinates we have the following expression
for $\mathfrak{S}_\e(\Psi_l)$
\begin{eqnarray}\label{eq:Secar} \nonumber
  \mathfrak{S}_\e(\Psi_l) & = & \e^2 \frac{C_0}{C_1} \mu_l w_0^{p-1} \pa_h w_0 \psi_l^h
  \\ \nonumber & - & 2 \e^3 \left( \z_i
  \G^b_a(E_i) - \zn H_{ab} + \zn H_\a^\a \d_{ab} \right)
  (\pa^2_{\oy_a \oy_b} \psi_l^h) \pa_h w_0 - \e^3 (\pa^2_{\oy_a
  \oy_a} \psi_l^h) \pa_h w_1  \\ & + & \e^3 \zn H_\a^\a \mu_l
  \psi_l^h g_0^h (1 - p w_0^{p-1}) - \e^3 p (p-1) w_0^{p-2} w_1
  \mu_l \psi_l^h g_0^h - \e^4 \mu_l (\pa^2_{\oy_a \oy_a} \psi_l^h)
  g_0^h \\ & + & \e^3 L_1 \psi_l + \e^4 L_3 \psi_l + \e^4 \mu_l L_1
  \psi_l + \e^5 \mu_l L_2 \psi_l, \nonumber
\end{eqnarray}
where $C_0$, $C_1$ are as in Subsection \ref{ss:ae}.
\end{lem}

\begin{pf}
As for the construction of the approximate solutions $u_{I,\e}$, we
can expand formally $\mathfrak{S}_\e(\Psi_l)$ in powers of $\e$ and
check carefully all the error terms, paying particular attention to
the ones involving derivatives in the variables $\oy_a$, which
produce larger and larger terms (as $l$ increases) in the Fourier
modes. When we differentiate with respect to the variables $\z$, the
quantities appearing will be considered as coefficients (depending
smoothly on $\z$, with exponential decay) of the functions $\psi_l$
or their derivatives in $\oy$.

We recall that the functions $w_0$ and $(g_i)_i$ in \eqref{eq:Psil}
are shifted in the $\z'$ variable by the (smooth) normal section
$\Phi(\oy)$. Hence, when differentiating with respect to $\oy$, the
derivatives of $\Phi$ might appear through the chain rule, see also
Subsection \ref{ss:pfpr31}. This fact will be assumed understood,
and it will not be mentioned anymore since it does not create any
serious difficulty, or any difference in the estimates.

By our construction of $\Psi_l$, all the terms multiplying powers of
$\e$ less or equal than $2$ reduce to $\e^2 \frac{C_0}{C_1} \mu_l
\left( - \D_{\z} (\psi_l^h \pa_h w_0) + \psi_l^h \pa_h w_0 \right) =
\e^2 p \frac{C_0}{C_1} \mu_l w_0^{p-1} \pa_h w_0 \psi_l^h$, so we
are left to consider the powers (of $\e$) of order $3$ and higher.
In the remainder of the proof, we use the symbol $\mathcal{A}_2(\e)$
to denote terms of order $1$, $\e$ or $\e^2$: since they all
generate a single term, we do not need to compute them separately.

We begin by considering the terms where derivatives in $\oy$ appear.
Since $\mathfrak{S}_\e$ is linear in $u$, we can deal with each
summand in $\Psi_l$ separately. Looking at $- \sqrt{\det g} \D_g
(\psi^h_l(\oy) \pa_h w_0)$, second derivatives in $\oy$ appear only
in the expression $- \sqrt{\det g} g^{ab} u_{ab}$, so from Lemma
\ref{l:expDgeu} and Remark \ref{r:l3} {\em (b)} we find that
\begin{eqnarray*}
  - \sqrt{\det g} \; \D_g (\psi^h_l(\oy) \pa_h w_0) & = & \mathcal{A}_2(\e)
  - 2 \e^3 \left( \z_i \G_a^b(E_i) - \zn H_{ab} + \zn H_\a^\a
   \right) (\pa^2_{\oy_a \oy_b} \psi^h_l) \pa_h w_0 \\
  & + & \e^3 L_1 \psi_l + \e^4 L_2 \psi_l,
\end{eqnarray*}
where $L_1, L_2$ are as in the statement of the lemma.

Similarly one finds
\begin{eqnarray*}
  - \sqrt{\det g} \; \D_g (\e \psi_l^h(\oy) g_3^h(\oy,\z)) & = & \mathcal{A}_2(\e)
  - \e^3 \pa^2_{\oy_a \oy_a} \psi^h_l \pa_h w_1 + \e^3 L_1 \psi_l + \e^4 L_2
\psi_l;
\end{eqnarray*}
\begin{eqnarray*}
  - \sqrt{\det g} \; \D_g (\e^2
\psi_l^h(\oy) g_2^h(\oy,\z)) & = & \mathcal{A}_2(\e) + \e^4 L_2
\psi_l + \e^3 L_1 \psi_l;
\end{eqnarray*}
\begin{eqnarray*}
  - \sqrt{\det g} \; \D_g (\e^2 \mu_l \psi_l^h(\oy) g_0^h(\oy,\z))
& = & \mathcal{A}_2(\e) - \e^4 (\pa^2_{\oy_a \oy_a} \psi^h_l) g_0 +
\e^4 \mu_l L_1 \psi_l + \e^5 \mu_l L_2 \psi_l.
\end{eqnarray*}
\begin{eqnarray*}
  - \sqrt{\det g} \; \D_g (\e^2 (\pa_{\oy_a} \psi^h_l(\oy)) g_1^h(\oy,\z))
& = & \mathcal{A}_2(\e) + \e^4 L_3 \psi_l.
\end{eqnarray*}

At this point we are left with the terms (of order $\e^3$ and
higher) which do not involve derivatives of $\psi_l$ in $\oy$: these
will appear as multiplicators of the summands in the expression of
$\Psi_l$. The ones involving $\pa_h w_0, g_1, g_2$ and $g_3$ are
included in the expression $\e^3 L_1 \psi_l$, so it remains to
consider $\e^2 \mu_l \psi_l^h g_0^h$. Recalling that  $\sqrt{\det g}
= 1 + \e \z_n H_\a^\a + O(\e^2)$ (see the proof of Lemma
\ref{l:expDgeu}), and expanding $- p u_{I,\e}^{p-1}$ as
\begin{eqnarray} \nonumber
  - p \left[ w_0^{p-1} + \e (p-1) w_0^{p-2} w_1 + \e^2 (p-1) w_0^{p-2} w_2
  + \frac 12 \e^2 (p-1) (p-2) w_0^{p-3} w_1^2 \right] + O(\e^3),
\end{eqnarray}
we obtain
\begin{eqnarray*}
  \sqrt{\det g} (1 - p u_{I,\e}^{p-1}) \e^2 \mu_l \psi_l^h g_0^h & = &
  \mathcal{A}_2(\e) + \e^3 \zn H_\a^\a \mu_l \psi_l^h g_0^h (1 - p w_0^{p-1})
  \\ & - & \e^3 p (p-1) w_0^{p-2} w_1
  \mu_l \psi_l^h g_0^h + \e^4 \mu_l L_0 \psi_l,
\end{eqnarray*}
where $L_0$ is a multiplication operator with coefficients also
satisfy \eqref{eq:estci}. This concludes the proof of the lemma.
\end{pf}

Next, using the above characterization, if $\hat{u}_2$ is a suitable
linear combination of the $\Psi_l$'s, we can estimate the scalar
products of $T_{\Sig_\e} \hat{u}_2$ (in $H_{\Sig_\e}$) with some
other elements belonging to the subspaces $H_1$, $\hat{H}_2$,
$\tilde{H}_2$ and $H_3$, see \eqref{eq:h1}-\eqref{eq:h23}.

\begin{lem}\label{l:prodhu2}
For some arbitrary real coefficients $(\a_l)_l$ and $(\b_l)_l$, we
consider functions $u_1 \in H_1$, $\hat{u}_2 \in \hat{H}_2$ and
$\tilde{u}_2 \in \tilde{H}_2$ of the form
$$
  u_1 = \sum_{j=0}^{\infty} \a_j \phi_j (\e y) u_{j,\e}(|\z|); \qquad
  \quad \hat{u}_2 = \sum_{l=0}^{\e^{- \d}} \b_l \Psi_l; \qquad
\quad \tilde{u}_2 = \sum_{\e^{-\d} + 1}^{\ov{C} \e^{-k}} \b_l
\psi_l^m (\e y) \hat v_{l,\e,m}(\z).
$$
We also let $u_3 \in H_3$. Then, for $\d \in \left( \frac{k}{2} +
\g, \frac 23 k - \g \right)$ and $\g$ sufficiently small, we have
the following relations
\begin{equation}\label{eq:hu2u1}
    (T_{\Sig_\e} \hat{u}_2, u_1)_{H_{\Sig_\e}} = o(\e^2) \left(
    \frac{1}{\e^k} \sum_{l=0}^{\e^{-\d}} |\mu_l| \b_l^2
    \right)^{\frac 12} \|u_1\|_{H_{\Sig_\e}};
\end{equation}
\begin{equation}\label{eq:hu2hu2}
  (T_{\Sig_\e} \hat{u}_2, \hat{u}_2)_{H_{\Sig_\e}} = C_0 (1 + o_\e(1))
  \frac{1}{\e^k} \sum_{l=0}^{\e^{-\d}} \e^2 \mu_l \b_l^2;
\end{equation}
\begin{equation}\label{eq:hu2tu2}
  (T_{\Sig_\e} \hat{u}_2, \tilde{u}_2)_{H_{\Sig_\e}} = O(\e^3)  \frac{1}{\e^k}
\left( \sum_{l=0}^{\e^{-\d}} (\mu_l^2 + \e^2 \mu_l^4) \b_l^2
\right)^{\frac 12} \left( \sum_{l=\e^{-\d} + 1}^{\ov{C} \e^{-k}}
\b_l^2 \right)^{\frac 12} = o(\e^{\frac 43})
\|\hat{u}_2\|_{H_{\Sig_\e}} \|\tilde{u}_2\|_{H_{\Sig_\e}};
\end{equation}
\begin{equation}\label{eq:hu2u3}
  (T_{\Sig_\e} \hat{u}_2, u_3)_{H_{\Sig_\e}} = O(1) \|u_3\|_{H_{\Sig_\e}} \left(
\frac{1}{\e^k} \sum_{l=0}^{\e^{-\d}} \left( \e^6 \mu_l^2 + \e^8
\mu_l^4 \right) \b_l^2 \right)^{\frac 12}.
\end{equation}
\end{lem}

\begin{pf} We recall that, by Lemma \ref{l:comp},
\eqref{eq:rsplnu2}, \eqref{eq:claim1} and \eqref{eq:claim2} there
holds
\begin{eqnarray}\label{eq:splits} \nonumber
   \|u_1\|_{H_{\Sig_\e}}^2 = \frac{1 + o_\e(1)}{\e^k} \sum_{j=0}^\infty \a_j^2;
\qquad \qquad
  \|\hat{u}_2\|_{H_{S_\e}}^2 = \frac{1 + o_\e(1)}{\e^k} \|\pa_1 w_0\|_{H^1(\R^{n+1}_+)}^2
  \sum_{l=0}^{\e^{- \d}} \b_l^2; \\ \\ \|\tilde{u}_2\|_{H_{S_\e}}^2 = \frac{1 + o_\e(1)}{\e^k}
 \sum_{l=\e^{-\d} + 1}^{\ov{C} \e^{-k}} \b_l^2. \nonumber
\end{eqnarray}

We show first \eqref{eq:hu2u1}. Since $u_1$ is even in $\z'$, when
we use the expression of $\mathfrak{S}_\e(\Psi_l)$ in
\eqref{eq:Secar} we have to consider only $- 2 \e^3\z_i \G^b_a(E_i)
\pa^2_{\oy_a \oy_b} \psi_l^h \pa_j w_0 = \e^3 L_2 \psi_l$ and the
errors $\e^s L_t \psi_l$, since the products of all the others terms
with $u_1$ will vanish by oddness. Therefore we leave this term as
it is, and we estimate the error terms only. So we get
$$
  (T_{\Sig_\e} \hat{u}_2, u_1)_{H_{\Sig_\e}} = \frac{1}{\e^k} \sum_{j,l} \a_j
  \b_l \int_K \int_{\R^{n+1}_+} u_{j,\e}(|\z|) \phi_j(\oy) \left( \e^3 L_2 \psi_l
+ \e^4 L_3 \psi_l + \e^4 \mu_l L_1 \psi_l
   + \e^5 \mu_l L_2 \psi_l \right) d \oy d\z.
$$
Reasoning as in Lemma \ref{l:abstr} (avoiding the scaling in $\e$,
which has been already taken care of) one can show that, for any
integer $m$
\begin{equation}\label{eq:gf}
  \int_K \int_{\R^{n+1}_+} \left( \sum_{l=0}^{\e^{-\d}} \b_l L_m \psi_l
\right)^2 \leq C \sum_{l=0}^{\e^{-\d}} (1 + |\mu_l|^m \b_l^2).
\end{equation}
From the H\"older inequality and the last three formulas we deduce
that
$$
  (T_{\Sig_\e} \hat{u}_2, u_1)_{H_{\Sig_\e}} \leq C
\|u_1\|_{H_{\Sig_\e}} \left[ \frac{1}{\e^k} \sum_{l=0}^{\e^{-\d}}
\left( \e^6 (1 + |\mu_l|^2) + \e^8 |\mu_l|^3 + \e^{10} |\mu_l|^4
\right) \b_l^2 \right]^{\frac 12}.
$$
Now, from the Weyl's asymptotic formula and from the fact that $\d
\in \left( \frac{k}{2} + \g, \frac 23 k - \g \right)$, one finds
that for $l \leq \e^{-\d}$ there holds $\e^2 |\mu_l|^2 = o_\e(1)
|\mu_l|$, that $\e^4 |\mu_l|^3 = o_\e(1)$ and that $\e^6 |\mu_l|^4 =
o_\e(1)$, so \eqref{eq:hu2u1} follows.

We turn now to \eqref{eq:hu2hu2}. It is convenient first to evaluate
some $L^2$ norms. Writing $\mathfrak{S}_\e(\Psi_l) = \e^2 p
\frac{C_0}{C_1} \mu_l w_0^{p-1} \pa_h w_0 \psi_l^h +
\tilde{\mathfrak{S}}_\e(\Psi_l)$, and $\Psi_l = \chi_\e(|\z|)
\psi_l^h \pa_h w_0 + \ov{\Psi}_l$, from \eqref{eq:gf} we find ($l$
runs between $0$ and $\e^{-\d}$)
\begin{equation}\label{eq:L2normPsil}
    \left\|\sum \b_l \Psi_l\right\|_{L^2}^2, \left\|\sum \b_l \psi_l^h \pa_h w_0
   \right\|_{L^2}^2
\leq \frac{C}{\e^k} \sum_l \left( 1 + \e^2 + \e^4 |\mu_l|^2 \right)
\b_l^2 \leq \frac{C}{\e^k} \sum_l \b_l^2;
\end{equation}
\begin{equation}\label{eq:L2normovPsil}
    \left\|\sum \b_l \ov{\Psi}_l\right\|_{L^2}^2 \leq \frac{C}{\e^k} \sum_l \left(
\e^2 + \e^4 |\mu_l|^2 \right) \b_l^2 \leq \frac{C}{\e^k} \e^2 \sum_l
(1 + \e^2 \mu_l^2) \b_l^2;
\end{equation}
\begin{equation}\label{eq:L2normfsePsil}
    \left\|\sum \b_l \mathfrak{S}_\e(\Psi_l) \right\|_{L^2}^2 \leq \frac{C}{\e^k} \sum_l \left(
\e^4 |\mu_l|^2 + \e^6 |\mu_l|^2  + \e^8 |\mu_l|^4 + \e^{10}
|\mu_l|^4 \right) \b_l^2 \leq \frac{C}{\e^k} \e^4 \sum_l \mu_l^2
\b_l^2;
\end{equation}
\begin{equation}\label{eq:L2normffsePsil}
    \left\|\sum \b_l \tilde{\mathfrak{S}}_\e(\Psi_l) \right\|_{L^2}^2 \leq
\frac{C}{\e^k} \sum_l \left( \e^6 |\mu_l|^2  + \e^8 |\mu_l|^4 +
\e^{10} |\mu_l|^4 \right) \b_l^2 \leq \frac{C}{\e^k} \e^6 \sum_l
(|\mu_l|^2 + \e^2 |\mu_l|^4 ) \b_l^2.
\end{equation}
Using the orthogonality of the $\psi_l$'s, \eqref{eq:coarea} and
recalling the definition of $C_1$ in Subsection \eqref{ss:ae}, we
find
\begin{equation}\label{eq:prodhu2hu2}
    (T_{\Sig_\e}(\Psi_l), \Psi_j)_{H_{\Sig_\e}} = \e^2 C_0
\mu_l \d_{lj} + (\tilde{\mathfrak{S}}_\e(\Psi_l), \psi_j^h \pa_h
w_0)_{L^2} + (\mathfrak{S}_\e(\Psi_l), \ov{\Psi}_j)_{L^2}.
\end{equation}
Multiplying by the coefficients $\b$'s, using the H\"older
inequality and \eqref{eq:L2normPsil}-\eqref{eq:L2normffsePsil} we
get
\begin{eqnarray*}
  (T_{\Sig_\e} \hat{u}_2, \hat{u}_2)_{H_{\Sig_\e}} & = & C_0 \sum_l
\e^2  \mu_l \b_l^2 + \frac{1}{\e^k} O(\e^3) \left[ \left( \sum_l
(\mu_l^2 + \e^2 \mu_l^4) \b_l^2 \right)^{\frac 12} \left( \sum_l
\b_l^2 \right)^{\frac 12} \right. \\
  & + & \left. \left( \sum_l \mu_l^2 \b_l^2
\right)^{\frac 12} \left( \sum_l (1 + \e^2 \mu_l^2) \b_l^2
\right)^{\frac 12} \right].
\end{eqnarray*}
Recalling the Weyl's asymptotic formula and the fact that $\d \in
\left( \frac{k}{2} + \g, \frac 23 k - \g \right)$, we obtain $\e^2
\mu_l^2 = o(\mu_l)$, $\e^4 \mu_l^4 = o(\mu_l)$ for $l \leq
\e^{-\d}$, so the last formula implies \eqref{eq:hu2hu2}.

To prove \eqref{eq:hu2tu2} we notice that, by the orthogonality of
the $\psi_l$'s, the term of order $\e^2$ in
$\mathfrak{S}_\e(\Psi_l)$, once multiplied by $\tilde{u}_2$ and
integrated, vanishes identically. Therefore, from the H\"older
inequality, \eqref{eq:splits} and \eqref{eq:L2normffsePsil} we find
$$
  (T_{\Sig_\e} \hat{u}_2, \tilde{u}_2)_{H_{\Sig_\e}} = O(\e^3) \frac{1}{\e^k}
\left( \sum_{l=0}^{\e^{-\d}} (\mu_l^2 + \e^2 \mu_l^4) \b_l^2
\right)^{\frac 12} \left( \sum_{l=\e^{-\d} + 1}^{\ov{C} \e^{-k}}
\b_l^2 \right)^{\frac 12},
$$
which is precisely \eqref{eq:hu2tu2}.

It remains to prove \eqref{eq:hu2u3}.  Using \eqref{eq:explaplu},
the formulas in the proof of Lemma \ref{l:expDgeu} and the fact that
(linearizing \eqref{eq:w0} at $w_0$) $- \D_{\z} (\pa_h w_0) + \pa_h
w_0 = p w_0^{p-1} \pa_h w_0$, one finds
\begin{eqnarray}\label{eq:normdual}\nonumber
  \sqrt{\det g_\e} (-\D_{g_\e} \Psi_l + \Psi_l) & = & p w_0^{p-1}
  \psi_l^h \pa_h w_0 + \e L_0 \psi_l + \e^2 (L_2 \psi_l +
  \mu_l L_0 \psi_l) + \e^3 L_2 \psi_l \\
   & + & \e^4(\mu_l L_2 \psi_l + L_3 \psi_l).
\end{eqnarray}
Hence from \eqref{eq:Secar} it follows that
\begin{eqnarray*}
  \mathfrak{S}_\e(\Psi_l) &=& \e^2 \frac{C_0}{C_1} \frac{\mu_l}{p} \sqrt{\det g_\e}
  (-\D_{g_\e} \Psi_l + \Psi_l) + \e^3 \mu_l L_0 \psi_l + \e^4 \mu_l (L_2 \psi_l +
  \mu_l L_0 \psi_l) \\ & + & \e^5 \mu_l L_2 \psi_l + \e^6 \mu_l (\mu_l L_2 \psi_l + L_3
\psi_l) + \tilde{\mathfrak{S}}_\e(\Psi_l).
\end{eqnarray*}
Since $u_3$ is orthogonal to $\hat{H}_2$ in $H_{\Sig_\e}$,
integrating by parts we have $\int_{\Sig_\e} u_3 (-\D_{g_\e} \Psi_l
+ \Psi_l) \sqrt{\det g_\e} dy d\z = 0$ for $l = 0, \dots, \e^{-\d}$.
Hence from  \eqref{eq:gf} and \eqref{eq:L2normffsePsil} we get
$$
  (T_{\Sig_\e} \hat{u}_2, u_3)_{H_{\Sig_\e}} =
O(1) \|u_3\|_{H_{\Sig_\e}} \left( \frac{1}{\e^k}
\sum_{l=0}^{\e^{-\d}} \left( \e^6 \mu_l^2 + \e^8 \mu_l^4 + \e^{12}
\mu_l^6 \right) \b_l^2 \right)^{\frac 12}.
$$
As shown before, $\e^2 \mu_l^2 = o_\e(1)$ for $l \leq \e^{-\d}$, so
we have $\e^{12} \mu_l^6 = o(\e^8 \mu_l^4)$, and the conclusion
holds.
\end{pf}

\

\noindent We have now the counterpart of Lemma \ref{l:prodhu2} with
$\tilde{u}_2$ replacing $\hat{u}_2$.

\begin{lem}\label{l:prodtu2}
For some arbitrary real coefficients $(\a_l)_l$ and $(\b_l)_l$, we
consider functions $u_1 \in H_1$, $\hat{u}_2 \in \hat{H}_2$ and
$\tilde{u}_2 \in \tilde{H}_2$ of the form
$$
  u_1 = \sum_{j=0}^{\infty} \a_j \phi_j (\e y) u_{j,\e}(|\z|); \qquad
  \quad \hat{u}_2 = \sum_{l=0}^{\e^{- \d}} \b_l \Psi_l; \qquad
\quad \tilde{u}_2 = \sum_{l=\e^{-\d} + 1}^{\ov{C} \e^{-k}} \b_l
\psi_l^m (\e y) \hat v_{l,\e,m}(\z).
$$
Suppose also that $u_3 \in H_3$. Then, for $\d \in \left(
\frac{k}{2} + \g, \frac 23 k - \g \right)$ and $\g$ sufficiently
small, we have the following relations
\begin{equation}\label{eq:tu2u1}
    (T_{\Sig_\e} \tilde{u}_2, u_1)_{H_{\Sig_\e}} = O(\e^{1-\g})
\|u_1\|_{H_{\Sig_\e}} \left( \frac{1}{\e^k} \sum_{l=\e^{-\d} +
1}^{\ov{C} \e^{-k}} \b_l^2 \right)^{\frac 12};
\end{equation}
\begin{equation}\label{eq:tu2tu2}
  (T_{\Sig_\e} \tilde{u}_2, \tilde{u}_2)_{H_{\Sig_\e}} \geq \frac{C^{-1}}{\e^k}
  \sum _{l=\e^{-\d} + 1}^{\ov{C} \e^{-k}} \e^2 \mu_l \b_l^2;
\end{equation}
\begin{equation}\label{eq:tu2u3}
  (T_{\Sig_\e} \tilde{u}_2, u_3)_{H_{\Sig_\e}} = O(\e^{1-\g})
  \|u_3\|_{H_{\Sig_\e}} \left( \frac{1}{\e^k} \sum _{l=\e^{-\d}
  + 1}^{\ov{C} \e^{-k}} \b_l^2 \right)^{\frac 12}.
\end{equation}
\end{lem}

\begin{pf} We show first \eqref{eq:tu2u1}. Since $u_1$ and
$\tilde{u}_2$, for any fixed $y$ are linear combinations of
spherical harmonics (in $\frac{\z}{|\z|}$) of different type, from
the arguments of Subsection \ref{ss:model} it follows that
$$
  (u_1, \tilde{u}_2)_{H_{S_\e}} = 0; \qquad \qquad \int_{S_\e}
  w_0^{p-1}(|\z|) u_1 \tilde{u}_2 d V_{\tilde{g}_\e} = 0,
$$
so we clearly have that $(T_{S_\e} u_1, \tilde{u}_2)_{H_{S_\e}} =
0$. Then \eqref{eq:tu2u1} follows immediately from Lemma
\ref{l:comp}.

%From \eqref{eq:explaplu}, the expansions in the proof of Lemma
%\ref{l:expDgeu} and the fact that $-\D_K^N \psi_l = \mu_l \psi_l +
%(\mathfrak{R} - \mathfrak{B}) \psi_l$, one finds that
%\begin{eqnarray}\label{eq:Setu2} \nonumber
%    \mathfrak{S}_\e(\psi_l^h(\e y) v_{l,\e,h}(\z)) & = & \psi_l^h \left(
%    - \D_\z v_{l,\e,h} + (1 + \e^2 \mu_l) v_{l,\e,h} - p w_0^{p-1}
%    v_{l,\e,h} \right) \\ & + & \e L_0 \psi_l + \e^2 L_1 \psi_l + \e^3 L_2
%\psi_l + \e^3 \mu_l L_0 \psi_l.
%\end{eqnarray}
%Multiplying $\mathfrak{S}_\e(\tilde{u}_2)$ by $u_1$ and integrating,
%since $v_{l,\e,h}$ is odd in $\z'$ and since $u_1$ is even, by
%cancelation, by the H\"older inequality, \eqref{eq:splits} and by
%\eqref{eq:gf} we obtain
%$$
%  (T_{\Sig_\e} \hat{u}_2, u_1)_{H_{\Sig_\e}} = O(1)
%\|u_1\|_{H_{\Sig_\e}} \left( \sum_{\e^{-\d} + 1}^{\ov{C} \e^{-k}}
%(\e^2 + \e^4 |\mu_l| + \e^6 \mu_l^2) \b_l^2 \right)^{\frac 12}.
%$$
%By the Weyl's asymptotic formula, and since $l \leq \ov{C}
%\e^{-k}$ ($\ov{C}$ is assumed small and, say, less than $1$), we
%have that $\e^2 \mu_l$ stays uniformly bounded, so \eqref{eq:tu2u1}
%follows.

To prove \eqref{eq:tu2tu2}, we reason as for the proof of Lemma
\ref{l:ntu2} to find
\begin{equation}\label{eq:eee}
    (T_{S_\e} \tilde{u}_2, w)_{H_{S_\e}} = \tilde{A}_1 + \tilde{A}_2 + \tilde{A}_3,
\end{equation}
where $w \in H_{S_\e}$ is arbitrary, and where
$$
  \tilde{A}_1(w) = \int_{S_\e}\sum\limits_{l=\e^{-\d}+1}^{\ov{C} \e^{-k}}
\left[ \left( - \D_{\z} + (1+\e^2\o_l) - p w_0^{p-1} \right)\bigg(
\sum\limits_{m=1}^n\b_l \psi_l^m (\e y) \hat
v_{l,\e}(|\z|)\,\frac{\z_m}{|\z|} \bigg) \right] w;
$$
$$
 \tilde{A}_2(w) = \e^2 \int_{S_\e}\sum\limits_{l=\e^{-\d}+1}^{\ov{C}
\e^{-k}}\bigg( \sum\limits_{m=1}^n\b_l \left((\mathfrak{B} -
  \mathfrak{R})\psi_l  \right)^m (\e y) \hat v_{l,\e}(|\z|)\,\frac{\z_m}{|\z|}
  \bigg) w;
$$
$$
  \tilde{A}_3(w) = \e^2 \int_{S_\e}\sum\limits_{l=\e^{-\d}+1}^{\ov{C}
\e^{-k}}\bigg( \sum\limits_{m=1}^n\b_l (\mu_l - \o_l) \psi_l^m
 (\e y) \hat v_{l,\e}(|\z|)\,\frac{\z_m}{|\z|}
  \bigg) w;
$$
As for \eqref{eq:estA2}, since $|\mu_l - \o_l|$ is uniformly bounded
one finds
\begin{equation}\label{eq:esttA2tA3}
  |\tilde{A}_2(w)| + |\tilde{A}_3(w)| \leq C \e^2
  \|\tilde{u}_2\|_{H_{S_\e}} \|w\|_{H_{S_\e}}
\end{equation}
for a fixed positive constant $C$. Taking $w = \tilde{u}_2$, by the
orthogonality of the $\psi_l$'s, by the fact that $T_{\e^2 \o_l}
v_{l,\e,m} = \s_{\e^2 \o_l, \e} v_{l,\e,m}$ (see Proposition
\ref{p:eicompe}) and by \eqref{eq:normeigenf}, with an integration
by parts we have
$$
  \tilde{A}_1(\tilde{u}_2) = \frac{1}{\e^k} \sum_{l = \e^{-\d}+1}^{\ov{C}
  \e^{-k}} \s_{\e^2 \o_l, \e} \b_l^2 \|v_{l,\e,1}\|_{\e^2 \o_l,
  \e} = \frac{1}{\e^k} \sum_{l = \e^{-\d}+1}^{\ov{C}
  \e^{-k}} \s_{\e^2 \o_l, \e} \b_l^2.
$$
From \eqref{eq:weyl3}, Proposition \ref{p:ta} and Proposition
\ref{p:eicompe}, which provide estimates on $\s_{\e^2 \o_l,\e}$, we
obtain
\begin{equation}\label{eq:tA1tu2}
    \tilde{A}_1(\tilde{u}_2) \geq \frac{C^{-1}}{\e^k} \sum _{l=\e^{-\d}
    + 1}^{\ov{C} \e^{-k}} \e^2 \mu_l \b_l^2
\end{equation}
for some fixed $C > 0$. Then \eqref{eq:tu2tu2} follows from
\eqref{eq:esttA2tA3}, \eqref{eq:tA1tu2}, Lemma \ref{l:ntu2} and
Lemma \ref{l:comp} (since $\e^2 \mu_l \gg \e^{1-\g}$ for $l >
\e^{-\d}$ and for $\g$ sufficiently small).

%$\mathfrak{S}_\e(\psi_l^h(\e y) v_{l,\e,h}(\z))$ by $\psi_j^i(\e y)
%v_{j,\e,i}(\z)$ and integrate. From the first term in the right-hand
%side of \eqref{eq:Setu2} by \eqref{eq:coarea} we obtain
%$$
%  \d_{l,j} \int_{\R^{n+1}} \left[ |\n_\z v_{l,\e,1}|^2 + \left(1 + \e^2 \mu_l -
%p w_0^{p-1} \right) v_{l,\e,1}^2 \right],
%$$
%where we have used the orthogonality of the $\psi_j$'s. Recalling
%the notation in Subsection \ref{ss:spa}, \eqref{eq:normeigenf} and
%the fact that $\mu_l = \o_l + O(1)$, the latter quantity becomes
%$$
%  \d_{lj} (T_{\e^2 \mu_l,\e} v_{l,\e,1}, v_{l,\e,1})_{\e^2 \mu_l,
%\e} = \d_{lj} \left[ (T_{\e^2 \o_l,\e} v_{l,\e,1}, v_{l,\e,1})_{\e^2
%\o_l, \e} + O(\e^2) \right] = \d_{lj} (\s_{\e^2 \o_l,\e} + O(\e^2)).
%$$
%Using \eqref{eq:muaae} one then has $\s_{\e^2 \o_l,\e} \geq C^{-1}
%\e^2 l^{\frac{2}{k}} \geq C^{-1} \e^2 \mu_l$, for $l \in \{
%\e^{-\d} + 1, \dots, \ov{C} \e^{-k} \}$. From these considerations
%and from \eqref{eq:gf} then one finds
%$$
%  (T_{\Sig_\e} \tilde{u}_2, \tilde{u}_2)_{H_{\Sig_\e}} \geq \sum _{\e^{-\d} + 1}^{\ov{C}
%\e^{-k}} \e^2 \mu_l \b_l^2 + O(1) \left( \sum_{\e^{-\d} +
%1}^{\ov{C} \e^{-k}} (\e^2 + \e^4 |\mu_l| + \e^6 \mu_l^2) \b_l^2
%\right)^{\frac 12} \left( \sum_{\e^{-\d} + 1}^{\ov{C} \e^{-k}}
%\b_l^2 \right)^{\frac 12}.
%$$
%Since for $l \in \{ \e^{-\d} + 1, \dots, \ov{C} \e^{-k} \}$ there
%holds $\e \ll \e^2 |\mu_l| \leq C$, \eqref{eq:tu2tu2} follows.

We turn now to \eqref{eq:tu2u3}. By \eqref{eq:esttA2tA3}, taking $w
= u_3$, it is sufficient to estimate $\tilde{A}_1 (u_3) +
\tilde{A}_3 (u_3)$. From $T_{\e^2 \o_l} v_{l,\e,m} = \s_{\e^2 \o_l,
\e} v_{l,\e,m}$ in $H_{\e^2 \o_l, \e}$, with an integration by parts
we find
$$
  \tilde{A}_1(u_3) + \tilde{A}_3 (u_3) =
  \int_{S_\e}\sum\limits_{l=\e^{-\d}+1}^{\ov{C} \e^{-k}}
\s_{\e^2 \o_l }\left[ \left( - \D_{\z} + (1+\e^2\mu_l) - p w_0^{p-1}
\right)\bigg( \sum\limits_{m=1}^n\b_l \psi_l^m (\e y) \hat
v_{l,\e}(|\z|)\,\frac{\z_m}{|\z|} \bigg) \right] u_3.
$$
From \eqref{eq:lapltg} and from the fact that $-\D_K^N \psi_l =
\mu_l \psi_l + (\mathfrak{R} - \mathfrak{B}) \psi_l$, one finds
$$
  \e^2 \mu_l \psi_l^m \hat v_{l,\e}(|\z|)\,\frac{\z_m}{|\z|} = - \e^2
  \D_K^N \psi_l^m \hat v_{l,\e}(|\z|)\,\frac{\z_m}{|\z|} + \e^2 ((\mathfrak{R} -
  \mathfrak{B}) \psi_l)^m \hat v_{l,\e}(|\z|)\,\frac{\z_m}{|\z|}.
$$
Therefore, integrating by parts we obtain
\begin{equation}\label{eq:eeee}
    \tilde{A}_1(u_3) + \tilde{A}_3 (u_3) = (\tilde{U}_2,
  u_3)_{H_{S_\e}} + \tilde{A}_4(u_3),
\end{equation}
where
$$
  \tilde{A}_4(u_3) = \e^2 \int_{S_\e}\sum\limits_{l=\e^{-\d}+1}^{\ov{C}
  \e^{-k}}\bigg( \sum\limits_{m=1}^n \s_{\e^2 \o_l, \e} \b_l
  \left((\mathfrak{B} - \mathfrak{R})\psi_l  \right)^m (\e y) \hat
  v_{l,\e}(|\z|)\,\frac{\z_m}{|\z|} \bigg) u_3,
$$
and where $\tilde{U}_2 = \sum_{\e^{-\d} + 1}^{\ov{C} \e^{-k}}
\s_{\e^2 \o_l, \e} \b_l \psi_l^m (\e y) \hat v_{l,\e,m}(\z) \in
H_2$. Now, as for $\tilde{u}_2$ it is possible to prove that there
exists a fixed $C > 0$ such that
$$
\|\tilde{U}_2\|_{H_{S_\e}}^2 \leq \frac{C}{\e^k}
\sum_{l=\e^{-\d}+1}^{\ov{C} \e^{-k}} \s_{\e^2 \o_l,\e}^2 \b_l^2 \leq
\frac{C}{\e^k} \sum_{l=\e^{-\d}+1}^{\ov{C} \e^{-k}} \b_l^2,
$$
where we used the fact that $\s_{\e^2 \o_l,\e}$ is uniformly bounded
for $l \leq \ov{C} \e^{-k}$. Since $u_3$ is orthogonal in
$H_{\Sig_\e}$ to $H_2$, from Lemma \ref{l:comp}, these observations
and the last two formulas it follows that
$$
  (\tilde{U}_2, u_3)_{H_{S_\e}} = O(\e^{1-\g}) \|\tilde{U}_2\|_{H_{S_\e}}
  \|u_3\|_{H_{S_\e}} \leq C \e^{1-\g} \left( \sum_{l=\e^{-\d}+1}^{\ov{C}
  \e^{-k}} \b_l^2 \right)^{\frac 12} \|u_3\|_{H_{S_\e}}.
$$
The arguments of the proof of Lemma \ref{l:ntu2} yield
$\tilde{A}_4(u_3) \leq C \e^4 \left( \sum_{l=\e^{-\d}+1}^{\ov{C}
\e^{-k}} \b_l^2 \right)^{\frac 12} \|u_3\|_{H_{S_\e}}$. Hence from
\eqref{eq:eee}, \eqref{eq:eeee} and Lemma \ref{l:comp} we find that
$$
  (T_{\Sig_\e} \tilde{u}_2, u_3)_{H_{\Sig_\e}} = (\tilde{U}_2,
  u_3)_{H_{S_\e}} + O(\e^{1-\g}) \left( \sum_{l=\e^{-\d}+1}^{\ov{C}
  \e^{-k}} \b_l^2 \right)^{\frac 12} \|u_3\|_{H_{S_\e}},
$$
which concludes the proof.
\end{pf}

\subsection{Applications}\label{ss:appl}

In this subsection we apply the estimates in Lemmas \ref{l:comp},
\ref{l:prodhu2} and \ref{l:prodtu2} to estimate the morse index of
$T_{\Sig_\e}$ as $\e$ tends to zero, and to characterize the
eigenfunctions of $T_{\Sig_\e}$ corresponding to {\em resonant}
eigenvalues.

From Proposition \ref{p:ta} we know that there exists a unique
positive number $\ov{\a}$ such that $\eta_{\ov{\a}} = 0$. If $C_k$
is the constant given in \eqref{eq:weyl1}, we also let
\begin{equation}\label{eq:Theta}
    \Theta = \left( \frac{\ov{\a}}{C_k} \right)^{\frac k2}\,Vol(K).
\end{equation}
Then we have the following result.

\begin{pro}\label{p:morse}
Let $\Theta$ be the constant given in \eqref{eq:Theta}, and let
$T_{\Sig_\e}$ be the operator given in \eqref{eq:TSige}. Then, as
$\e$ tends to zero, the Morse index of $T_{\Sig_\e}$ is asymptotic
to $\Theta \e^{-k}$.
\end{pro}

%
%\begin{rem}\label{r:morse}
%We point out that by Corollary \ref{c:eigenv} and Proposition
%\ref{p:morse} the Morse indices of $T_{\Sig_\e}$ and $T_{S_\e}$
%satisfy the same asymptotics. Looking at the explicit structure of
%the spectrum of $T_{S_\e}$, this is not apparent from Lemma
%\ref{l:comp}. In proving the estimate of Proposition \ref{p:morse}
%we use crucially the accurate estimates of Lemmas \ref{l:prodhu2}
%and \ref{l:prodtu2}.
%\end{rem}

\begin{pf}
%{\sc of Proposition \ref{p:morse}}
For any $m \in \N$, the $m$-th eigenvalue $\l_m$ of $T_{\Sig_\e}$,
and the $m$-th eigenvalue $\tilde{\l}_m$ of $T_{S_\e}$ can be
evaluated via the classical Rayleigh quotients
\begin{equation}\label{eq:Ray}
    \l_m = \inf_{\text{dim} M_m = m} \sup_{u \in M_m}
   \frac{(T_{\Sig_\e} u, u)_{H_{\Sig_\e}}}{(u, u)_{H_{\Sig_\e}}};
   \qquad \qquad \tilde{\l}_m = \inf_{\text{dim} M_m = m} \sup_{u \in M_m}
   \frac{(T_{S_\e} u, u)_{H_{S_\e}}}{(u, u)_{H_{S_\e}}}
\end{equation}
where $M_m$ is a vector subspace of $H_{\Sig_\e}$. Choosing $M_m =
\tilde{M}_m$ to be the span of the first $m$ eigenfunctions of
$T_{S_\e}$, from the above formula for $\l_m$ and from Lemma
\ref{l:comp} we get
$$
  \l_m \leq \sup_{u \in \tilde{M}_m}
   \frac{(T_{\Sig_\e} u, u)_{H_{\Sig_\e}}}{(u, u)_{H_{\Sig_\e}}} =
   \sup_{u \in \tilde{M}_m} \frac{(T_{S_\e} u, u)_{H_{S_\e}} + O(\e^{1-\g})
   (u, u)_{H_{S_\e}}}{(1+O(\e^{1-\g}))(u, u)_{H_{S_\e}}} \leq
   \tilde{\l}_m + O(\e^{1-\g}).
$$
Reasoning in the same way we also find $\tilde{\l}_m \leq \l_m +
O(\e^{1-\g})$, and hence it follows that
\begin{equation}\label{eq:llllll}
    |\l_m - \tilde{\l}_m| \leq C \e^{1-\g} \qquad \qquad \hbox{ for
  all } m \in \N \hbox{ and for } \e \hbox{ small},
\end{equation}
where $C > 0$ is a fixed constant.

Now we let $N_1(\e)$ denote the number of eigenvalues $\tilde{\l}_m$
less or equal than $- \e^{\frac{1-\g}{2}}$, and by $N_2(\e)$ the
number of eigenvalues $\tilde{\l}_m$ less or equal than
$\e^{\frac{1-\g}{2}}$. From Proposition \ref{p:dec} it follows that
$N_1(\e)$ is the number of the $\eta_{l,\e}$'s which are smaller
than $- \e^{\frac{1-\g}{2}}$. Reasoning as in Corollary
\ref{c:eigenv} one finds that, as $\e$ tends to zero
$$
  N_1(\e) \simeq \left( \frac{\ov{\a}}{C_k} \right)^{\frac{k}{2}}
  Vol(K) \e^{-k}.
$$
On the other hand, still by Proposition \ref{p:dec} we have that
$N_2(\e) = N_{2,1}(\e) + N_{2,2}(\e)$, where $N_{2,1}(\e)$ is the
number of $\eta_{l,\e}$'s which are smaller than
$\e^{\frac{1-\g}{2}}$, and $N_{2,\e}$ the number of $\s_{l,\e}$'s
which are smaller than $\e^{\frac{1-\g}{2}}$. From \eqref{eq:weyl1},
\eqref{eq:weyl3} and Proposition \ref{p:eicompe} we obtain, for $\e$
small
$$
  N_{2,1}(\e) \simeq \left( \frac{\ov{\a}}{C_k} \right)^{\frac{k}{2}}
  Vol(K) \e^{-k}; \qquad \qquad N_{2,2}(\e) \simeq \left( \frac{1}{C_{N-1,k}}
  \right)^{\frac{k}{2}} Vol(K) \e^{\frac{k(1-\g)}{4}-k} = o(\e^{-k}).
$$
From the last formula we deduce that also
$$
  N_2(\e) \simeq \left( \frac{\ov{\a}}{C_k} \right)^{\frac{k}{2}} Vol(K) \e^{-k}.
$$
Since by \eqref{eq:llllll} the Morse index of $T_{\Sig_\e}$ is
between $N_1(\e)$ and $N_2(\e)$, the conclusion follows.
\end{pf}

\

\noindent We can now characterize the eigenfunctions of
$T_{\Sig_\e}$ corresponding to eigenvalues close to zero.

\begin{pro}\label{p:ei}
For $\e$ sufficiently small, let $\l$ be an eigenvalue of
$T_{\Sig_\e}$ such that $|\l| \leq \e^{\varsigma}$, for some
$\varsigma > 2$, and let $u \in H_{\Sig_\e}$ be an eigenfunction of
$T_{\Sig_\e}$ corresponding to $\l$ with $\|u\|_{H_{\Sig_\e}} = 1$.
In the above notation, let $u = u_1 + u_2 + u_3$, with $u_i \in
H_i$, $i = 1, 2, 3$. Then, if $u_1 = \sum_{j=0}^\infty \a_j
\phi_j(\e y) u_{j,\e}(|\z|)$, one has
\begin{equation}\label{eq:eigenf}
  \left\| u - \sum_{\left\{j : |\eta_{j,\e}| \leq \e^{\frac{1-\g}{2}} \right\}}
  \a_j \phi_j u_{j,\e} \right\|_{H_{\Sig_\e}} \to 0 \qquad \quad
\hbox{ as } \e \to 0.
\end{equation}

\end{pro}

\begin{pf}
We show that $u_2, u_3$ tend to zero as $\e$ tends to zero. This
clearly implies $\|u - u_1\|_{H_{\Sig_\e}} \to 0$. Once this
verified, \eqref{eq:eigenf} can be proved as in \cite{malm2}
Proposition 4.1.

To prove that $u_3$ tends to zero as $\e \to 0$, we take the scalar
product of the eigenvalue equation $T_{\Sig_\e} u = \l u$ with
$u_3$. Using the above arguments (in particular Lemma \ref{l:comp})
we easily find
$$
  \frac{1}{C \ov{C}^{\frac{2}{k}}} \|u_3\|^2_{H_{\Sig_\e}} +
O(\e^{1-\g}) \|u\|_{H_{\Sig_\e}} \|u_3\|_{H_{\Sig_\e}} \leq
(T_{\Sig_\e} u, u_3)_{H_{\Sig_\e}} = \l(u, u_3)_{H_{\Sig_\e}} = \l
\|u_3\|_{H_{\Sig_\e}}^2.
$$
This implies $\|u_3\|_{H_{\Sig_\e}}^2 = O(\e^{1-\g})
\|u\|_{H_{\Sig_\e}} \|u_3\|_{H_{\Sig_\e}}$, and hence
$\|u_3\|_{H_{\Sig_\e}} \leq C \e^{1-\g} \|u\|_{H_{\Sig_\e}} \leq C
\e^{1-\g}$.

Next we take the scalar product of the eigenvalue equation with
$u_2$. From Lemmas \ref{l:prodhu2} and \ref{l:prodtu2} we find
\begin{eqnarray*} \nonumber
% \nonumber to remove numbering (before each equation) \nonumber
  (T_{\Sig_\e} u_2, u_2)_{H_{\Sig_\e}} & \geq & \frac{C_0 (1 + o_\e(1))}{\e^k}
\sum_{l=0}^{\e^{-\d}}
  \e^2 \mu_l \b_l^2 +  \frac{O(1)}{\e^k} \left( \e^5 \sum_{l=0}^{\e^{-\d}}
(\mu_l^2 + \e^2 \mu_l^4) \b_l^2 \right)^{\frac 12} \left( \e
\sum_{l=\e^{-\d}+1}^{\ov{C} \e^{-k}} \b_l^2 \right)^{\frac 12}  \\
  & + & \frac{C^{-1}}{\e^k} \sum_{l=\e^{-\d}+1}^{\ov{C} \e^{-k}}
\e ^2 \mu_l \b_l^2.
\end{eqnarray*}
Since $\e^5 \mu_l^2 + \e^7 \mu_l^4 = o_\e(1) |\mu_l|$ for $l \leq
\e^{-\d}$ and $\e = o(\e^2 \mu_l)$ for $l > \e^{-\d}$ (recall that
$\d \in \left( \frac{k}{2} + \g, k - \g \right)$), it follows that
\begin{equation}\label{eq:Tu2u2u2}
 (T_{\Sig_\e} u_2, u_2)_{H_{\Sig_\e}} \geq C^{-1} \frac{1}{\e^k}
 \sum_{l=0}^{\ov{C} \e^{-k}} \e^2 \mu_l \b_l^2
\end{equation}
for a fixed positive constant $C$. Finally, still from Lemmas
\ref{l:prodhu2}-\ref{l:prodtu2}, from the fact that $\e^4 |\mu_l| +
\e^6 |\mu_l|^3 = o_\e(1)$ for $l \leq \e^{-\d}$ and $\e^{2-2\g} =
o(\e^2 \mu_l) \gg 1$ for $l > \e^{-\d}$ (taking $\g$ sufficiently
small) we have also that
\begin{equation}\label{eq:asasas}
    (T_{\Sig_\e} u_2, u_1 + u_3)_{H_{\Sig_\e}} = o_\e(1) (\|u_1\|_{H_{\Sig_\e}} +
  \|u_3\|_{H_{\Sig_\e}}) \left( \frac{1}{\e^k} \sum_{l=0}^{\ov{C} \e^{-k}}
\e^2 |\mu_l| \b_l^2 \right)^{\frac 12}.
\end{equation}
From \eqref{eq:Tu2u2u2} and \eqref{eq:asasas} and the fact that
$T_{\Sigma_\e}$ is self-adjoint we deduce that
\begin{eqnarray*}
  \frac{C^{-1}}{\e^k} \sum_{l=0}^{\ov{C} \e^{-k}} \e^2 \mu_l
  \b_l^2 + o_\e(1) \left( \sum_{l=0}^{\ov{C} \e^{-k}} \e^2 |\mu_l| \b_l^2
  \right)^{\frac 12} (\|u_1\|_{H_{\Sig_\e}} + \|u_3\|_{H_{\Sig_\e}}) & \leq &
(T_{\Sig_\e} u, u_2)_{H_{\Sig_\e}} = \l(u, u_2)_{H_{\Sig_\e}} \\
  & \leq & C \e^{\varsigma} \|u\|_{H_{\Sig_\e}} \|u_2\|_{H_{\Sig_\e}}.
\end{eqnarray*}
Also, from Lemma \ref{l:decae}, testing the eigenvalue equation on
$\sum_{l \leq l_0} \b_l \Psi_l$, where $l_0$ is the biggest integer
such that $\mu_{l_0} < 0$, one finds
$$
  \frac{1}{\e^k} \e^2 \sum_{l \leq l_0} \b_l^2 |\mu_l| = O(\e^3)
\|u\|_{H_{\Sig_\e}}.
$$
The last two formulas imply that $\frac{1}{\e^k} \sum_{l=0}^{\ov{C}
\e^{-k}} \b_l^2 = o_\e(1)$, namely that $\|u_2\|_{H_{\Sig_\e}}$
tends to zero as $\e$ tends to zero. This concludes the proof.
\end{pf}

\subsection{Proof of Theorem \ref{t:m}}\label{ss:pf}

Once Propositions \ref{p:morse} and \ref{p:ei} have been
established, the proof goes as in \cite{malm}, Section 8 (see also
\cite{mal2} Section 5) and therefore we will limit ourselves to
sketch the main steps.

\

\noindent  First of all, using Kato's theorem, see \cite{ka}, pag.
445, one can prove that the eigenvalues of $T_{\Sig_\e}$ are
differentiable with respect to $\e$, and if $\l$ is such an
eigenvalue, then there holds
\begin{equation}\label{eq:dere}
  \frac{\partial \l}{\partial \e} = \left\{ \hbox{eigenvalues of }
  Q_\l \right\},
\end{equation}
where $Q_\l : H_\l \times H_\l \to \R$ is the quadratic form given
by
\begin{eqnarray}\label{eq:ql}
  Q_\l(u,v) & = & (1-\l) \frac{2}{\e} \int_{\Sig_\e} \n u
  \cdot \n v - p (p-1) \int_{\Sig_\e} u v
  u_{I,\e}^{p-2} \left(\frac{\partial \ov{u}_{I,\e}}{\partial
  \e}\right) \left( \e \, \cdot \right).
\end{eqnarray}
Here $H_\l \subseteq H_{\Sig_\e}$ stands for the eigenspace of
$T_{\Sig_\e}$ corresponding to $\l$ and the function $\ov{u}_{I,\e}
: \O \to \R$ is defined by the scaling $\ov{u}_{I,\e}(x) =
u_{I,\e}(\e x)$, where $u_{I,\e}$ is as in Section \ref{s:as}.
Notice that, since $\l$ might have multiplicity bigger than $1$,
when we vary $\e$ this eigenvalue can split into a multiplet, which
is allowed by formula \eqref{eq:dere}.

Taking $\l$ as in Proposition \ref{p:ei}, we can apply
\eqref{eq:dere}, and evaluate the quadratic form in \eqref{eq:ql} on
the couples of eigenfunctions in $H_\l$, which are characterized by
\eqref{eq:eigenf}. Reasoning as in \cite{mal2}, Proposition 5.1 one
can prove the following result.

\begin{pro}\label{p:dere2}
Let $\l$ be as in Proposition \ref{p:ei}. Then for $\e$ small one
has
$$
\frac{\partial \l}{\partial \e} = \frac{1}{\e} (\ov{F} + o_\e(1)),
$$
where $\ov{F}$ is a positive constant depending on $N, k$ and $p$.
\end{pro}

\noindent Now we are in position to prove the following proposition,
which states the invertibility of $T_{\Sig_\e}$ for suitable values
of $\e$.

\begin{pro}\label{p:gap}
For a suitable sequence $\e_j \to 0$, the operator $J''_\e(u_{I,\e})
: H^1(\O_\e) \to H^1(\O_\e)$ is invertible and the inverse operator
satisfies $\left\| J''_{\e_j}(u_{I,\e\-j})^{-1}
\right\|_{H^1(\O_{\e_j})} \leq \frac{C}{\min \{\e_j^{k},
\e_j^{\varsigma}\}}$, for all $j \in \N$.
\end{pro}

\begin{pf}
{From} Proposition \ref{p:morse} we have that, letting $N_\e$ denote
the Morse index of $T_{\Sig_\e}$, there holds $N_\e \simeq \left(
\frac{\ov{\a}}{C_k}\right)^{\frac k2} Vol(K) \e^{-k}$. For $l \in
\N$, let $\e_l = 2^{-l}$. Then we have
\begin{equation}\label{eq:diffe}
  N_{\e_{l+1}} - N_{\e_l} \simeq \left(
\frac{\ov{\a}}{C_k}\right)^{\frac k2} Vol(K) (2^{k(l+1)} - 2^{kl})
\simeq \left( \frac{\ov{\a}}{C_k}\right)^{\frac k2} Vol(K) (2^k-1)
\e_l^{-k}.
\end{equation}
By Proposition \ref{p:dere2}, the eigenvalues $\l$ of $T_{\Sig_\e}$
with $|\l| \leq \e^{\varsigma}$ are strictly monotone functions of
$\e$ so by the last equation the number of eigenvalues which cross
$0$, when $\e$ decreases from $\e_l$ to $\e_{l+1}$, is of order
$\e_l^{-k}$. Now we define
$$
A_l = \left\{ \e \in (\e_{l+1}, \e_l) \; : \; \ker T_{\Sig_\e} \neq
\emptyset \right\}; \qquad \qquad B_l = (\e_{l+1}, \e_l) \setminus
A_l.
$$
By Proposition \ref{p:dere2} and \eqref{eq:diffe} we deduce that
card$(A_l) < C \e_l^{-k}$, and hence there exists an interval $(a_l,
b_l)$ such that
\begin{equation}\label{eq:albl}
  (a_l, b_l) \subseteq B_l; \qquad \qquad |b_l - a_l| \geq C^{-1}
  \frac{\hbox{meas}(B_l)}{\hbox{card}(A_l)} \geq C^{-1} \e_l^{k+1}.
\end{equation}
{From} Proposition \ref{p:dere2}, then it follows that every
eigenvalue of $T_{\Sigma_{\frac{a_l+b_l}{2}}}$ in absolute value is
bigger than $C^{-1} \min \{ \e^k, \e^{\varsigma} \}$ for some $C >
0$. By Lemma \ref{l:expcomp} then the same is true for the
eigenvalues of $J''_\e(u_{I,\e})$ so the conclusion follows taking
$\e_j = \frac{a_j+b_j}{2}$.
\end{pf}

\begin{rem}\label{r:final}
The arguments in the proof of Proposition \ref{p:ei} can be easily
adapted to the case in which $|\l| \leq C^{-1} \e^2$ with $C$ is
sufficiently large. Therefore the result of Proposition \ref{p:gap}
can be improved to $\left\| J''_{\e_j}(u_{I,\e_j})^{-1}
\right\|_{H^1(\O_{\e_j})} \leq \frac{C}{\min \{\e_j^{k},
\e_j^{2}\}}$, for all $j \in \N$.
\end{rem}

\

\noindent Below, $\| \cdot \|$ denotes the standard norm of
$H^1(\O_\e)$. For the values of $\e$ such that $J''_\e(u_{I,\e})$ is
invertible, it is sufficient to apply the contraction mapping
theorem. Writing $\e = \e_j$, we find a solution $\tilde{u}_\e$ of
\eqref{eq:tpe} in the form $\tilde{u}_\e = u_{I,\e} + w$, with $w
\in H^1(\O_\e)$ small in norm. Since $J''_{\e}(u_{I,\e})$ is
invertible we have that $J'_{\e}(u) = 0$ if and only if $w = -
  \left( J''_{\e}(u_{I,\e}) \right)^{-1} \left[
  J'_{\e}(u_{I,\e}) + G(w) \right]$,
where
$$
G(w) = J'_\e(u_{I,\e} + w) - J'_\e(u_{I,\e}) - J''_\e(u_{I,\e})[w].
$$
Note that
$$
G(w)[v] = - \int_{\O_\e} \left[ (u_{I,\e} + w)^p - u_{I,\e}^p - p
u_{I,\e}^{p-1} w \right] v; \qquad v \in H^1(\O_\e).
$$
Reasoning as in the last section of \cite{malm2}, we find the
following estimates, which are based on elementary inequalities
\begin{equation}\label{eq:Gw}
  \|G(w)\| \leq
  \begin{cases}
    C \|w\|^p & \text{ for } p \leq 2,  \\
    C \|w\|^2 & \text{ for } p > 2;
  \end{cases} \qquad \qquad \|w\| \leq 1;
\end{equation}
\begin{equation}\label{eq:Gw1w2}
  \| G(w_1) - G(w_2) \| \leq
  \begin{cases}
    C \left( \|w_1\|^{p-1} + \|w_2\|^{p-1} \right) \|w_1 - w_2\|
    & p \leq 2,  \\
    C \left( \|w_1\| + \|w_2\| \right) \|w_1 - w_2\| & p > 2;
  \end{cases} \qquad \qquad \|w_1\|, \|w_2\| \leq 1.
\end{equation}
Defining $F_\e : H^1(\O_\e) \to H^1(\O_\e)$ as
$$
  F_\e(w) = - \left( J''_{\e}(u_{I,\e}) \right)^{-1}
  \left[ J'_\e(u_{I,\e}) + G(w) \right], \qquad w \in H^1(\O_\e),
$$
we will show that $F_\e$ is a contraction in some closed ball of
$H^1(\O_\e)$. From \eqref{eq:uke1}, Proposition \ref{p:gap} (with
Remark \ref{r:final}) and \eqref{eq:Gw}-\eqref{eq:Gw1w2} we get
\begin{equation}\label{eq:fw}
  \|F_{\e}(w)\| \leq
  \begin{cases}
    C \e^{-(k+1)} \left( \e^{I + 1 - \frac{k}{2}} + \|w\|^p \right) &
    \text{ for } p \leq 2,  \\ C \e^{-(k+1)} \left( \e^{I + 1 - \frac{k}{2}}
    + \|w\|^2 \right) & \text{ for } p > 2;
  \end{cases} \qquad \qquad \|w\| \leq 1;
\end{equation}
\begin{equation}\label{eq:fw1w2}
  \| F_{\e}(w_1) - F_{\e}(w_2) \| \leq
  \begin{cases}
    C \e^{-(k+1)} \left( \|w_1\|^{p-1} + \|w_2\|^{p-1} \right)
    \|w_1 - w_2\| & p \leq 2,  \\
    C \e^{-(k+1)} \left( \|w_1\| + \|w_2\| \right) \|w_1 - w_2\|
    & p > 2; \end{cases} \qquad  \|w_1\|, \|w_2\| \leq 1.
\end{equation}
Now we choose integers $d$ and $k$ such that
\begin{equation}\label{eq:dk}
  d >
  \begin{cases}
    \frac{k+1}{p-1} & \text{ for } p \leq 2,  \\
    k+1 & \text{ for } p > 2;
  \end{cases} \qquad \qquad I  > d - 1 + \frac{3}{2}k,
\end{equation}
and we set
$$
\mathcal{B} = \left\{ w \in H^1(\O_\e) \; : \; \|w\| \leq \e^d
\right\}.
$$
{From} \eqref{eq:fw}-\eqref{eq:fw1w2} we deduce that $F_{\e}$ is a
contraction in $\mathcal{B}$ for $\e$ small, so the existence of a
critical point $\tilde{u}_\e$ of $J_\e$ near $u_{I,\e}$ follows. All
the properties listed in Theorem \ref{t:m}, including the positivity
of the solutions, follow from the construction of $u_{I,\e}$ and
standard arguments. As in \cite{malm2}, when $p$ is supercritical
one can use truncations and $L^\infty$ estimates to apply the above
argument working in the function space $H^1(\O_\e) \cap
L^\infty(\O_\e)$.

\begin{rem}\label{r:pac}
With the arguments given in Section \ref{s:real} we could obtain
sharp estimates on the Morse index of $T_{\Sig_\e}$ and on the
eigenfunctions corresponding to resonant eigenvalues. In particular
about the latter we showed that the components in $H_2, H_3$ are
small, and that in $H_1$ the Fourier modes are localized near some
precise frequencies. This allowed us to prove Proposition
\ref{p:gap} using Kato's theorem.

Even if we did not work the computations out, it seems it should be
possible to give a more rough characterization of these
eigenfunctions (in particular on the $H_2$ component) and to prove a
(non sharp) estimate on the derivatives of the eigenvalues, still
obtaining invertibility. This might slightly simplify the proof of
existence, although most of the delicate estimates will be shifted
from the analysis of $T_{\Sig_\e}$ to that of the quadratic form
$Q_\l$ defined in \eqref{eq:ql}.
\end{rem}

\end{document}